%% file: rmsf-master.tex
\documentclass[11pt]{article}
\usepackage{anysize}
\marginsize{3cm}{3cm}{3cm}{3cm}
\usepackage[utf8]{inputenc}
\usepackage{amsmath}
\usepackage{amsfonts}
\usepackage{amssymb}
\usepackage{diagrams}
\usepackage{theorem}
\usepackage{stmaryrd}
\usepackage{eucal}
\usepackage{graphicx}
\usepackage{indentfirst}
\usepackage{hyperref}
\usepackage{baskervald}
\usepackage{mathrsfs}
\usepackage{chngpage,calc}

\usepackage[T1]{fontenc}

\normalfont

\makeatletter
\begingroup
\gdef\th@dotted{\normalfont\itshape
  \def\@begintheorem##1##2{%
        \item[\hskip\labelsep \theorem@headerfont ##1\ ##2.]}%
\def\@opargbegintheorem##1##2##3{%
   \item[\hskip\labelsep \theorem@headerfont ##1\ ##2\ (##3).]}}
\endgroup
\makeatother

\theoremstyle{dotted}

\newtheorem{thm}{Theorem}[section]
\newtheorem{lemma}[thm]{Lemma}
\newtheorem{prop}[thm]{Proposition}
\newtheorem{corr}[thm]{Corollary}

\newtheorem{introthm}{Theorem}%
\newtheorem{introprop}[introthm]{Proposition}%

\makeatletter
\begingroup
\gdef\th@upshape{\normalfont
  \def\@begintheorem##1##2{%
        \item[\hskip\labelsep \theorem@headerfont ##1\ ##2.]}%
\def\@opargbegintheorem##1##2##3{%
   \item[\hskip\labelsep \theorem@headerfont ##1\ ##2\ (##3).]}}
\endgroup
\makeatother

\theoremstyle{upshape}

\newtheorem{opr}[thm]{Definition}
\newtheorem{rem}[thm]{Remark}
\newtheorem{ntn}[thm]{Notation}
\newtheorem{exm}[thm]{Example}
\newtheorem{const}[thm]{Construction}
\newtheorem{conv}[thm]{Convention}

\newcommand{\proof}[1][Proof. ]{\smallskip\noindent{\bf #1}}
\def\endproof{\hfill\ensuremath{\square}\par\medskip}

\makeatletter
\@addtoreset{equation}{section}
\makeatother

\newcommand{\bR}{{\mathbb R}}

\newcommand{\bN}{{\mathbb N}}

\newcommand{\bL}{{\mathbb L}}

\newcommand{\cB}{{\mathcal B}}

\newcommand{\cC}{{\mathcal C}}
\newcommand{\cD}{{\mathcal D}}
\newcommand{\cJ}{{\mathcal J}}
\newcommand{\cF}{{\mathcal F}}
\newcommand{\cG}{{\mathcal G}}
\newcommand{\cN}{{\mathcal N}}
\newcommand{\cM}{{\mathcal M}}
\newcommand{\cL}{{\mathcal L}}
\newcommand{\cZ}{{\mathcal Z}}

\newcommand{\cS}{{\mathcal S}}

\newcommand{\cW}{{\mathcal W}}

\newcommand{\cX}{{\mathcal X}}
\newcommand{\cY}{{\mathcal Y}}

\newcommand{\cE}{{\mathcal E}}

\newcommand{\cR}{{\mathcal R}}

\renewcommand{\phi}{\varphi}

\newcommand{\Hom}{\operatorname{Hom}}

\newarrow {Into} {boldhook}--->
\newarrow {Mapsto} |--->
\newcommand{\Set}{{\rm \bf Set}}

\newcommand{\SSet}{{\rm \bf SSet}}
\newcommand{\Cat}{{\rm \bf Cat }}

\newcommand{\Mod}{{\rm \bf Mod}}

\newcommand{\CRing}{{\rm \bf CRing}}

\newcommand{\Ob}{{\rm Ob} \,}
\newcommand{\Mor}{{\rm Mor} \,}
\newcommand{\Ho}{{\rm Ho} \,}

\newcommand{\op}{{\sf op}}

\newcommand{\bc}{\mathbf c}

\renewcommand{\lim}{\operatorname{\varprojlim}}
\newcommand{\colim}{\operatorname{\varinjlim}}

\newcommand{\Map}{\operatorname{Map}}
\newcommand{\Sect}{{\sf Sect}}

\newcommand{\Lax}{{\sf Lax}}
\newcommand{\Cart}{{\sf Cart}}

\newcommand{\Lat}{\mathscr{Lat}}
\newcommand{\Mat}{\mathscr{Mat}}
\newcommand{\sR}{\mathscr{R}}
\newcommand{\sL}{\mathscr{L}}
\newcommand{\sA}{\mathscr{A}}
\newcommand{\sS}{\mathscr{S}}

\newcommand{\Alpha}{\mathrm{A}}

\newcommand{\Fun}{\operatorname{\mathsf{Fun}}}

\newcommand{\cf}{cf}

\author{Edouard Balzin}
\title{Reedy Model Structures in Families}
\date{}

\begin{document}

	\maketitle
    \thispagestyle{empty}
  \begin{abstract}
  Given a family of model categories $\cal E \to \cal R$ over a Reedy category,
  we outline a set of conditions which lead to the existence of a Reedy model
  structure on the category of sections ${\sf Sect}(\cal R, \cal E)$. We prove
  that for a wide class of examples, this model structure serves as a strictification of the
  $(\infty,1)$-category of sections of the higher-categorical family associated to
  $\cal E \to \cal R$.
  \end{abstract}

\tableofcontents
\include{rmsf-intro}

\include{rmsf-fibrations}
\include{rmsf-reedy}

\include{rmsf-comparison}
\include{rmsf-appendices}

\include{rmsf-biblio}
\end{document}

%% file: rmsf-intro.tex
\section*{Introduction}
\addcontentsline{toc}{section}{\protect\numberline{}Introduction}%
Let $\cM$ be a model category. Given a small category $\cC$, the associated
functor category $\Fun(\cC,\cM)$ can be endowed with objectwise weak equivalences:
take those natural transformations $ X \to Y$ such that for each $c \in \cC$,
the value $  X(c) \to Y(c)$ is a weak equivalence in $\cM$.
For a general \cite{QHA} choice of $\cM$ (assuming small (co)limits, see Definition \ref{defmodelcat}) and $\cC$,
there is no known way to complete the class of objectwise weak equivalences on $\Fun(\cC,\cM)$
into a model structure.
One option to treat this issue is to assume
additional requirements about $\cM$, such as being cofibrantly generated, or further, combinatorial.
Another way concerns the assumptions on $\cC$.

The notion of a Reedy category is attributed to \cite{REEDY}. To fix a definition, a
\emph{Reedy category}  is a small category $\cR$
possessing two subcategories $\cR_- \subset \cR \supset \cR_+$ with
$\Ob \cR_- = \Ob \cR = \Ob \cR_+$,
and a degree function
$deg: \cR \to \bN$ taking values in natural numbers, such that
	\begin{enumerate}
    \item the isomorphisms of $\cR$ are identities,
		\item the non-identities of $\cR_-$ lower the value of $deg$, and the non-identities of $\cR_+$ raise the value of $deg$,
    \item any morphism $f: x \to y$ factors uniquely as $f = f^+ \circ f_-$, where
    $f_- \in \cR_-$ and $f^+ \in \cR_+$.
  \end{enumerate}
Many elementary diagram categories are naturally Reedy: for example, the two-arrow category
$a \longleftarrow b \longrightarrow c$ appearing in the calculus of pullbacks,
can be endowed with at least two Reedy structures; in the first one, both maps raise the degree,
in the other, both maps lower the degree. Similarly,
the category
$$
[n] = 0 \longrightarrow ... \longrightarrow n
$$
admits two different Reedy structures: for the first structure,
$deg (0) = 0$ and all maps advance the degree, for the second, $deg (0) = n$ and all
maps lower the degree. Taking all $[n]$ for $n \in \bN$ defines a
full subcategory $\Delta \subset \Cat$ of the category of small categories. The surjections,
the injections and $deg ([n]) = n$ make $\Delta$ (and $\Delta^\op$) into a Reedy category.
Furthermore, for any category $\cC$, denote by $\Delta \cC$ the category of functors $\bc:[n] \to \cC$,
with morphisms $\bc \to \bc'$ given by functors $[n] \to [m]$ in the overcategory
$\Cat /\cC$. The Reedy structure on $\Delta$ induces a Reedy structure
on $\Delta \cC$ \cite[22.10]{DHKS}. This example shows that Reedy categories exist
in a great variety.

The classical result (see e.g. \cite{HOVEY}) asserts that given a Reedy category
$\cR$ and a model category $\cM$, there exists an explicitly-constructed model structure
on $\Fun(\cR,\cM)$, with objectwise weak equivalences. The existence of a model structure
on functors from a Reedy category permits the computation of various homotopy limits,
even for a general category $\cC$: the assignment $$(\bc:[n] \to \cC) \mapsto \bc(n)$$
defines a (homotopy cofinal) functor \cite[22.11]{DHKS} $p_t: \Delta \cC \to \cC$. We can use
the pull-back along $p_t$ to embed $\Fun(\cC,\cM)$ in $\Fun(\Delta \cC,\cM)$, and then use
the Reedy model structure on $\Fun(\Delta \cC,\cM)$ for homotopy colimit computations.

\subsubsection*{Sections over a Reedy category}
In this paper, we study families of model categories and their associated categories
of sections. In light of the Grothendieck construction \cite{SGA1,VIST}, a natural way to treat
a family of categories indexed by $\cC$ is to consider a functor $\cE \to \cC$ together
with additional conditions (that we survey in Section \ref{chaptergrothfib}).

Given a family $\cE \to \cC$, we can associate to it the category
of sections $\Sect(\cC,\cE)$, consisting of those functors
$S: \cC \to \cE$ such that $p \circ S = id_\cC$. The case of functors $\cC \to \cM$
can be recovered by considering the trivial family, given by projection
$pr_\cC: \cM \times \cC \to \cC$: the sections of $pr_\cC$ are identified with $\Fun(\cC,\cM)$.

For an example of a nontrivial model-categorical family, denote by $\CRing$ the category of commutative rings and consider
a small diagram $A: \cC \to \CRing$. For a ring $R$, denote by
$C^\bullet (\Mod_R)$ the
category of chain complexes of modules over $R$, that can be endowed with projective
model structure \cite{HOVEY}.
Define $\Mod_A$ to be the category consisting of pairs $(c,M)$ where $c$ is in $\cC$ and
$M$ is a complex of $A(c)$-modules. A morphism $(c,M) \to (c',N)$ in $\cE$ is given by
a map $f:c \to c'$ in $\cC$ and a chain map $M \to N|_{A(c)}$ of complexes of $A(c)$-modules, where
we restrict $N$ along the map $A(f):A(c) \to A(c')$; by adjunction, it is the same data as a map of
complexes of $A(c')$-modules $A(c') \otimes_{A(c)} M \to N$.
This example of algebras and modules can be generalised \cite{HAG2,HAPRASMA},
with the abstraction being the notion of a Quillen presheaf \cite{H-S,BARLEFT}, which
is the same data as
a functor from a small category $\cC$ to model categories and Quillen adjunctions.

When $\cC = \Delta$,
the category of sections $\Sect(\Delta,\Mod_A)$ plays an important role in
the questions of cohomological descent: inside $\Sect(\Delta,\Mod_A)$, there is
a full subcategory denoted $\Sect_{\bL cart}(\Delta,\Mod_A)$ consisting of sections $S$ that send
$[n] \to [m]$ in $\Delta$ to a quasi-isomorphism $A(m) \otimes^\bL_{A(n)} S(n) \to S(m)$. As
explained in \cite{HAG2}, the category $\Sect_{\bL cart}(\Delta,\Mod_A)$ is the category
of descent data in derived algebraic geometry, see also Remark \ref{remdescent}.

To work with sections homotopically, one would like to, just as in the case of functors
to a model category, have a model structure. We thus pose the following problem: given a family
$p:\cE \to \cR$ over a Reedy category $\cR$,
formulate the requirements on $p$ for
the category $\Sect(\cR,\cE)$ to have a model structure, that specialises
to the Reedy model structure in the case of a trivial family.

Our findings are summarised by Theorem \ref{reedymodelstructuretheorem} in the main text:
\begin{introthm}
\label{reedymodelstructuretheoremintro}
Let $\cR$ be a Reedy category and $\cE \to \cR$ an \emph{admissible model semifibration}.
Then the category of sections $\Sect(\cR, \cE)$ has a model structure with objectwise
weak equivalences.
\end{introthm}
Let us explain the conditions imposed on $\cE \to \cR$.

The functor $p:\cE \to \cR$ cannot be completely general: as it turns out, the appropriate condition
on $p$ is that of a \emph{semifibration}, which is a
mixture of the conditions of Grothendieck fibration and opfibration.
In detail, the semifibration property consists in requiring
that the restriction $\cE|_{\cR_-} \to \cR_-$ is a Grothendieck prefibration; the restriction
$\cE|_{\cR_+} \to \cR_+$ is dually assumed to be a preopfibration. Denoting $\cE(c):= p^{-1}(c)$, these two conditions
imply that for each $f_-: x \to y$ in $\cR_-$, there is a functor $f_-^*: \cE(y) \to \cE(x)$,
and dually for each $f^+: x \to y$ in $\cR_+$, there is a functor
$f^+_!: \cE(x) \to \cE(y)$. The final requirement of a semifibration relates these
two classes of functors, see Definition \ref{semifibrationdefinition} and
Proposition \ref{janusmate} for detail.

Given an object $x \in \cR$ of a Reedy category $\cR$, denote $Lat(x) \subset \cR /x$ the
``latching category'', consisting of all $\cR_+$-maps $y \to x$ without the identity
map. Similarly, the ``matching category'' $Mat(x) \subset x \backslash \cC$ consists
of all $\cR_-$-maps $x \to y$ minus the identity.
The semifibration property of $\cE \to \cR$ implies the existence of two restriction functors
$L_x: \cE|_{Lat(x)} \to \cE(x)$ and $R_x: \cE|_{Mat(x)} \to \cE(x)$. Given a section
$S: \cR \to \cE$, the latching and the matching object functors are defined as follows:
$$
\Lat_x S: =\colim_{Lat(x)} L_x S|_{Lat(x)}, \, \, \, \Mat_x S := \lim_{Mat(x)} R_x S|_{Mat(x)}
$$
(we assume that all needed (co)limits exist).
Given a section
$S: \cR \to \cE$ of the semifibration $\cE \to \cR$, its value $S(x)$ fits into a diagram
$$
\Lat_x S \to S(x) \to \Mat_x S
$$
and just as in the classical Reedy case, such diagrams for different $x$ completely control the behaviour of $S$.

The semifibration $\cE \to \cR$ being \emph{model} means that each fibre $\cE(x)$ is a model category,
and that we also require the transition functors $f^+_!:\cE(x) \to \cE(y)$
along $\cR_+$-maps to preserve cofibrations and
trivial cofibrations (we do not require the preservation of colimits, which permits
to consider such functors as tensor products, and consequently the examples of
algebras and TQFT \cite{LU}).
A dual condition of preserving fibrations and trivial fibrations should hold for $f^*_-: \cE(y) \to \cE(x)$.
For each map
of sections $S \to T$, there is a naturally induced diagram
\begin{diagram}[small]
\Lat_x S & \rTo & S(x) & \rTo & \Mat_x S \\
\dTo      &     &   \dTo &    &   \dTo    \\
\Lat_x T & \rTo & T(x) & \rTo & \Mat_x T. \\
\end{diagram}
We define $S \to T$ to be a \emph{Reedy cofibration}
if for each $x \in \cR$, the map $\Lat_x T \coprod_{\Lat_x S} S(x) \to T(x)$ is a cofibration in $\cE(x)$.
Dually, $S \to T$ is a \emph{Reedy fibration} if for each $x \in \cR$, the
map $S(x) \to \Mat_x S \prod_{\Mat_x T} T(x)$ is a fibration in $\cE(x)$. The Reedy weak equivalences are defined objectwise.

The defined classes of maps give a model structure on $\Sect(\cR,\cE)$, provided
that the model semifibration $\cE \to \cR$ satisfies the \emph{admissibility} condition, Definition \ref{admissibilitydef}.
The admissibility guarantees that the functor $S \mapsto \Lat_x S$ sends (trivial) Reedy cofibrations to
(trivial) cofibrations, and dually for $\Mat_x$. The examples of admissible families are covered by Lemma
\ref{admissibilitycriterion}.
A Quillen presheaf is admissible,
since its transition functors preserve (co)limits.
A Grothendieck fibration in model categories and right derivable functors over $[n]$ is also admissible,
which follows from the simple structure of the matching categories of $[n]$.

The proof of Theorem \ref{reedymodelstructuretheoremintro} is similar to the one of the classical case.
The explicit
description of the model structure permits to verify various properties,
such as being cofibrantly generated under reasonable assumptions, Propositions
\ref{propreedycofibrantgenerationimplicit} and \ref{propreedycofibrantgenerationexplicit}.

\subsubsection*{Comparison with higher-categorical sections}

The second part of our work establishes a relation with the higher category theory
(that we model using the language of quasicategories \cite{LUHTT,CISBOOK}).
To begin, recall that given a category $\cM$ with weak equivalences $\cW$,
its $(\infty,1)$-localisation is an $(\infty,1)$-functor $F: \cM \to L_\cW \cM$,
such that for any $(\infty,1)$-category $\cZ$,
the induced $(\infty,1)$-functor $F^*:\Fun(L_\cW \cM,\cZ) \to \Fun(\cM,\cZ)$ is full and faithful,
and its essential image consists of $(\infty,1)$-functors that send $\cW$ to equivalences
of $\cZ$. For model categories, the $(\infty,1)$-localisation enjoys various
special properties \cite{HINICH}. In particular, Cisinski \cite{CISBOOK} has shown that
given any model category $\cM$ and a small category $\cC$, the natural
$(\infty,1)$-functor $\Fun(\cC,\cM) \to \Fun(\cC, L_\cW \cM)$ induces an
$(\infty,1)$-equivalence
$$
L \Fun(\cC,\cM) \to \Fun(\cC, L_\cW \cM)
$$
where on the left, we localise with respect to the objectwise weak equivalences.

Thanks to the various literature on quasicategories \cite{LUHTT,CISBOOK,GEPNERNIK,HINICH,MGEE1} there is now a well-developed theory of Grothendieck fibrations
of $(\infty,1)$-categories. In some cases \cite{HINICH,MGEE1}, it is known how to $(\infty,1)$-localise
a family $\cE \to \cR$ to get a properly behaved higher-categorical family $L \cE \to \cR$.
We thus ask how the localisation $L \Sect(\cR,\cE)$ of the model structure of
Theorem \ref{reedymodelstructuretheoremintro} is compared with the
$(\infty,1)$-sections $\Sect(\cR,L \cE)$.

We treat the comparison issue in the following generality.
A functor $\cE \to \cR$ is a \emph{left model Reedy fibration} if it is a
Grothendieck opfibration, a Grothendieck fibration over $\cR_-$, and is an admissible
model semifibration. Localising $\cE$ along the union $\cup_c \cW(c)$ of the
fibrewise weak equivalences yields an $(\infty,1)$-functor $L \cE \to \cR$ that
is a coCartesian fibration (in the sense of \cite[Definition 2.4.2.1]{LUHTT}) and a Cartesian fibration
over $\cR_-$. The natural functor $\cE \to L \cE$ induces the $(\infty,1)$-functor
$\Sect(\cR,\cE) \to \Sect(\cR,L \cE)$. The following result is Theorem \ref{thmcomparison}:
\begin{introthm}
\label{thmcomparisonintro}
Let $\cE \to \cR$ be a left model Reedy fibration. Then
the induced infinity-functor $L \Sect(\cR,\cE) \to \Sect(\cR,L \cE)$
is an equivalence of quasicategories.
\end{introthm}
The origins of Theorem \ref{thmcomparisonintro} lie in the paper \cite{H-S}. However,
the proof of \cite[Th\'eor\`eme 18.2]{H-S} is incomplete: the authors do not
provide proofs for the needed higher-categorical arguments, and the model-categorical
considerations of \cite[Th\'eor\`eme 18.2]{H-S} contain a mistake (it is assumed
that the latching object functor preserves fibrant objects). Nonetheless, our proof of Theorem \ref{thmcomparisonintro}
shows that the infinity-categorical
Reedy induction can be carried out
along the general lines of \cite{H-S}.
There are other comparison results \cite{HARPAZCOMP,SPITZCOMP} that work under more assumptions on the family $\cE \to \cR$.
For example, \cite{HARPAZCOMP} bypasses the Reedy induction, yet is only valid
for Quillen presheaves in the combinatorial model setting.

The example of Quillen presheaves has a particularly good strictification result, conjectured
in \cite{H-S}
and outlined in Proposition \ref{propquillenpresheafstrict}:
\begin{introprop}
Let $\cM \to \cC$ be a Quillen presheaf over a small category $\cC$. Then localising $\cM$ along the
fibrewise weak equivalences yields an equivalence of $(\infty,1)$-categories
$$
L\Sect(\cC,\cM) \stackrel \sim \longrightarrow \Sect(\cC, L \cM).
$$
where we treat $\Sect(\cC,\cM)$ as a category equipped with objectwise weak equivalences, with no model
structure.
\end{introprop}

\subsubsection*{Organisation of the paper} Section \ref{chaptergrothfib} covers
various categorical preliminaries that are related to the semifibrations and
Reedy categories. We outline the induction for sections using the notions
of Noether categories, Definition \ref{noethercatdefinition}. The notion
of a semifibration is introduced in Definition \ref{semifibrationdefinition}. Section \ref{chaptergrothfib}
includes more material than is required for Section \ref{chapterreedymodstr}, however,
we decided to keep many propositions for future reference. This comment
also applies to the Appendix, which gives a variation of the argument leading to
Theorem \ref{reedymodelstructuretheoremintro}
in a specialised setting.

Section \ref{chapterreedymodstr} establishes the notion of a model semifibration
and admissibility, Definition \ref{admissibilitydef}, and proves Theorem \ref{reedymodelstructuretheoremintro}.
We discuss the cofibrant generation of Reedy model structures, proving Propositions
\ref{propreedycofibrantgenerationimplicit} and \ref{propreedycofibrantgenerationexplicit}.

In Section \ref{subsectionreedycomp} we address the comparison theorem.
We start by recalling
a result of Lurie \cite[Proposition A.2.9.14]{LUHTT} asserting the higher-categorical
Reedy induction for the functors to a quasicategory. We then generalise it to
the case of a coCartesian fibration $\cE \to \cR$ of quasicategories over a Reedy category that is
also Cartesian over $\cR_-$ and is suitably bicomplete, Proposition \ref{propinfinityreedyinduction}.
After studying some aspects of relative categories in families, we prove a model-categorical
counterpart, Proposition \ref{propinductionformodelcategoriesofsections}, of Proposition \ref{propinfinityreedyinduction}.
Both these propositions lead to the proof of Theorem \ref{thmcomparisonintro}.
We conclude Section \ref{subsectionreedycomp} by discussing the strictification of
Quillen presheaves.

\subsubsection*{Acknowledgements} I would like to thank Carlos Simpson and
Andr\'e Hirschowitz for their pioneering work on the subject. This work
has profited from numerous discussions with Andrea Gagna, Dmitry Kaledin,
Carlos Simpson, and in particular with Yonatan Harpaz who has been of immense
help with his valuable insights and noticing a mistake in an earlier version of this paper.
The author would also like to thank Aaron Mazel-Gee and
Jay Shah for useful correspondences.

This work has been completed during the author's employment as a Hadamard Lecturer
in \'Ecole Polytechnique. I thank both FMJH and CMLS for financial support and
hospitable working atmosphere. Lastly, I thank Paul Taylor for his diagrams package.

%% file: rmsf-fibrations.tex
\section{Semifibrations}
\label{chaptergrothfib}

\subsection{Cartesian arrows, prefibrations, sections}

Let $p:\cE \to \cC$ be a functor. For $c \in \cC$, denote by $\cE(c)$ the fibre
category $p^{-1} (c)$ over $c$. It thus consists of all $X \in \cE$ with
$p(X) = c$ and all the maps $X \to X'$ with $p(X \to X') = id_c$.

\begin{opr}
	\label{cartesianmorphismdefinition}

  A morphism $\alpha: X \to Y$ of $\cE$

  \begin{itemize}

    \item is \emph{$p$-cartesian}, or simply cartesian, if for any other map
    $\beta: X' \to Y$ with $p (\beta) = p (\alpha)$ there exists a unique
    morphism $\gamma: X' \to X$ in $\cE(p(X))$ which
    factors $\beta$ as $\alpha \circ \gamma$.

    \item is \emph{$p$-opcartesian}, or simply opcartesian,
    if for any other map $\delta: X \to Y'$ with $p (\delta) = p (\alpha)$
    there exists a unique morphism $\eta: Y \to Y'$ in $\cE(p(Y))$ which factors
    $\delta$ as $\eta \circ \alpha$.
	\end{itemize}

	A $p$-cartesian or $p$-opcartesian morphism $\alpha: X \to Y$ is
  covering the morphism $f: c \to c'$ iff $p(\alpha) = f$.
\end{opr}

In our definition of cartesian and opcartesian morphisms, we are faithful to
the original terminology of \cite{SGA1}. Today, a different definition of
(op)cartesian maps is presented in many sources \cite{VIST,LU},
with the definition of \cite{SGA1} referred to as
``locally cartesian'' morphism.

\begin{opr}
	\label{prefibrationdefinition}
	A functor $p: \cE \to \cC$ is a
	\begin{itemize}
		\item prefibration iff for any $f: x \to y$ of $\cC$ and $Y \in \cE(y)$
    there exists a cartesian morphism $\alpha: X \to Y$ covering $f$,
    that is, $p(\alpha) = f$.
		\item preopfibration iff for any $f: x \to y$ of $\cC$ and $X \in \cE(x)$
    there exists a cartesian morphism $\delta: X \to Z$ covering $f$,
    that is, $p(\alpha) = f$.
	\end{itemize}
\end{opr}

\begin{lemma} If $p : \cE \to \cC$ is a prefibration, then
$p^\op : \cE^\op \to \cC^\op$ is a preopfibration.  \endproof \end{lemma}

\begin{ntn}
	If $p: \cE \to \cC$ is a prefibration, $f: x \to y$ is a morphism and
  $Y \in \cE(y)$, we shall usually denote a chosen cartesian lift by
  $f^* Y \to Y$. The same applies when $p$ is a preopfibration,
  where for $X \in \cE(x)$, we denote by $X \to f_! X$
  the chosen opcartesian lift.
\end{ntn}

\begin{opr}
	\label{prefibsmalldiscrete}
	A prefibration or preopfibration $q:\cE \to \cC$ is small if both $\cC$
  and $\cE$ are small categories. A prefibration or preopfibration $q$
  is $\emph{discrete}$ if for each $c \in \cC$, the category $\cE(c)$
  has no non-identity maps (in other words, it is isomorphic to a set).
\end{opr}

\begin{lemma}
	Let $p: \cE \to \cC$ be a discrete prefibration.
  Then the composition of cartesian morphisms of $\cE$ is cartesian.
  The dual is true for preopfibrations.
\end{lemma}
\proof Clear, due to the lack of fibre maps. \endproof

In general, not any pre(op)fibration has the property described in the previous
lemma. Those which have it, are called Grothendieck (op)fibrations.

\begin{opr}
	\label{fibrationdefinition}
	A prefibration $p:\cE \to \cC$ is, furthermore, a \emph{Grothendieck
  fibration} iff the composition of cartesian maps is cartesian. The
  definition for \emph{Grothendieck opfibrations} is dual.
\end{opr}

Discrete pre(op)fibrations are thus automatically (op)fibrations.
The examples of non-discrete fibrations are, however, abundant.

\begin{rem}
	It is not necessarily the case that the category $\cE$ is ``bigger''
  than $\cC$. For example, the functor $\cC \to \cC \coprod \cD$
  is a fibration and an opfibration.
\end{rem}

\begin{rem}
	In what follows, (op)fibrations will be considered as special cases of
  pre(op)fibra-tions, with additional remarks where necessary.
  Otherwise, any definition or a result given for a pre(op)fibration
  implies the same for an (op)fibration.
\end{rem}

\begin{const}
	\label{grothendieckconstruction}
	Given a functor $E$ from $\cC$ to categories, we produce an opfibration,
  which we denote $\int E \to \cC$ and call the \emph{Grothendieck construction}
  \cite{VIST} of $E$. An object of $\int E$ is a pair $(c,X)$ of $c \in \cC$
  and $X \in E(c)$, and a morphism $(c,X) \to (c',X')$ consists of
  $f: c \to c'$ together with a map $\alpha: E(f)(X) \to X'$ in $E(c')$.

	Dually, for a contravariant category-valued functor $F$ defined on $\cC$,
  its Grothendieck construction is a fibration $\int F \to \cC$ with same pairs
  $(c,Y)$ serving as objects, but with maps given by pairs of $f: c \to c'$ and
  $\beta: Y \to F(f) Y'$ in $F(c)$.
\end{const}

Grothendieck construction motivates the following perspective.
Consider a prefibration $p:\cE \to \cC$. Let $f: c \to c'$ be a morphism in
$\cC$ and $Y \in \cE(c')$. Choose a cartesian morphism $\alpha:f^*Y \to Y $
covering $f$. This specifies an object $f^*Y \in \cE(c)$.
By the universal property of cartesian maps, the assignment $Y \mapsto f^*Y$
defines a functor $f^*:\cE(c') \to \cE(c)$, which is called transition
functor along $f$. Due to the universal property of cartesian arrows,
for each composable pair $f,g$,
there exists a `coherence' natural transformation
$f^* \circ g^* \to (g \circ f)^* $, which is an isomorphism if $p$ is a
Grothendieck fibration. For any composable triple of arrows $f, g, h$,
any choice of coherence morphisms leads to the following commutative diagram:

\begin{equation}
\label{prefibcoherence}
\begin{diagram}[small]
f^* g^* h^* & \rTo & (gf)^* h^* \\
\dTo 	&					&	\dTo\\
f^* (hg)^* & \rTo  & (hgf)^* . \\
\end{diagram}
\end{equation}
For a preopfibration, the whole picture is dual. In the literature
(see \cite{SGA1} and \cite{VIST} for the case of Grothendieck
fibrations), such choice of an assignment $f \mapsto f^*$ together with
coherence isomorphisms is called a cleavage. One may thus wonder if there
is a way to obtain (lax) category-valued functors from (pre)fibrations.

\begin{opr}
	Let $p:\cE \to \cC$ and $q: \cE' \to \cC$ be two functors.
	\begin{itemize}

    \item A morphism of $p$ and $q$ is a functor $F: \cE \to \cE'$ commuting
    with the functors to $\cC$, that is, $q \circ F = p$.

    \item A section of $p$ is a functor $S: \cC \to \cE$ such that
    $p \circ S = id_\cC$. In other words, it is a morphism from
    $id_\cC: \cC \to \cC$ to $p: \cE \to \cC$.
	\end{itemize}
	Given two morphisms $F,F': \cE \to \cE'$, a morphism between them is a natural
  transformation $\alpha:F \to F'$ such that for each $X$ in the domain
  $\cE$, $\alpha_X$ projects to $id_{p(X)}$.
\end{opr}

We denote by $\Lax(\cE,\cE')$ the category of morphisms between $p$ and $q$,
with the functors to $\cC$ being implicit. By $\Sect(\cC,\cE)= \Lax(\cC,\cE)$
we denote the category of sections of $p$.

\begin{opr}
	Let $p:\cE \to \cC$ and $q: \cE' \to \cC$ be two prefibrations (respectively
  preopfibrations). A morphism $F: \cE \to \cE'$ is called a
  \emph{cartesian morphism} if it takes (op)cartesian morphisms of $\cE$ to
  (op)cartesian morphisms of $\cE'$.
\end{opr}

We denote by $\Cart(\cE,\cE')$ the full subcategory of $\Lax(\cE,\cE')$
consisting of cartesian morphisms.

\begin{const}
	Take a fibration $p: \cE \to \cC$, and for each $c \in \cC$, denote by
  $\cC / c$ the category of objects over $c$ \cite{ML}. The forgetful functor
  $\cC / c \to \cC$ is an fibration. Then the assignment
  $c \mapsto \Cart(\cC / c, \cE)$ defines a contravariant category-valued
  functor on $\cC$. When $\cC$ is small, this construction is inverse up to
  an equivalence \cite{VIST} to (Grothendieck)
  Construction \ref{grothendieckconstruction}.

	If $p$ is only a prefibration, the assignment
  $c \mapsto E(c)=\Cart(\cC / c, \cE)$ defines a lax contravariant functor
  from $\cC$ to categories. Indeed, for each $f: c \to c'$, we get a functor
  $f^* : E(c') \to E(c)$, and as before, one can witness the existence of
  natural transformations $f^*g^* \to (gf)^*$ and of the diagram
  like (\ref{prefibcoherence}).
\end{const}

The cited result \cite{VIST} implies that any fibration
(and, similarly, an opfibration)
$p: \cE \to \cC$ can be, up to an equivalence, replaced by an fibration
$\tilde p: \tilde \cE \to \cC$, for which the assignment $c \mapsto \cE(c)$
can be made into a strict functor by a choice of transition functors along
maps in $\cC$. We call the fibrations (similarly, fibrations) with the latter
property \emph{strictly cleavable}.


\begin{opr}
	\label{isofibrationdefinition}
	A functor $p:\cE \to \cC$ is an \emph{isofibration} if for any
  isomorphism $f:c \stackrel \sim \to d$ of $\cC$ and an object $Y$
  with $p(Y) = d$ there exists an isomorphism
  $\alpha: X \stackrel \sim \to Y$ with $p \alpha = f$.
\end{opr}

A Grothendieck (op)fibration is automatically an isofibration,
but a pre(op)fibration is not. In particular, in an arbitrary prefibration,
a cartesian lift of an isomorphism is not necessarily an isomorphism.

\begin{conv}
	\label{convprefibrationareisofibration}
 From now on, any prefibration or preopfibration we consider
  is assumed to be also an isofibration. For an isofibration $p: \cE \to \cC$
  and $c \in \cC$, the notation $\cE(c)$ will denote $p^{-1}(c)$,
  the strict categorical fibre of $p$ over $c$. Note that in this case,
	the strict fibre is equivalent to the \emph{essential fibre} of $p$
  over $c$: the objects of the latter are pairs of $d \in \cD$ and $\alpha: p(d) \cong c$
  in $\cC$, and morphisms $(d,\alpha) \to (d',\beta)$ are given by
  $f:d \to d'$ such that $\beta p(f) = \alpha$.
  In particular, $p(f)$ is an isomorphism.
\end{conv}

\begin{exm}
	\label{laxfunctorgrothendieckconstruction}
	Let $L: \int \cE \to \int \cE'$ be a morphism between two Grothendieck constructions of covariant functors $\cE, \cE': \cC \to \Cat$. For each $c \in \cC$, $L$ specifies a functor $L_c : \cE(c) \to \cE'(c)$. For each morphism $f: c \to c'$, we get a $2$-square
	\begin{diagram}[small]
	\cE(c)  & \rTo^{L_c} & \cE'(c) \\
	\dTo<{\cE(f)}			& \stackrel{L_{f}}{\Leftarrow} & \dTo>{\cE'(f)} \\
	\cE(c')  & \rTo_{L_{c'}} & \cE'(c'). \\
	\end{diagram}
	The natural transformation appears because the image under $L$ of an opcartesian map $X \to \cE(f) X$ ($X \in \cE(c)$) may not be opcartesian. Factoring $LX \to L \cE(f) X$,
	$$
	LX \to \cE'(f) LX \to L \cE(f) X,
	$$
	gives $\cE'(f) LX \to L \cE(f) X$; for each $X \in \cE(c)$, all such maps assemble into $L_f$. For two composable arrows $f: c \to c'$, $g: c' \to c''$, there is a pasting property relating $L_f, L_g$ and $L_{g f}$: the pasting of this diagram
	\begin{diagram}[small]
	\cE(c) & \rTo^{\cE(f)} & \cE(c') & \rTo^{\cE(g)} & \cE(c'') \\
	\dTo<{L_c}  &	\stackrel{L_f}{\Downarrow}	& \dTo_{L_{c'}} &	\stackrel{L_g}{\Downarrow}	&  \dTo>{L_{c''}}  \\
	\cE'(c) & \rTo_{\cE'(f)} & \cE'(c') & \rTo_{\cE'(g)} & \cE'(c'')\\
	\end{diagram}
	equals $L_{gf}$.

	For fibrations, there is a difference on the level of 2-diagrams.
  Consider $\cF, \cF': \cC^\op \to \Cat$ and take a
  lax morphism $M: \int \cF \to \int \cF'$ of fibrations
  over $\cC$. For $f: c \to c'$, we then obtain a diagram
	\begin{diagram}[small]
	\cF(c)  & \rTo^{M_c} & \cF'(c) \\
	\uTo<{\cF(f)}			& \stackrel{M_{f}}{\Rightarrow} & \uTo>{\cF'(f)} \\
	\cF(c')  & \rTo_{M_{c'}} & \cF'(c') \\
	\end{diagram}
	with $M_f$ given by arrows of the form $M \cF(f) Y \to \cF'(f) M Y$.
\end{exm}

\subsection{Limits and adjunctions}

Consider a Grothendieck prefibration $\cE \to \cC$ over a base $\cC$. Let us
study the question when the category of sections $\Sect(\cC, \cE)$ admits limits
 or colimits. As a related question, given a pullback square of fibrations
\begin{diagram}[small]
F^* \cE & \rTo & \cE \\
\dTo		&			&	\dTo \\
\cD 	& \rTo^F & \cC \\
\end{diagram}
we ask if the natural restriction functor
$F^* : \Sect(\cC,\cE) \to \Sect(\cD,\cE)$ admits an adjoint.


\subsubsection*{Basic results}

\begin{opr}
	A functor $\cE \to \cC$ is \emph{fibrewise-complete} if every fibre $\cE(c)$ is complete. Likewise, $\cE \to \cC$ is \emph{fibrewise-cocomplete} if every fibre $\cE(c)$ is cocomplete.

	A fibration, opfibration, prefibration or preopfibration is fibrewise complete
	or cocomplete if it is true on the level of the underlying functor.
\end{opr}

\begin{prop}
	\label{prefibrationcolimits}
	Let $\cE \to \cC$ be a prefibration which is fibrewise cocomplete. Then the category $\Sect(\cC,\cE)$ is cocomplete, with colimits calculated fibrewise. The dual result concerns limits in the category of sections of a complete preopfibration.
\end{prop}
\proof Let $S_\bullet: I \to \Sect(\cC,\cE)$ be a diagram of sections,
$$
(i,c) \in I \times \cC \, \, \, \mapsto \, \, \, S_i(c) \in \cE(c).
$$
We then define $(\colim_I S_\bullet)(c) = \colim_I S_i (c)$, that is, the colimit of $S_\bullet (c): I \to \cE(c)$ in the fibre $\cE(c)$. Take a morphism $f:c \to d$, it then suffices to construct
\begin{equation}
\label{colimitsglue}
(\colim_I S_\bullet)(c) \to f^* (\colim_I S_\bullet)(d)
\end{equation}
for some choice of a cartesian morphism $f^*(\colim_I S_\bullet)(d) \to (\colim_I S_\bullet)(d)$. If we choose cartesian morphisms for each $i \in I$, obtaining the diagram
$$
f^* S_\bullet(d): I \to \cE(c), \, \, \, i \mapsto f^* S_i(d),
$$
then we have the canonical morphism
$$
\colim_I f^* S_\bullet(d) \to f^* (\colim_I S_\bullet(d))
$$
induced by the colimit property. Combining it with the map $\colim_I S_\bullet(c) \to \colim_I f^* S_\bullet(d)$ induced by the section structure of $S_\bullet$, we get the map (\ref{colimitsglue}). One can check that the induced maps are compatible with the composition of morphisms in $\cC$ in a suitable way. We leave it to the reader: everything follows, in essence, from the universality of maps from a colimit.

Let $X \in \Sect(\cC,\cE)$ be a section, and denote by $c^*X : I \to \Sect(\cC,\cE)$ the constant diagram valued at $X$. Given a map $S_\bullet \to c^* X$, we want to construct an adjoint map $\colim_I S_\bullet \to X$. First, we can construct, fibre by fibre, the maps
$$
\colim_I S_\bullet(c) \to X(c).
$$
For a morphism $f: c \to d$, we can then draw the diagram
\begin{diagram}[small]
\colim_I S_\bullet(c)	&	\rTo	&	\colim_I f^* S_\bullet(d)	&	\rTo	&	f^*\colim_I S_\bullet(d) \\
\dTo								&			&				\dTo							&	&	\dTo \\
X(c)								&	\rTo	&	f^* X(d)		& \rTo^=						&	f^*X(d) \\
\end{diagram}
The left square commutes because $S_\bullet \to c^* X$ is a morphism of
sections, the right square commutes due to the universal property of
colimits. We thus see that the family of fibrewise maps gives a
morphism of sections $\colim_I S_\bullet \to X$. The verification in the other
direction is similar.
\endproof

Given a pullback square of prefibrations,
\begin{diagram}[small]
F^* \cE & \rTo & \cE \\
\dTo		&			&	\dTo \\
\cD 	& \rTo^F & \cC, \\
\end{diagram}
the assignment $S \mapsto  S \circ F$ defines a functor $F^*: \Sect(\cC,\cE) \to \Sect(\cD,\cE)$. One would tentatively write, then, the left adjoint $F_!$ to $F^*$ as a certain colimit over the comma category $F / c$. However, the fibration structure does not permit for sensible formulas to appear. What remains true is the following:
\begin{prop}
	\label{pushforwardalongopfibration}
	Let  $\cE \to \cC$ be a fibrewise-cocomplete prefibration, and
	\begin{diagram}[small]
	F^* \cE & \rTo & \cE \\
	\dTo		&			&	\dTo \\
	\cD 	& \rTo^F & \cC \\
	\end{diagram}
	be a pullback square. Assume that $F: \cD \to \cC$ is an opfibration.
	Then $F^*: \Sect(\cC,\cE) \to \Sect(\cD,\cE)$ admits a left adjoint
	$F_!$, which can be calculated as
	$$
	F_! T(c) = \colim_{\cD(c)} T|_{\cD(c)}.
	$$
\end{prop}
\proof Straightforward and similar to Proposition \ref{prefibrationcolimits}.
Note that $F: \cD \to \cC$ being an opfibration implies that the natural
functor $\cD(c) \to F/ c$ admits a left adjoint and is hence cofinal. \endproof

\subsubsection*{Locally Noether categories}



In what follows, we shall use the words ``sequence'' and ``chain'' interchangeably.
\begin{opr}
	Let $\cC$ be a category, and $c \in \cC$ be an object. We say that $c$ is \emph{$k$-bounded from the right} for some $k \in \bN$ if any sequence of $n$ morphisms starting with $c$,
	$$
	c \longrightarrow c_1 \longrightarrow ... \longrightarrow c_n
	$$
	contains at least $n - k$ isomorphisms so long as $n>k$. Dually, $c$ is $k$-{bounded from the left} if any sequence of $n$ morphisms ending with $c$,
	$$
	c_n \longrightarrow c_{n-1} \longrightarrow ... \longrightarrow c_1 \longrightarrow c
	$$
	contains at least $n - k$ isomorphisms so long as $n>k$
\end{opr}

We shall often say ``bounded'' without being precise about the direction when it leads to no confusion.

\begin{opr}
	\label{noethercatdefinition}
	A category $\cC$ is called \emph{locally Noetherian}, or simply a \emph{Noether category} if for each object $c \in \cC$ there exists a number $k$, such that $c$ is $k$-bounded from the right.

	Dually, a category $\cC$ is called \emph{locally Artinian}, or simply an \emph{Artin category} if for each object $c \in \cC$ there exists a number $k$, such that $c$ is $k$-bounded from the left.
\end{opr}

\begin{rem} \label{noethertoartin} Evidently, if $\cC$ is a Noether category,
then $\cC^\op$ is an Artin category. We shall henceforth stick with
the Noether case in our considerations, but all the results can of course
be dualised for the Artin case.
\end{rem}

For a Noether category $\cC$ and $c \in \cC$, denote by $|c| \geq 0$ the minimal such $k$ so that $c$ is $k$-bounded from the right.

\begin{lemma}
	\label{noethergradingproperties}
	For $c, c' \in \cC$, if $|c| < |c'|$, then $\cC(c,c') = \emptyset$. If $|c| = |c'|$ and there is a map $c \to c'$, then it is an isomorphism. In particular, any endomorphism of $c$ is an isomorphism.
\end{lemma}
\proof Let $c' \to c'_1 \to ... \to c'_{|c'|}$ be a chain starting with $c$ of
length $|c'|$ such that no map in the sequence is an isomorphism.
If there is a map $c \to c'$ in $\cC$, composing with it would yield a
sequence of maps of length $|c'|+1$ starting from $c$.

Thus, if $|c| < |c'|$, we have a sequence of non-invertible maps of length $|c'|+1$ starting from $c$, out of which at least $|c'|$ maps are non-invertible, and this is impossible. If $|c| = |c'|$, having a map $c \to c'$ becomes only possible if it is an isomorphism. \endproof

We thus have a degree function $c \mapsto |c|$, which can be considered as a contravariant functor $|-|: \cC^\op \to \bN$ to the category $\bN$ of natural numbers and unique morphisms in positive direction.


\begin{ntn}
	For a Noether category $\cC$, denote by $\cC_n$ the subcategory of objects $c$
	such that $|c| \leq n$. There is an induced filtration
	$\cC_0 \subset \cC_1 \subset ... \subset \cC_n \subset ... \subset \cC$. Denote
	also by $\cG_n$ the subcategory of $\cC$ consisting of $c$ with $|c| = n$.
	Lemma \ref{noethergradingproperties} implies that $\cG_n$ is a groupoid.
\end{ntn}

Let $\cE \to \cC$ be a prefibration. For $x \in \cC$, if $\cD$ is a subcategory
$x \backslash \cC$, then the prefibration property implies the existence of
a functor $Res_x: \cE|_\cD \to \cE(x)$. An object $Y \in \cE|_\cD$ living over
$f: x \to y$ of $\cD$ is sent to $f^* Y$ where $f^*Y \to Y$ is a cartesian map.
The choice of $Res_x$ is unique up to a unique isomorphism.

Let $S$ be a section over $\cC_{n-1}$. Consider the limit $\lim_{c \backslash \cC_{n-1}} Res_c S$ where $c \in \cG_n$. Since the maps $c \to c'$ are isomorphisms for $|c| = |c'|$, we naturally have $\cE(c) \cong \cE(c')$ (see Convention \ref{convprefibrationareisofibration}) and we get a canonically determined map $\lim_{c \backslash \cC_{n-1}} Res_c S \to \lim_{c' \backslash \cC_{n-1}} Res_{c'} S$.

\begin{opr}
\label{noethermatchingsystem}
	Let $\cE \to \cC$ be a prefibration over a Noether category $\cC$ and $S \in \Sect(\cC_{n-1},\cE)$. The \emph{$n$-th matching system of $S$}, denoted $\Mat_n S$, is the section
	$$
	\Mat_n S: \cG_n \to \cE|_{\cG_{n}}, \, \, \, \, c \mapsto \lim_{c \backslash \cC_{n-1}} Res_c S \, \in \, \cE(c)
	$$
	of the prefibration $\cE \to \cG_{n}$, assuming that all the necessary limits exist.
\end{opr}

The assignment $S \mapsto \Mat_n S$ defines a functor $\Mat_n: \Sect(\cC_{n-1}, \cE) \to \Sect(\cG_n, \cE)$.



\begin{prop}
	\label{noetherfunctorextension}
	There is a comma square
	\begin{diagram}[small]
	\Sect(\cC_{n},\cE) & \rTo & \Sect(\cG_{n},\cE) \\
	\dTo							&	\Leftarrow	&			\dTo>= \\
	\Sect(\cC_{n-1},\cE) & \rTo_{\Mat_n} & \Sect(\cG_{n},\cE) \\
	\end{diagram}
	making $\Sect(\cC_{n},\cE)$ into the comma category $\Sect(\cG_n,\cE) / \Mat_n$. In other words, the assignment
	$$
	Y \in \Sect (\cC_{n},\cE)  \, \, \, \, \mapsto \, \, \, \, (Y|_{\cC_{n-1}},Y|_{\cG_n}, Y|_{\cG_n} \to \Mat_n Y|_{\cC_{n-1}}) \in  \Sect(\cG_n,\cE) / \Mat_n
	$$
	is an equivalence of categories.
\end{prop}
\proof Assume that we are given a section $S$ on $\cC_{n-1}$ and a map
$X \to \Mat_n S$ of sections in $\Sect(\cG_n , \cE)$. We show how to construct
a new section $\tilde S: \cC_n \to \cE$.  For an object $c \in \cC_n$ of
$|c|=n$, there are two kinds of maps: $c \to c'$ with $|c'| = n$ and
$c \to c''$ with $|c''| < n$. The first ones are isomorphisms of $\cG_n$
and are included in $X$ as part of the data. The map $X \to \Mat_n S$
then provides us with morphisms $X(c) \to S(c'')$ in a manner compatible
with $\cG_n$.\endproof

Let $I$ be a small category and denote by $X_\bullet: \in \Sect(\cR,\cE)^I$ a diagram of sections,
$$
(x, i) \, \, \, \mapsto \, \, \, X_i(x).
$$

If the fibre $\cE(x)$ admits limits, we may compute the limit of the functor $i \mapsto X_i(x)$, which we denote $\lim_I (X_\bullet(x))$. We would now like to conclude if the limit of $X_\bullet$, denoted $\lim_I X_\bullet$, exists globally in $\Sect(\cC,\cE)$.

\begin{prop}
	\label{noetherlimits}
	Let $\cC$ be a Noether category and $\cE \to \cC$ a Grothendieck prefibration
	with complete fibres. Then the category of sections $\Sect(\cC,\cE)$ admits
	limits, and moreover, for each $X_\bullet \in \Sect(\cC,\cE)^I$ and
	an object $x$ with $|x|=n$, there is the following pullback square:
	\begin{equation}
	\label{limitinductiondiagram}
	\begin{diagram}[small]
	(\lim_I X_\bullet) (y) & \rTo & \lim_I (X_\bullet (y))		\\
	\dTo						&				&			\dTo						\\
	\Mat_n (\lim_I X_\bullet) (y) & \rTo & \lim_I (\Mat_n X_\bullet) (y).	\\
	\end{diagram}
	\end{equation}
	where $\Mat_n X_\bullet: \cG_n \times I \to \cE$ is the functor $(y,i) \mapsto (\Mat_n X_i)(y)$.
\end{prop}
\proof For each $x$ with $|x|= 0$ we define $(\lim_I X_\bullet) (x) = \lim_I (X_\bullet (x))$, that is we take the limit in the corresponding fibre $\cE(x)$. Since there are no maps out of objects of degree zero, and $\cE(x) \cong \cE(x')$ for $x \cong x'$, we get a well-defined section $\cC_0 \to \cE$.

Having specified $(\lim_I X_\bullet)$ on $\cC_{n-1}$, the diagram (\ref{limitinductiondiagram}) tells us precisely how to define the value $(\lim_I X_\bullet) (y)$ for $y \in \cG_n$. The right vertical arrow exists as a limit of the natural map $X_\bullet(y) \to (\Mat_n X_\bullet)(y)$. The bottom horizontal arrow exists because, by induction, there are natural maps $(\lim_I X_\bullet) (x) \to X_i(x)$ for $x \in \cC_{n-1}$.
These maps induce $\Mat_n (\lim_I X_\bullet) (y) \to (\Mat_n X_i) (y)$ and then,
consequently, we get a map to $\lim (\Mat_n X_\bullet) (y)$.



To verify that the constructed section $Y = \lim_I X_\bullet$ is the limit in $\Sect(\cC,\cE)$, proceed by induction (which is trivial in degree zero) and consider a map $c^* Z \to X_\bullet$, where $c^* Z$ is the constant $I$-section valued at $Z: \cC_{n} \to \cE$. For each $y$ with $|y|=n$, we then get the following diagram:
\begin{diagram}[small]
Z(y)					&		&	\rTo	& &	\lim_I (X_\bullet (y)) \\
\dTo					&			&				&		&	\dTo				\\
\Mat_n Z (y)	& \rTo & \Mat_n Y(y) & \rTo&  \lim_I (\Mat_n X_\bullet) (y)	\\
\end{diagram}
which is commutative because it is simply a factoring of the commutative diagram
\begin{diagram}[small]
Z(y)			 	&	\rTo	&	\lim_I (X_\bullet (y)) \\
\dTo				 &			&			\dTo				\\
\Mat_n Z (y) &  \rTo &  \lim_I (\Mat_n X_\bullet) (y)	\\
\end{diagram}
where the factoring $\Mat_n Z (y) \to \Mat_n Y(y) \to \lim_I (\Mat_n X_\bullet) (y)$
exists due to the limit property of $Y$ on $\cC_{n-1}$. We thus get the
commutative square
\begin{diagram}[small]
Z(y)				&	\rTo	&	\lim_I (X_\bullet (y)) \\
\dTo			&				&			\dTo				\\
\Mat_n Y(y) & \rTo & \lim_I (\Mat_n X_\bullet) (y)	\\
\end{diagram}
which, by the pullback property of the diagram (\ref{limitinductiondiagram}),
supplies us with $Z(y) \to Y(y)$, as desired. \endproof

Proposition \ref{noetherfunctorextension} can be relativised.
Recall the following notions \cite[Definition 1.33]{GROTH}:

\begin{opr}
	\label{openimmersion}
	\label{closedimmersion}
	\label{noethercosieveinclusion}
	A functor $F: \cD \to \cC$ is
	\begin{itemize}
		\item An open immersion if it is full, faithful, injective on objects,
		and for each $f:c \to F(d)$ of $\cC$ there exists a (unique) map
		$\tilde f:d' \to d$ in $\cD$ covering $f$.
		\item An closed immersion if it is full, faithful, injective on objects,
		and for each $f:F(d) \to c$ of $\cC$ there exists a (unique) map
		$\tilde f:d \to d'$ in $\cD$ covering $f$.
	\end{itemize}
\end{opr}

Recall that, for $c \in \cC$, a cosieve is a subcategory
$S \subset c \backslash \cC$ closed under postcomposition:
$f:c \to c' \in S$ implies that $g f$ is in $S$ for any
$g: c' \to c''$ of $\cC$.

\begin{lemma}
\label{lemmaopenimmersiontfae}
	For a functor $F: \cD \to \cC$ injective on objects, the following
	are equivalent
	\begin{itemize}
		\item $F$ is a closed immersion,
		\item $F$ is a faithful isofibration (Definition
		\ref{isofibrationdefinition}), and for each $d \in \cD$, the essential
		image of $d \backslash \cD$ in $F(d) \backslash \cC$ is a cosieve.
		\item $F$ is a fully faithful Grothendieck opfibration with discrete fibres.
	\end{itemize}
	The dual is true for an open immersion.
\end{lemma}
\proof Clear.   \endproof

In particular, let $c \in \cC$ be an object not contained in the image of $F$. Then $\cC(F(d),c) = \emptyset$ for any $d \in \cD$. Thus, at most, there are only morphisms going out of $c$ to $\cD$.

Let $\cC$ be a Noether category and $F: \cD \to \cC$ a closed immersion.
In what follows, we identify $\cD$, which is also a Noether category, with its
image in $\cC$.

\begin{ntn}
	\label{noetherinductivenotation}
	Define $\cD_n$ to be the subcategory consisting of $\cD$ and all the objects $c \in \cC$ not belonging to $\cD$ with $|c| \leq n$. Denote by $F_n: \cD \to \cD_n$ the inclusion functor. There is also an inclusion $\cD_n \to \cC$ which we leave unnamed. Finally, denote by $\cG_n$ the subcategory of $\cD_n$ consisting of those objects $c$ which do not belong to $\cD_{n-1}$.

	For an object $c \in \cC$ (usually assumed to be outside in $\cG_n$) we can define the category $c \backslash \cD_{n-1}$ as the usual comma category for the inclusion $\cD_{n-1} \to \cC$: its objects are maps $c \to d$ in $\cC$, where $d$ belongs to $\cD_{n-1}$.
\end{ntn}

As usual for comma categories and prefibrations, we get the restriction functor $Res_c: \cE|_{c \backslash \cD_{n-1}} \to \cE(c)$.

\begin{prop}
	\label{noetherrightadjoint}
	Let $F:\cD \to \cC$ be a closed immersion of Noether categories and
	$\cE \to \cC$ be a prefibration with complete fibres.
	Then any section $X \in \Sect(\cD,\cE)$ admits a right Kan extension
	$Ran_F X \in \Sect(\cC, \cE)$ which restricts to right Kan extensions
	$Ran_{F_n} X \in \Sect(\cD_n, \cE)$ of $X$ along $F_n: \cD \to \cD_n$.
	Moreover, $F^* Ran_F X \cong X$ and for  any $x \in \cG_n$,
	\begin{equation}
	\label{noetherraninductive}
	(Ran_{F_n} X) (x) = \lim_{x \backslash \cD_{n-1}} Res_x \circ Ran_{F_{n-1}} X
	\end{equation}
	where we implicitly restrict $Ran_{F_{n-1}} X$ to $x \backslash \cD_{n-1}$ along the evident projection.
\end{prop}
\proof
We construct $Ran_{F_n} X$ for each value of $n$ by induction. For $n=0$, the only objects of $x \in \cD_0$ which are not in $\cD$ are those which admit no non-invertible maps out of themselves, since $|x|=0$. We thus pose $(Ran_{F_0} X) (x)$ to be a terminal object of $\cE(x)$. The formula (\ref{noetherraninductive}) then explains how to carry on the induction: for $x, y \in \cD_n$ which are not in $\cD_{n-1}$, the maps $x \to y$,
if exist, are invertible, and the construction of $(Ran_{F_n} X) (x) \to (Ran_{F_n} X) (y)$ is thus as trivial as in Proposition \ref{noetherfunctorextension}. Finally, each object (or a morphism, or a composition of morphisms) of $\cC$ belongs to some $\cG_n$, which permits us to define $Ran_F X$ on the whole of $\cC$.

Since $F$ is a closed immersion, $F^* Ran_F X$ is verified, using
(\ref{noetherraninductive}), to be isomorphic to $X$.
Let $T \in \Sect(\cD, \cE)$ be a section
and assume we have a map $\alpha: F^* T \to X$. We would like now to obtain
a (canonical) morphism $\beta: T \to Ran_F X$. Assume by induction
(which is again trivially initiated for objects of zero degree)
that we obtained this map for all $c \in \cD_{n-1}$ in a compatible fashion.
Let now $x$ be an object of $\cG_n$. There is a diagram in $\cE(x)$ of the form
$$
T(x) \to \lim_{x \backslash \cD_{n-1}} Res_x \circ T \to \lim_{x \backslash \cD_{n-1}} Res_x \circ Ran_{F_{n-1}} X = Ran_{F_n} X (x)
$$
where, if needed, both $T$ and $Ran_{F_{n-1}} X$ are restricted to
$x \backslash \cD_{n-1}$. The first map exists due to the section structure
of $T$, the second map is given by the inductive assumption, and together
they provide $T(x) \to Ran_{F_n} X (x) = Ran_{F} X (x)$. The described
assignment is a bijection, as verified quite easily by applying $F^*$.\endproof

The assignment $X \mapsto Ran_F X$ thus defines a fully faithful functor $F_* : \Sect(\cD,\cE) \to \Sect(\cC,\cE)$ right adjoint to $F^*$.

Consider a closed immersion $F: \cC' \to \cC$ and an object $c \in \cC$. One can form the following pullback square in $\Cat$
\begin{diagram}[small]
c \backslash \cC' & \rTo^{\pi'} & \cC'  \\
\dTo<{F_c}			&					&	\dTo>{F} \\
c \backslash \cC & \rTo_\pi & \cC
\end{diagram}
with $c \backslash \cC'$ coinciding with the usual comma category $c \backslash F$. Moreover, one can verify that each category in this diagram is Noether, with all functors preserving the degrees and
the vertical ones, $F$ and $F_c$, being closed immersions
(the functors $\pi$ and $\pi'$, while being discrete Grothendieck
fibrations, are merely faithful).

If we are given a fibrewise complete prefibration over $\cC$, then there is the following induced $2$-diagram
\begin{diagram}[small]
\Sect(c \backslash \cC',\cE) & \lTo^{\pi'^*} & \Sect(\cC',\cE)  \\
\dTo<{F_{c,*}}			&		\Leftarrow			&	\dTo>{F_*} \\
\Sect(c \backslash \cC,\cE) & \lTo_{\pi^*} & \Sect(\cC,\cE).
\end{diagram}

\begin{prop}
	\label{noetherbasechange}
	In the diagram above, the map $\pi^* F_* \to F_{c,*} \pi'^*$ is an isomorphism.
\end{prop}

We prove it by induction, forming, for each $c \in \cC$, denote by $\cC'_n$ and $(c \backslash \cC')_n$ the induction categories as in Notation \ref{noetherinductivenotation}, with $\pi_n: (c \backslash \cC')_n \to \cC'_n$. being the projection functor. One can see that, moreover, $(c \backslash \cC')_n \cong c \backslash \cC'_n$.
Then Proposition \ref{noetherbasechange} will follow from

\begin{prop}
	\label{noetherbasechangeinductive}
	Let $F:\cC' \to \cC$ be a closed immersion of Noether categories
	and $\cE \to \cC$ be a prefibration with complete fibres.
	Then for each $n$ the $2$-square
	\begin{diagram}[small]
	\Sect(c \backslash \cC',\cE) & \lTo^{\pi'^*} & \Sect(\cC',\cE)  \\
	\dTo<{Ran_{F_{c,n}}}			&		\Leftarrow			&	\dTo>{Ran_{F_n}} \\
	\Sect(c \backslash \cC'_n,\cE) & \lTo_{\pi_n^*} & \Sect(\cC'_n,\cE).
	\end{diagram}
	commutes up to an isomorphism.
\end{prop}
\proof
We shall proceed by induction on $n$. For $n=0$, the extension to objects of
degree zero outside of $\cC'$ or $c \backslash \cC'$ is given by terminal
objects, hence the isomorphism is trivial. Take now an object of
$c \backslash \cC'_n$, represented by a map $c \to d$ with $d$ outside of
t i$\cC'$ and the degree of $|c \to d| = |d|$ equal to $n$. We can then write that
$$
\pi^*_n Ran_{F_n} X (c \to d) = Ran_{F_n} X (d) = \lim _{d \backslash \cC'_{n-1}} Res_d \pi_{n-1}^* Ran_{F_{n-1}} X
$$
with $\pi_{n-1}$ here being the functor $d \backslash \cC'_{n-1} \to \cC'_{n-1}$, and also that
$$
F_{c,*} \pi'^* X(c \to d) = \lim_{(c \to d) \backslash (c \backslash C'_{n-1})} Res_{c \to d} Ran_{F_{c,n-1}} \pi'^* X  \cong \lim_{d \backslash C'_{n-1}} Res_{d} Ran_{F_{c,n-1}} \pi'^* X
$$
where in the middle term one more restriction is implicit. By induction, $$\pi_{n-1}^* Ran_{F_{n-1}} X \to  Ran_{F_{c,n-1}} \pi'^* X$$ is an isomorphism, which induces the isomorphism between the two limit expressions above.  \endproof

\subsection{Factorisation systems and semifibrations}

\begin{opr}
	\label{factorisationsystemdefinition}
	\label{factorisationcategorydefinition}
	A \emph{factorisation system} on a category $\cC$ consists of a
	pair of subcategories $\sL,\sR \subset \cC$ containing all
	isomorphisms of $\cC$, such that any morphism $f: c \to c'$ in $\cC$
	can be decomposed as
	\begin{equation}
	\label{factorisationequation}
	f:c \stackrel l \longrightarrow c'' \stackrel r \longrightarrow c'
	\end{equation}
	with $l \in \Mor \sL$ and $r \in \Mor \sR$. This factorisation must be
	moreover unique up to unique isomorphism.

	In this work, a \emph{factorisation category} will denote a triple $(\cC,\sL,\sR)$ of a category together with a factorisation system $(\sL,\sR)$.
\end{opr}

When clear, we shall simply refer to a factorisation category $(\cC,\sL,\sR)$ as $\cC$. Due to the isomorphism condition $\sL$ and $\sR$ contain all the objects of $\cC$. We shall often refer to $\sL$ as the left class of maps, and to $\sR$ as the right class of maps.

\begin{opr}
	A strict \emph{factorisation functor} $F: (\cC',\sL',\sR') \to (\cC,\sL,\sR)$ is a functor $\cC' \to \cC$ such that $F(\sL') \subset \sL$ and $F(\sR') \subset \sR$. We shall occasionally denote by $F_L: \sL' \to \sL$ and $F_R: \sR' \to \sR$ the induced functors.
\end{opr}

An important class of factorisation categories is given by Reedy categories. To repeat,
\begin{opr}
	A \emph{Reedy category} $\cR$ is a factorisation category $(\cR,\cR_-,\cR_+)$ together with a degree function $deg: \cR \to \bN$ taking values in natural numbers, such that
	\begin{enumerate}
		\item the isomorphisms of $\cR$ are identities,
		\item the non-identities of $\cR_-$ lower the value of $deg$,
		\item the non-identities of $\cR_+$ raise the value of $deg$,
	\end{enumerate}
\end{opr}

It is implied that $\cR_-$ is locally Noether and $\cR_+$ is locally Artin.
For the literature concerning Reedy categories, see \cite{LUHTT,HIRSCHHORN,HOVEY,RIEHL}.
We assume the degree function to be taking values in natural numbers. While this suffices
for most practical examples, our choice excludes from consideration the case $\cR = \beta$
for an arbitrary ordinal $\beta$. The latter has some importance for the theory-building
of Section \ref{subsectionreedycomp}, but is a relatively mild case and will be treated
by hand once needed.
Henceforth, we shall also be implicit about the degree
function in our notation.

The following definition concerns the way factorisation functors interact with the factorisations (\ref{factorisationequation}).

\begin{opr}
	\label{factorisationrightclosedfunctor}
	\label{factorisationleftclosedfunctor}
	Let $F:(\cC',\sL',\sR') \to (\cC,\sL,\sR)$ be a factorisation functor. We say
	that $F$ is \emph{right-closed} if for any $\cC$-map of the form
	$c \to F(c')$, the $(\sL,\sR)$-factorisation of this map can be chosen as
	$$
	c \stackrel l \longrightarrow F(c'') \stackrel{F(r)}{\longrightarrow} F(c')
	$$
	with $r: c'' \to c'$ belonging to $\sR'$.
	Dually, $F$ is \emph{left-closed}, if for any $\cC$-map of the form
	$F(c') \to c$, the $(\sL,\sR)$-factorisation of this map can be chosen as
	$$
	F(c') \stackrel{F(l)}{\longrightarrow} F(c'') \stackrel r \longrightarrow c
	$$
	with $l: c' \to c''$ belonging to $\sL'$.
\end{opr}

\subsubsection*{Semifibrations}

\begin{opr}
	\label{semifibrationdefinition}
	Let $(\cC,\sL,\sR)$ be a factorisation category. A functor $p: \cE \to \cC$ is called a \emph{semifibration} over $\cC$ if it is an isofibration and the following conditions are satisfied.
	\begin{enumerate}
		\item For any $l: c \to c'$ in $\sL$ and $Y$ with $p(Y) = c'$ there exists a cartesian (Definition \ref{cartesianmorphismdefinition}) lift $\lambda:Y' \to Y$ of $l$.
		\item For any $r: x \to y$ in $\sR$ and $X$ with $p(X) = x$ there exists an opcartesian lift $\rho:X \to X'$ of $r$.
		\item For any $\alpha: X \to Y$ of $\cE$ such that $p(\alpha)$ decomposes as $$p(X) \stackrel r \longrightarrow c \stackrel l \longrightarrow p(Y)$$ with $r \in \sR$ and $l \in \sL$, we require that $\alpha$ factors as
		\begin{equation}
		\label{janustriplefactorisation}
		\alpha: X \stackrel \rho \longrightarrow X' \stackrel \phi \longrightarrow Y' \stackrel \lambda \longrightarrow Y
		\end{equation}
		with $\rho: X \to X'$, being an opcartesian morphism over $r$, $\lambda: Y' \to Y$ being a cartesian morphism over $l$, and $p(\phi) = id_c$.

	\end{enumerate}
\end{opr}

\begin{lemma}
	The third condition of Definition \ref{semifibrationdefinition} is equivalent to the following: for any $\alpha: X \to Y$ of $\cE$ such that $p(\alpha)$ decomposes as $$p(X) \stackrel r \longrightarrow c \stackrel l \longrightarrow p(Y)$$ with $r \in \sR$ and $l \in \sL$, we require that
	$$
	\alpha: X \stackrel \rho \longrightarrow X' \stackrel \phi \longrightarrow Y' \stackrel \lambda \longrightarrow Y
	$$
	with $\rho: X \to X'$, being a morphism over $l$, $\lambda: Y' \to Y$ a morphism over $r$, and $p(\phi) = id_c$.
\end{lemma}

\proof Follows from the universality of op(cartesian) arrows. \endproof

Given a semifibration $p:\cE \to \cC$, If $f: c \to c'$ is a map in $\cL$, then there is a functor $f^*: \cE(c') \to \cE(c)$ naturally induced by cartesian lifts. If $g: x \to y$ is a map in $\cR$, we equally have $g_! : \cE(x) \to \cE(y)$ induced by opcartesian lifts.

\begin{prop}
	\label{janusmate}
	\label{mate}
	Let $p: \cE \to \cC$ be a semifibration over $(\cC,\sL,\sR)$. Then
	\begin{enumerate}

		\item The factorisation (\ref{janustriplefactorisation}) is natural and unique up to unique isomorphism,
		\item Let
		\begin{diagram}[small,nohug]
		x			&		\rTo^f		&		y \\
		\dTo<g	&						&		\dTo>h \\
		z			&		\rTo_k		&		t	\\
		\end{diagram}
		be a commutative diagram with $f,k \in \sL$ and $g,h \in \sR$. We then have a two-square
		\begin{diagram}[small,nohug]
		\cE(x)			&		\lTo^{f^*}		&		\cE(y) \\
		\dTo<{g_!}	&		\Rightarrow		&		\dTo>{h_!} \\
		\cE(z)			&		\lTo_{k^*}		&		\cE(t)	\\
		\end{diagram}
		with the natural transformation $g_! f^* \to k^* h_!$ induced canonically.

	\end{enumerate}
\end{prop}

\proof
The first assertion is clear given the universal properties of cartesian and
opcartesian morphisms.

For the second assertion, take $Y \in \cE(y)$. Then we get the diagram in $\cE$
\begin{diagram}[small,nohug]
Y		& \lTo^{cart} & f^* Y & \rTo^{ocart} & g_! f^* Y \\
&\rdTo_{ocart}&			&						&					\\
&						& h_! Y &	\lTo_{cart} & k^* h_! Y\\
\end{diagram}
with maps labeled as $cart$ being cartesian from the fibration structure over $\sL$, and likewise $ocart$ being opcartesian from the opfibration structure over $\sR$. Then, since $hf = kg$, the composition $f^*Y \to Y \to h_! Y$ lies over $x \stackrel g \to z \stackrel k \to t$, and so, by (3) of Definition \ref{semifibrationdefinition}, it can be decomposed as
$$
f^* Y \to g_! f^*Y  \to k^* h_! Y \to h_! Y
$$
and we get a morphism $g_! f^*Y \to k^* h_! Y$ as desired.
\endproof

Since $\cC$ is a factorisation category, any morphism $x \stackrel g \to z \stackrel k \to t$ with $g$ in $\sR$ and $k$ in $\sL$ can be completed to a diagram
\begin{diagram}[small,nohug]
x			&		\rTo^f		&		y \\
\dTo<g	&						&		\dTo>h \\
z			&		\rTo_k		&		t	\\
\end{diagram}
as in Proposition \ref{janusmate} above. So the base-change property for the transition functors can be obtained if one assumes one of the following.
\begin{lemma}
	\label{semifibexm}
	Let $(\cC,\sL,\sR)$ be a factorisation category and $\cE \to \cC$ be an
	\begin{itemize}
		\item either a fibration over $\cC$ which is a preopfibration over $\sR$,
		\item or an opfibration over $\cC$ which is a prefibration over $\sL$,
	\end{itemize}
	then $\cE \to \cC$ is a semifibration.
\end{lemma}
\proof In the first case, for the diagram
\begin{diagram}[small,nohug]
Y		& \lTo^{cart} & f^* Y & \rTo^{ocart} & g_! f^* Y \\
&\rdTo_{ocart}&			&						&					\\
&						& h_! Y &	\lTo_{cart} & k^* h_! Y\\
\end{diagram}
as before we get that the composition $f^* Y \to Y \to h_! Y$ factors through the cartesian map $k^* h_! Y \to h_! Y$ (as implied by the stronger universal property of cartesian maps in this case \cite{VIST}), so we get a map $f^* Y \to k^* h_! Y$. This map in turn is factored by the opcartesian map $f^*Y \to g_! f^* Y$, and we obtain the $Y$-part $g_! f^* Y \to k^* h_! Y$ of the base-change natural transformation. It can then be used to construct the factorisation of Definition \ref{semifibrationdefinition}. The second case is dual. \endproof

A more general statement in this direction is:

\begin{lemma}
	\label{semifibrationrecognition}
	Let $\cE \to \cC$ be a prefibration over a factorisation category $(\cC,\sL,\sR)$, such that the restriction $\cE|_\sR \to \sR$ is also a preopfibration, and such that the composition of cartesian lifts covering $x \stackrel r \to z \stackrel l \to y$ (with $r$ in $\sR$ and $l$ in $\sL$) is cartesian. Then $\cE \to \cC$ is a semifibration over $(\cC,\sL,\sR)$.
\end{lemma}
\proof In the proof of Lemma \ref{semifibexm}, we need the strong cartesian universal property exactly for arrows covering compositions like $x \stackrel r \to z \stackrel l \to y$. \endproof

\subsection{Limits and adjoints in categories of sections}

In a moment, we shall prove the following:

\begin{prop}
	\label{janussectionslimits}
	Let $(\cC,\sL,\sR)$ be a factorisation category and $\cE \to \cC$ be a semifibration with fibres which are complete and admit arbitrary coproducts. Assume that the category $\Sect(\sL,\cE|_\sL)$ has limits. Then so does the category $\Sect(\cC,\cE)$. Moreover, the restriction functor $\Sect(\cC,\cE) \to \Sect(\sL,\cE)$ preserves limits.

	Dually, if $\cE \to \cC$ has cocomplete fibres and fibrewise products, and $\Sect(\sR,\cE|_\sR)$ admits colimits, then so does the category $\Sect(\cC,\cE)$, and the restriction functor $\Sect(\cC,\cE) \to \Sect(\sR,\cE)$ preserves colimits.
\end{prop}

This can be interpreted as saying that in order to calculate limits in
$\Sect(\cC,\cE)$, it is sufficient to do so in $\Sect(\cL,\cE)$. From now on,
we shall concentrate on the limit part, the colimit part being dual.

\begin{lemma}
\label{lemmarestrictiontocommapreserveslimits}
	Let $c \in \cC$ and consider the undercategory $c \backslash \sL$. Then the functor $u_c^*: \Sect(\sL,\cE) \to \Sect(c \backslash \sL,\cE)$, which is induced along the natural forgetful functor $u_c: c \backslash \sL \to \sL$, preserves limits.
\end{lemma}

\proof The functor $u^*_c$ admits a left adjoint
$$
u^c_!:  \Sect( c \backslash \sL, \cE) \to \Sect ( \sL, \cE)
$$
given by the formula
$(u^c_! X)(c') = \coprod_{\sL (c,c')} X(c')$.
\endproof


For any object $c \in \cC$, the semifibration structure provides us with the restriction functor
$$
Res_c : \cE|_{c \backslash \cL } \to \cE(c).
$$

\proof[Proof of Proposition \ref{janussectionslimits}.]
Let $X_\bullet: I \to \Sect(\cC,\cE)$ be a diagram,
$$
i \in I \, \, \mapsto \, \, ( \, \, c \mapsto X_i(c) \, \, ),
$$
and we would like to construct its limit $Y = \lim_I X_\bullet  \in \Sect(\cC,\cE)$. We write the following expression
$$
Y(c) = \lim X_\bullet (c) = \lim_{c \backslash \sL} Res_c ( \lim^{c \backslash \sL}_I X_\bullet |_{c \backslash \sL} )
$$
where $ \lim^{c \backslash \sL}_I X_\bullet |_{c \backslash \sL}$ is the limit of $X_\bullet |_{c \backslash \sL}$ taken in $\Sect(c \backslash \sL, \cE)$, and we shall henceforth drop
the restriction notation for $X_\bullet$.

Because the category $c \backslash \sL$ has an initial object,
$$
\lim_{c \backslash \sL} Res_c ( \lim^{c \backslash \sL}_I X_\bullet) \cong ( \lim^{c \backslash \sL} X_\bullet)(c \stackrel{id}{\to} c) \cong ( \lim^{\sL} X_\bullet)(c),
$$
so, thanks to Lemma \ref{lemmarestrictiontocommapreserveslimits}, our formula
is just another way to write the limit in $\Sect(\sL,\cE)$.

Suppose $r: c \to d$ is a $\sR$-map. We then need to construct $Y(r): Y(c) \to Y(d)$. The semifibration structure implies the necessity to construct an $\cE(d)$-map $r_! Y(c) \to Y(d)$
for some opcartesian map $Y(c) \to r_! Y(c)$. We note that for each
$\sL$-morphism $l: d \to d'$ the factorisation system
of $\cC$ implies the existence of a unique diagram
\begin{equation}
\label{squarelimits}
\begin{diagram}[small,nohug]
c 			& \rTo^r 		& d \\
\dTo<k &					&	\dTo>l \\
c'			& \rTo^t 		& d' \\
\end{diagram}
\end{equation}

with vertical arrows in $\sL$ and horizontal arrows in $\sR$. In terms of undercategories, we can say that there is an induced functor
$$
F: d \backslash \sL \to c \backslash \sL, \, \, \, \,  (l: d \to d') \, \mapsto \, (k: c \to c').
$$

As usual, given any functor $G: c \backslash \sL \to \cM$ we have a natural map between limits $\lim_{c \backslash \sL} G \to \lim_{d \backslash \sL} F^* G$, provided they exist.
Thus, we see that in order to construct a map $f_1$ in
$$
r_! \lim_{c \backslash \sL} Res_c ( \lim^{c \backslash \sL}_I X_\bullet) \stackrel{f_1}{\longrightarrow} \lim_{d \backslash \sL} Res_d ( \lim^{d \backslash \sL}_I X_\bullet)
$$
we can attempt instead to construct another map $f_2$ in
$$
r_! \lim_{d \backslash \sL}F^* Res_c ( \lim^{c \backslash \sL}_I X_\bullet) \stackrel{f_2}{\longrightarrow} \lim_{d \backslash \sL} Res_d ( \lim^{d \backslash \sL}_I X_\bullet).
$$
In turn, due to the universal property of limits, we may instead try to find a map $f_3$ in
$$
\lim_{d \backslash \sL} r_! F^* Res_c ( \lim^{c \backslash \sL}_I X_\bullet) \stackrel{f_3}{\longrightarrow} \lim_{d \backslash \sL} Res_d ( \lim^{d \backslash \sL}_I X_\bullet).
$$
We can now leave out $\lim_{d \backslash \sL}$ and attempt to construct instead
the morphism $f_4$ of functors
$$
r_! F^* Res_c ( \lim^{c \backslash \sL}_I X_\bullet) \stackrel{f_4}{\longrightarrow} Res_d ( \lim^{d \backslash \sL}_I X_\bullet).
$$
Using the notation of Diagram (\ref{squarelimits}), on $l: d \to d'$, the map $f_4$ would yield
$$
r_! k^* ( \lim^{c \backslash \sL}_I X_\bullet) (c \stackrel k \to c') \stackrel{f_4(l)}{\longrightarrow} l^* (\lim^{d \backslash \sL}_I X_\bullet)(d \stackrel l \to d').
$$
Remembering the base-change (Proposition \ref{janusmate}) morphism
$r_! k^* \to l^* t_!$, and the equalities
$$(\lim^{c \backslash \sL}_I X_\bullet)(c \stackrel k \to c') = (\lim^\sL_I X_\bullet) (c')$$
of Lemma \ref{lemmarestrictiontocommapreserveslimits} and the like for $d,d'$,
we see that instead of $f_4$ we may construct maps
$$
l^* t_! (\lim^\sL_I X_\bullet) (c') \stackrel{f_5(l)}{\longrightarrow} l^* (\lim^\sL_I X_\bullet) (d')
$$
or even simpler, $ t_! (\lim^\sL_I X_\bullet) (c') \to  (\lim^\sL_I X_\bullet) (d')$. Examining  $t_! (\lim^\sL_I X_\bullet) (c')$, we witness, naturally, that there are maps
$$
t_! (\lim^\sL_I X_\bullet) (c') \to t_! X_i (c')  \to X_i (d')
$$
with first arrow being a $t_!$ of the limit projection, and the second
given by the section structure of $X_i$.
We assemble these maps together to get $f_5(l)$ for each $l: d \to d'$,
and in turn, $f_4, f_3, f_2$ and $f_1$.

\smallskip

This defines $Y(r): \lim X_\bullet(c) \to \lim X_\bullet(d)$ for $\sR$-maps
of $\cC$. The factorisation structure on $\cC$ and a tedious verification
(which goes through given all the maps in the reasoning above are canonical
in one way or another) then permits to see that $c \mapsto Y(c)$ is indeed
a section of $\cE \to \cC$ that has the required universal property.
\endproof

Recall that $F$ is a right-closed factorisation functor (Definition \ref{factorisationrightclosedfunctor}) if for any $c \to F(c')$ there is a factorisation
$$
c \stackrel l \longrightarrow F(c'') \stackrel{F(r)}{\longrightarrow} F(c')
$$
with $r: c'' \to c'$ belonging to $\sR' \subset \cC'$. This implies that for each map $r:c_1 \to c_2$ of $\sR$ we have the following diagram
\begin{diagram}[small]
c_1 \backslash \sL' & \rTo^{F_{c_1}} & c_1 \backslash \sL \\
\uTo<{r_{\sL'}}		&							&		\uTo>{r_{\sL}} \\
c_2 \backslash \sL' & \rTo^{F_{c_2}} & c_2 \backslash \sL
\end{diagram}
with functors $r_{\sL'},r_{\sL}$ given by factoring the morphisms. One has to be careful about the pullbacks of $\cE \to \cC$ to this diagram. If we denote by $\pi_1 : c_1 \backslash \sL \to \cC$, $\pi_2 : c_2 \backslash \sL \to \cC$ the evident projections, then the factorisations
\begin{equation}
\label{squareadjunctions2}
\begin{diagram}[small,nohug]
c_1 			& \rTo^r 		& c_2\\
\dTo<k &					&	\dTo>l \\
d_1			& \rTo^{ t } 		&  d_2  \\
\end{diagram}
\end{equation}
which define $r_\sL$ as the assignment $l \mapsto k$, imply that there is a natural transformation $\tau:
 \pi_1 r_\sL \to \pi_2$ with components, given by maps like $t$ in the diagram above, lying in $\sR$.

\begin{lemma} \label{janus2functoriality} Let $p:\cE \to \cC$ be a semifibration over $(\cC,\sL,\sR)$ and $F,G: \cD \to \cC$ be two functors taking values in $\sL$, and $\tau: F \to G$ be a natural transformations with components in $\sR$. Then
	\begin{enumerate}
		\item both $F^* \cE \to \cD$ and $G^* \cE \to \cD$ are prefibrations,
		\item the assignment $(X,d,F(d)=p(X)) \mapsto (\tau(d)_! X,d)$ has the property that $p(\tau(d)_! X)=G(d)$ and defines a (lax) morphism of fibrations $\tau_!: F^* \cE \to G^* \cE$ over $\cD$,
		\item there is an induced functor
		$\tau_!: \Sect(\cD,F^* \cE) \to \Sect(\cD, G^* \cE)$ on the categories of
		sections. Moreover, for each $X \in \Sect(\cC, \cE)$, there is a natural
		(in $X$) map $\tau_! F^* X \to G^* X$.
		\item Let $H: \cD' \to \cD$ be a functor, and assume that there are right adjoints,
		$$
		H^*_F: \Sect(\cD,F^* \cE) \leftrightharpoons \Sect(\cD',F^* \cE) : H^F_*,
		$$
		$$
		H^*_G: \Sect(\cD,G^* \cE) \leftrightharpoons \Sect(\cD',G^* \cE) : H^G_*
		$$
		for the restriction functors $H^*_F, H^*_G$. Then there is a natural map
		$$
		\tau_! H^F_* \longrightarrow H^G_* \tau'_!
		$$
		where $\tau'_!:  \Sect(\cD',H^*F^* \cE) \to \Sect(\cD', H^*G^* \cE)$ is the functor induced as in previous paragraph.
	\end{enumerate}
\end{lemma}
\proof The first statement is clear. For the second we are left to prove
that the assignment $X \mapsto \tau(d)_! X$ is indeed a morphism of
prefibrations. For a map $f: d \to d'$, we can draw the following square
\begin{equation}
\label{janus2functorialitydiagram}
\begin{diagram}[small,nohug]
F d' 			& \rTo^{ \tau(d') } 		& Gd' \\
\dTo<{Ff} &					&	\dTo>{Gf} \\
F d			& \rTo^{ \tau(d) } 		&  Gd \\
\end{diagram}
\end{equation}
Using the fibrewise-cartesian factoring on $F^* \cE$, it remains to see what
happens to the cartesian maps $Ff^* Y \to Y$, $p(Y) = Fd'$. We observe
that the base-change for the diagram above implies the existence of the map
$$
\tau(d)_! Ff^* Y \longrightarrow Gf^* \tau(d')_! Y.
$$
Choosing (or, rather, remembering) a cartesian map
$Gf^* \tau(d')_! Y \to \tau(d')_! Y$, we get the composition
$$
\tau(d)_! Ff^* Y \longrightarrow Gf^* \tau(d')_! Y \to Gf^* \tau(d')_! Y \to \tau(d')_! Y
$$
needed for constructing the functor $\tau_!: F^*\cE \to G^* \cE$.

The functor $\tau_!$ of the third statement is simply induced by the
post-composition with the functor $\tau_!$ of the second statement.
The existence of the natural family of maps $\tau_! F^* X \to G^* X$ happens
for the following reason: on an object $d \in \cD$,
the map $\tau(d)_! X(F(d)) \to X(G(d))$ is supplied by the section
structure of $X$ along the $\sR$-map $\tau(d): F(d) \to G(d)$.

For the fourth statement, consider the diagram
\begin{diagram}[small,nohug]
\Sect(\cD,F^* \cE)	& \rTo^{ H^*_F } & \Sect(\cD',F^* \cE)\\
\dTo<{\tau_!} &					&	\dTo>{\tau'_!} \\
\Sect(\cD,G^* \cE)& \rTo^{ H^*_G } 	& \Sect(\cD',G^* \cE) \\
\end{diagram}
and observe by explicit check that it commutes up to an isomorphism. Hence the sought-after map
$$
\tau_! H^F_* \longrightarrow H^G_* \tau'_!
$$
is given by the usual base-change argument. \endproof

\begin{rem}
The functor $\tau_!: F^*\cE \to G^* \cE$ takes a cartesian maps to
cartesian whenever the base-change map for (\ref{janus2functorialitydiagram})
is an isomorphism.
\end{rem}

We would now like to prove a statement about adjoints similar to Proposition
\ref{janussectionslimits}. Namely, given a semifibration $\cE \to \cC$ and
a right-closed functor $F:\cD \to \cC$, we would like to deduce the existence
of a right adjoint to the pullback functor $F^*: \Sect(\cC,\cE) \to \Sect(\cC,\cE)$ from assuming the existence of one for $F_\sL^*: \Sect(\sL,\cE) \to \Sect(\sL',\cE)$. However, we
shall need to require some additional properties, which will make
harmless the passage to comma categories.

\begin{opr}
	\label{pointwiseradj}
	In the situation above, we say that pull-back $F^*_{\sL}$ admits a \emph{pointwise right adjoint} if
	\begin{enumerate}
		\item the functor $F_\sL^*: \Sect(\sL,\cE) \to \Sect(\sL',\cE)$ admits a right adjoint $F_{\sL,*}$,
		\item for each $c \in \sL$, the pull-back $F^*_{c}: \Sect(c \backslash \sL,\cE) \to \Sect(c \backslash \sL',\cE)$ along the induced functor $F_{c}:c \backslash \sL' \to c \backslash \sL$, admits a right adjoint $F_{c,*}$ and moreover the natural base-change map $\pi^* F_{\sL,*} \to  F_{c,*} \pi'^*$ arising from the square
		\begin{diagram}[small]
		\sL' & \rTo^{F_{\sL}} &  \sL \\
		\uTo<{\pi'}		&							&		\uTo>{\pi} \\
		c \backslash \sL' & \rTo^{F_{c}} & c\backslash \sL
		\end{diagram}
		is an isomorphism.
	\end{enumerate}

	In other words, this means that $ F_{c,*} \pi'^* X$ can be computed as
	$F_{\sL,*} X$ and then restricted again to the comma category.


\end{opr}

\begin{prop}
	\label{semifibrationrightadjoint}
	Let $F: \cC' \to \cC$ be a right-closed factorisation functor, and $\cE \to \cC$ a fibrewise complete semifibration over $\cC$. Assume that the functor $F_\sL^*: \Sect(\sL,\cE) \to \Sect(\sL',\cE)$ admits a pointwise right adjoint $F_{\sL,*}$ in the sense of Definition \ref{pointwiseradj}. Then the functor $F^*: \Sect(\cC,\cE) \to \Sect(\cC,\cE)$
	admits a right adjoint $F_*$  such that the induced $2$-diagram
	\begin{diagram}[small]
	\Sect(\cD,\cE) & \rTo^{F_{*}} & \Sect(\cC,\cE) \\
	\dTo 		&			\Rightarrow				&		\dTo  \\
	\Sect(\sL',\cE) & \rTo^{F_{\sL,*}} & \Sect(\sL,\cE), \\
	\end{diagram}
	(with vertical arrows given by restrictions), is in fact commutative up to an isomorphism.
\end{prop}

We can thus make conclusions about $F_*$ by passing to the left categories and using the functor $F_{\sL,*}$.

\proof We shall proceed in a manner similar to Proposition
\ref{janussectionslimits}. For $c \in \cC$ and $X \in \Sect(\cC',\cE)$, put
$$
Y(c) := F_* X (c) = \lim_{c \backslash \sL} Res_c F_{c,*} (X|_{c \backslash \sL'})
$$
where $F_c : c \backslash \sL' \to c \backslash \sL$ is the functor induced from $F$. Indeed, $Y(c) \cong F_{\sL*}X (c)$, but we will need such a presentation for $Y$ for the proof to work.

Assume given a map $r: c_1 \to c_2$. We need to construct
$$
r_! \lim_{c_1 \backslash \sL} Res_{c_1} F_{c_1,*} (X|_{c_1 \backslash \sL'}) \stackrel{f_1}{\longrightarrow} \lim_{c_2 \backslash \sL} Res_{c_2} F_{c_2,*} (X|_{c_2 \backslash \sL'})
$$

Since $F$ is right-closed, we have the following diagram
\begin{diagram}[small]
c_1 \backslash \sL' & \rTo^{F_{c_1}} & c_1 \backslash \sL \\
\uTo<{r_{\sL'}}		&							&		\uTo>{r_{\sL}} \\
c_2 \backslash \sL' & \rTo^{F_{c_2}} & c_2 \backslash \sL
\end{diagram}
with functors $r_{\sL'},r_{\sL}$ given by factoring the morphisms. One has to be careful about the pullbacks of $\cE \to \cC$ to this diagram. If we denote by $\pi_1 : c_1 \backslash \sL \to \cC$, $\pi_2 : c_2 \backslash \sL \to \cC$ the evident projections, then the factorisations
\begin{equation}
\label{squareadjunctions3}
\begin{diagram}[small,nohug]
c_1 			& \rTo^r 		& c_2\\
\dTo<k &					&	\dTo>l \\
d_1			& \rTo^{ t } 		&  d_2  \\
\end{diagram}
\end{equation}
imply that there is a natural transformation $\tau:                                                                                                                                                                                                                                                                                                                                                                                                                                                                                          \pi_1 r_\sL \to \pi_2$ with components, given by maps like $t$ in the diagram above, lying in $\sR$.

We can thus attempt instead to construct another map $f_2$ in
$$
r_! \lim_{c_2 \backslash \sL} r^*_{\sL} Res_{c_1} F_{c_1,*} (X|_{c_1 \backslash \sL'}) \stackrel{f_2}{\longrightarrow} \lim_{c_2 \backslash \sL} Res_{c_2} F_{c_2,*} (X|_{c_2 \backslash \sL'}).
$$
In turn, due to the universal property of limits, we may instead try to find a map $f_3$ in
$$
\lim_{c_2 \backslash \sL} r_!   r^*_{\sL} Res_{c_1} F_{c_1,*} (X|_{c_1 \backslash \sL'}) \stackrel{f_3}{\longrightarrow} \lim_{c_2 \backslash \sL} Res_{c_2} F_{c_2,*} (X|_{c_2 \backslash \sL'}).
$$
We can now leave out $\lim_{c_2 \backslash \sL}$ and construct instead the morphism $f_4$ of functors
$$
r_!   r^*_{\sL} Res_{c_1} F_{c_1,*} (X|_{c_1 \backslash \sL'})  \stackrel{f_4}{\longrightarrow} Res_{c_2} F_{c_2,*} (X|_{c_2 \backslash \sL'}).
$$
Using the notation of the diagram (\ref{squareadjunctions3}) coming from the factorisation on $\cC$, the map $f_4$ would yield
$$
r_!  k^* F_{c_1,*} (X|_{c_1 \backslash \sL'})(c_1 \stackrel k \to d_1 ) \stackrel{f_4(l)}{\longrightarrow} l^* F_{c_2,*} (X|_{c_2 \backslash \sL'})(c_2 \stackrel l \to d_2).
$$
Remembering the base-change morphism $r_! k^* \to l^* t_!$,
we see that instead of $f_4$ we may construct maps
$$
t_! F_{c_1,*} (X|_{c_1 \backslash \sL'})(c_1 \stackrel k \to d_1  ) \stackrel{f_5(l)}{\longrightarrow} F_{c_2,*} (X|_{c_2 \backslash \sL'})(c_2 \stackrel l \to d_2).
$$
We note that $F_{c_1,*} (X|_{c_1 \backslash \sL'})(c_1 \stackrel k \to d_1  ) = r^*_\sL F_{c_1,*}  (X|_{c_1 \backslash \sL'})(c_2 \stackrel l \to d_2  )$, where $r^*_\sL$ is now the pullback on sections, and see that we are looking for $f_5$ in
$$
\tau_! r^*_\sL F_{c_1,*} (X|_{c_1 \backslash \sL'}) \stackrel{f_5 }{\longrightarrow} F_{c_2,*} (X|_{c_2 \backslash \sL'})
$$
with $\tau_!$ induced from $\tau:                                                                                                                                                                                                                                                                                                                                                                                                                                                                                          \pi_1 r_\sL \to \pi_2$ by Lemma \ref{janus2functoriality}.

There is a base-change map $$r^*_\sL F_{c_1,*} \to F'_{c_2,*} r^*_{\sL'}$$ with components lying the category $\Sect(c_2 \backslash \sL, (\pi_1 r_\sL)^* \cE)$. The prime over the functor $F'_{c_2,*}$ denotes that it is adjoint for the sections of the prefibration $(\pi_1 r_\sL)^* \cE$ and not $\pi_2^* \cE$. Now, apply $\tau_!$ and get
$$
\tau_! r^*_\sL F_{c_1,*} \to \tau_! F'_{c_2,*} r^*_{\sL'} \to F_{c_2,*} \tau'_! r^*_{\sL'}
$$
with the second arrow existing due to the fourth statement of Lemma \ref{janus2functoriality}, with $\tau': \pi'_1 r'_{\sL'} \pi'_2$ be the natural transformation between the evident projections $\pi_1 : c_1 \backslash \sL' \to \cC'$, $\pi_2 : c_2 \backslash \sL' \to \cC'$ and $r_{\sL'}$.

Examining what is remaining we see that to get $f_5$, we may as well construct $f_6$ in
$$
F_{c_2,*} \tau'_! r^*_{\sL'} X|_{c_1 \backslash \sL'} \stackrel{F_{c_2,*} f_6}{\longrightarrow} F_{c_2,*} X|_{c_2 \backslash \sL'},
$$
or, removing $F_{c_2,*}$,
$$
\tau'_! r^*_{\sL'} X|_{c_1 \backslash \sL'} \stackrel{ f_6}{\longrightarrow}  X|_{c_2 \backslash \sL'},
$$
This map, is, however, simply there by the third statement of Lemma \ref{janus2functoriality}, since $X$ is a factual section of a semifibration. If we consider the factorisation diagram defining $r_{\sL'}$,
\begin{equation}
\begin{diagram}[small,nohug]
c_1 			& \rTo^r 		& c_2\\
\dTo<a &					&	\dTo>b \\
F(d_1)			& \rTo^{ F(e) } 		&  F(d_2)  \\
\end{diagram}
\end{equation}
then the map $f_6(b)$ corresponds to $F(e)_! X(F(d_1)) \to X(F(d_2))$. We thus get $f_6$ and reverse all the discussion to get $f_1$. \endproof

\begin{corr}
	Let $F:\cD \to \cC$ be a factorisation right-closed functor such that its restriction $F_\sL: \sL' \to \sL$ is a closed immersion of Noether categories. Then for any fibrewise-complete semifibration $\cE \to \cC$, there is an adjunction $F^*: \Sect(\cC,\cE) \rightleftarrows \Sect(\cD,\cE): F_*$, and the right adjoint can be calculated by restricting to the left parts of the factorisation systems.
\end{corr}

\proof The right adjoint for $F_\sL: \sL' \to \sL$ exists thanks to Proposition \ref{noetherrightadjoint} and is pointwise due to Proposition \ref{noetherbasechange}. \endproof

%% file: rmsf-reedy.tex
\section{Reedy model structures}
\label{chapterreedymodstr}

\subsection{Model categories and localisation}

\begin{opr}
	A \emph{homotopical, or relative, category}, or a \emph{localiser}, is a pair $(\cM,\cW)$ of a category $\cM$ and a subcategory $\cW$
	containing all objects of $\cM$, called the category of weak equivalences.
\end{opr}

The definition of a model category used in this work is the following one:
\begin{opr}
\label{defmodelcat}
	A category $\cM$ carries a \emph{model structure}, or is a
  \emph{model category}, if there are given three subcategories $(\cW,\cC,\cF)$
  containing all objects of $\cM$, called respectively the subcategory of weak
  equivalences, cofibrations and fibrations, such that the following list of
  axioms is satisfied.
	\begin{enumerate}
		\item[M1] (Property of $\cM$) the category $\cM$ admits small limits and colimits.
		\item[M2] The subcategory $\cW$ satisfies $3$-for-$2$: given two composable
    maps $f, g$, if any two elements of the set $\{f, g, gf \}$ are morphisms
    of $\cW$, then so is the third.
		\item[M3] The subcategories $\cW,\cC,\cF$ are stable by retracts: given a commutative diagram
		\begin{diagram}[small]
		A & \rTo^{i_1} & X & \rTo^{r_1} & A \\
		\dTo<{f} & & \dTo<{g} & & \dTo<{f} \\
		B & \rTo^{i_2} & Y & \rTo^{r_2} & B \\
		\end{diagram}
		with $r_1 i_1 = id_A$ and $r_2 i_2 = id_B$, if $g$ belongs to $\cW$ (respectively to $\cC$,$\cF$), then so does $f$.
		\item[M4] In a commutative diagram
		\begin{diagram}[small,nohug]
		A	& \rTo^a	& X	 \\
		\dTo<i	 &			&	\dTo>f		\\
		B & \rTo^b & Y	 \\
		\end{diagram}
		with $i$ in $\cC$ and $f$ in $\cF$, whenever any of $i,f$ is also in $\cW$, there exists a map $p: B \to X$ with $pi = a$ and $fp = b$.
		\item[M5] Any morphism $p: X \to Y$ can be factored as $X \stackrel i \to Z \stackrel f \to Y$ with $i$ in $\cC$ and $f$ in $\cF \cap \cW$, and as  $X \stackrel j \to Z' \stackrel g \to Y$, with $j$ in $\cC \cap \cW$ and $g$ in $\cF$.
	\end{enumerate}
\end{opr}

A functor $F: \cM \to \cN$ of model categories is called \emph{left-derivable}
if it preserves cofibrations and trivial cofibrations, and \emph{left Quillen}
if in addition it is a left adjoint. The notions of \emph{right-derivable} and
\emph{right Quillen} functors are defined dually.

\begin{opr}
	Let $(\cM,\cW)$ be a homotopical category. The \emph{localisation} of $\cM$
  along $\cW$ \cite{DHKS,HOVEY}, is the category which we denote $\cW^{-1} \cM$
  or $\Ho \cM$, together with a functor $p: \cM \to \cW^{-1} \cM$
  such that any functor $F: \cM \to \cN$ which sends
  $\cW$ to isomorphisms of $\cN$, factors through $p$. The factorisation
  is unique up to a canonical isomorphism.
\end{opr}

\begin{prop}[\cite{DHKS,HIRSCHHORN,HOVEY}]
	For a model category $\cM$, the localisation $\Ho \cM$ of $\cM$ along $\cW$ exists and is of the same (set-theoretical) size as $\cM$.
\end{prop}

The higher-categorical localisation is revisited in Section \ref{subsectionreedycomp}.

\subsection{Semifibrations over Reedy categories}
\label{subsectionmodelsemifib}

\begin{opr}
	Let $\cR$ be a Reedy category. A \emph{model semifibration} over $\cR$ is a functor $\cE \to \cR$ such that it is a semifibration over $(\cR,\cR_-,\cR_+)$, each fibre $\cE(x)$ is a model category, and
	\begin{enumerate}
		\item the transition functors along $\cR_-$ are right derivable,
		\item the transition functors along $\cR_+$ left derivable.
	\end{enumerate}
\end{opr}

We shall prove that under some admissibility conditions, the category of sections $\Sect(\cR, \cE)$ carries a model structure.

\bigskip
Recall \cite{HIRSCHHORN,HOVEY} that for each object $x \in \cR$, we have associated latching and matching categories $Lat(x)$ and $Mat(x)$. Let $\cE \to \cR$ be a semifibration. Then for each $x \in \cR$, there are natural restriction functors $L_x:\cE|_{Lat(x)} \to \cE(x)$ and $R_x: \cE|_{Mat(x)} \to \cE(x)$. Indeed, by $L_x$, an object $X \in \cE|_{Lat(x)}$ living over $f:y \to x$ is sent to its opcartesian image $f_! X \in \cE(x)$, and dually for $R_x$.

There are choices involved in constructing $L_x$ and $R_x$; both functors are unique up to a natural isomorphism.

\begin{opr}[cf Definition \ref{noethermatchingsystem}]
\label{defreedylatchingmatching}
	For $S \in \Sect(\cR,\cE)$ and $x$ in $\cR$, we define the latching object of
  $S$ at $x$ to be the following colimit:
	$$
	\Lat_x S  := \colim_{Lat(x)} L_x \circ S|_{Lat(x)}.
	$$
	The matching object of $S$ at $x$ is defined to be the following limit:
	$$
	\Mat_x S  := \lim_{Mat(x)} R_x \circ S|_{Mat(x)}.
	$$

\end{opr}

The latching and matching object constructions are suitably functorial. Denote by $\cR_{< n}$ the subcategory of objects of degree less than $n$. We see that $\Lat_x$ and $\Mat_x$ define functors $\Sect(\cR_{<n} , \cE) \to \cE(x)$.
Now, consider a section $S: \cR_{< n} \to \cE$. Then for each $z$ of degree (up to) $n$, the map $\Lat_z S \to \Mat_z S$ is canonically determined. To see this, we need to supply, for each degree-raising map $g:x \to z$ and each degree-lowering map $k:z \to t$, a map $g_! S(x) \to k^* S(t)$. Since $\cR$ is a Reedy category, we have the following square
\begin{diagram}[small,nohug]
x			&		\rTo^f		&		y \\
\dTo<g	&						&		\dTo>h \\
z			&		\rTo_k		&		t	\\
\end{diagram}
in which the vertical maps raise the degree and the horizontal maps lower the degree. Proposition \ref{janusmate} then implies that we have a natural transformation $g_! f^* \to k^* h_!$. The sought-after map is then defined as the composition
$$
g_! S(x) \to g_! f^* S(y) \to k^* h_! S(y) \to k^* S(t)
$$
with $S(x) \to f^* S(y)$ and $h_! S(y) \to S(t)$ existing because $S$ is a section on $\cR_{< n}$. Combining different maps $g_! S(x) \to k^* S(t)$, we get the map from the colimit to the limit, that is, $\Lat_z S \to \Mat_z S$.

For a section $S: \cR \to \cE$ defined on the whole of $\cR$, we are supplied with maps $\Lat_x S  \to S(x) \to \Mat_x S $ in the fibre $\cE(x)$ which can be seen to factor the canonical map $\Lat_x S \to \Mat_x S$.

\begin{prop}
\label{propsectionlatchingmatching}
	Let $\cE \to \cR$ be a semifibration and $S: \cR_{< n} \to \cE$ be a section defined on objects of degree less than $n$. Then an extension of $S$ to a section on objects $x$ of degree $n$ is equivalent to factoring the canonical maps $\Lat_x S \to \Mat_x S$ as $\Lat_x S \to S(x) \to \Mat_x S$ for each $x$.
\end{prop}
\proof Given any non-trivial map $x \to y$ between two objects of degree $n$, we factor it as $x \to z \to y$, and the corresponding map $S(x) \to S(y)$ is constructed as $$S(x) \to \Mat_x S \to S(z) \to \Lat_y S \to S(y),$$ with the middle maps well defined as $deg z < n$. \endproof

As mentioned above, the assignments $S \mapsto \Lat_x S$ and $S \mapsto \Mat_x S$ define functors from $\Sect(\cR,\cE)$ to $\cE(x)$. Thus, given a map $f: S \to T$ of two sections $S,T \in \Sect(\cR, \cE)$, we get, naturally, two following squares
\begin{diagram}[small]
\Lat_x S 	& \rTo	& S(x)	& \rTo	&	\Mat_x S  \\
\dTo			&			&	\dTo		&			&		\dTo		\\
\Lat_x T  & \rTo & T(x)	& \rTo	& \Mat_x T  \\
\end{diagram}

\begin{opr}
	\label{reedymodelstructureclassesofmapsdefinition}
	A map of sections $f: S \to T$ is a
	\begin{enumerate}
		\item \emph{Reedy cofibration} if the map $ \Lat_x T  \coprod_{\Lat_x S } S(x) \to T(x)$ is a cofibration in $\cE(x)$ for each $x \in \cR$.
		\item \emph{Reedy fibration} if the map $ S(x) \to \Mat_x S  \prod_{\Mat_x T } T(x)$ is a fibration in $\cE(x)$ for each $x \in \cR$.
		\item \emph{Reedy weak equivalence} if it is a fibrewise weak equivalence.
	\end{enumerate}
\end{opr}

\begin{opr}
	\label{admissibilitydef}
	Let $\cR$ be a Reedy category and $\cE \to \cR$ a model semifibration. We call $\cE \to \cR$ \emph{left-admissible} if the following holds. Let $\alpha: A \to B$ be a map
	of sections such that $\Lat_x B  \coprod_{\Lat_x A } A(x) \to B(x)$ is a (trivial) cofibration for each $x \in \cR$. Then for any $y \in \cR$, the map $\Lat_y(\alpha) : \Lat_y A \to \Lat_y B$ is also a (trivial) cofibration.

	Dually, one can define a \emph{right-admissible} model semifibration. A semifibration is called \emph{admissible} if it is both left- and right-admissible.
\end{opr}

We will show that under the left admissibility condition,
the trivial Reedy cofibrations are exactly the maps $A \to B$ such that $\Lat_x B  \coprod_{\Lat_x A } A(x) \to B(x)$ is a (trivial) cofibration for each $x \in \cR$.
The left admissibility condition can be interpreted as ensuring that the functors $\Lat_x: \Sect(\cR,\cE) \to \cE(x)$
are left derivable.

\begin{lemma}
\label{admissibilitycriterion}
A model semifibration $\cE \to \cR$ is left-admissible if for each $x \in \cR$, either of the following
holds:
\begin{enumerate}
	\item the restriction $\cE|_{Mat(x)} \to Mat(x)$ is a fibration (and not
	simply a prefibration) whose transition functors preserve limits,
	\item the category $Mat(x)$ is a disjoint union of categories
	possessing initial objects.
\end{enumerate}
A dual result is valid for right-admissibility.
\end{lemma}
\proof For each $y$ in $\cR$, there are two possible scenarios.
\begin{enumerate}
	\item The functor $\cE|_{Lat(y)} \to Lat(y)$ is an opfibration (not merely a preopfibration) whose transition functors preserve colimits. The restriction to the fibre $\cE(y)$ provides us with a functor $L:\Sect(Lat(y),\cE) \to \Fun(Lat(y),\cE(y))$. A section $S \in \Sect(Lat(y),\cE)$ is sent to $$L(S): (z \stackrel f \to y ) \in Lat(y) \, \, \, \mapsto \, \, \, L(S)(f) \cong f_! S(z) \in \cE(y).$$
	The category $ \Fun(Lat(y),\cE(y))$ has a well-known Reedy structure \cite{HOVEY}. Let us compute the value $Lat_f L(S)$ of the latching object functor at $f: z \to y$ (we write $Lat$ to distinguish from $\Lat$ which we used for sections). Abusing slightly the notation, one can see that, naturally in $S$,
	$$
	Lat_f L(S) \cong
	\colim_{g: t \to z \in Lat(z)} L(S)(t \stackrel g \to z \stackrel f \to y) \cong \colim_{t \to z \in Lat(z)} (fg)_! S(t)$$  $$\cong \colim_{t \to z \in Lat(z)} f_! g _! S(t) \cong f_! \Lat_z (S),
	$$
	where the last two isomorphisms are consequences of the given admissibility condition. One can use this and similar computations to verify that the image of the map $\alpha: A \to B$ in $\Fun(Lat(y),\cE(y))$ is a (trivial) Reedy cofibration. Given that $\Lat_y A \cong \colim_{Lat(y)} L(A)$, the necessary result follows from the classical case \cite{HOVEY}.

	\item The category $Lat(y)$ is a disjoint union of categories with terminal objects. Any colimit over a such category is a coproduct of evaluations at terminal objects of the components. If we suppose that the assertion of the lemma was proven by induction for lesser degrees, the map $\Lat_y(\alpha): \Lat_y A \to \Lat_y B$ will be represented as a coproduct $\coprod X_i \to \coprod Y_i$, where each map $X_i \to Y_i$ is a (trivial) cofibration. Thus $\Lat_y (\alpha)$ is also a (trivial) cofibration. \endproof

\end{enumerate}

\begin{thm}
	\label{reedymodelstructuretheorem}
	Let $\cR$ be a Reedy category and $\cE \to \cR$ an admissible model semifibration. Then the category of sections $\Sect(\cR, \cE)$ carries a model structure given by Reedy cofibrations, Reedy fibrations and Reedy weak equivalences of Definition \ref{reedymodelstructureclassesofmapsdefinition}.
\end{thm}

\begin{rem}
	The admissibility condition is used in the proof of Lemma \ref{reedycofibrationsgivepointwise}. Without admissibility, many aspects of
	the Reedy proof do indeed go through, but one has no control over the
	intersection of the classes of maps, and consequently, over factorisations
	and (co)fibrant replacements.

	The condition $(1.)$ and its dual of Lemma \ref{admissibilitycriterion} has
	the property that if it is true for an object $x$, then it is also true for all
	objects $y$ in its matching (or latching) category.

	The condition $(2.)$ and its dual of Lemma \ref{admissibilitycriterion} are
	related to the notion of fibrant and cofibrant constants \cite{DHKS}.
	The reason that it appears in our setting is somehow dual to that of
	\cite{DHKS}; see Lemma \ref{reedycofibrationsgivepointwise} for details.
\end{rem}

\begin{lemma} The Reedy weak equivalences are stable under retracts and
satisfy the ``three-for-two'' axiom.
\end{lemma}
\proof Clear, by considering what happens in each fibre. \endproof

\begin{lemma}
	\label{retracts}
	Let $f: S \to T$ be a map of sections such that $f$ satisfies one of the properties below:
	\begin{itemize}
		\item For each $x \in \cR$, the map $\Lat_x T  \coprod_{\Lat_x S } S(x) \to T(x)$ is a cofibration,
		\item For each $x \in \cR$, the map $  \Lat_x T  \coprod_{\Lat_x S } S(x) \to T(x)$ is a trivial cofibration,
		\item For each $x \in \cR$, the map $  S(x) \to \Mat_x S  \prod_{\Mat_x T } T(x)$ is a fibration,
		\item For each $x \in \cR$, the map $ S(x) \to \Mat_x S  \prod_{\Mat_x T } T(x)$ is a trivial fibration.
	\end{itemize}
	Then any retract of $f$ also satisfies such a property.
\end{lemma}
\proof Let
\begin{diagram}[small]
A & \rTo^{i_1} & X & \rTo^{r_1} & A \\
\dTo<{f} & & \dTo<{g} & & \dTo<{f} \\
B & \rTo^{i_2} & Y & \rTo^{r_2} & B \\
\end{diagram}
be a retract diagram in $Sect(\cR,\cE)$. The assignment $A \mapsto \Lat_x A$ is functorial in $A$, so it preserves retracts. Then, for $x \in \cR$, there is a diagram $D_1$
\begin{diagram}[small]
A(x)  & \rTo^{i_1(x)}  & X(x)  & \rTo^{r_1(x)}  & A(x)  \\
\uTo  &                & \uTo  &                & \uTo  \\
\Lat_x A & \rTo^{\Lat_x i_1} & \Lat_x X & \rTo^{\Lat_x r_1} & \Lat_x A \\
\dTo  &                & \dTo  &                & \dTo  \\
\Lat_x B & \rTo^{\Lat_x i_2} & \Lat_x Y & \rTo^{\Lat_x r_2} & \Lat_x B \\
\end{diagram}
which can be viewed as a retract diagram in $Fun(I,\cE(x))$, where $I$ is the category $0 \leftarrow 1 \rightarrow 2$. There is also a retract diagram $D_2$
\begin{diagram}[small]
B(x)  & \rTo^{i_2(x)}  & Y(x)  & \rTo^{r_2(x)}  & B(x)  \\
\end{diagram}

For a category $\cD$, let $Ret(\cD)$ be the category of retract diagrams: its objects are pairs of arrows $C \stackrel{i}{\to} D \stackrel{r}{\to} C$ with $r \circ i = id_C$. For any small category $\cJ$ the constant diagram functor $c^*_\cJ: \cD \to Fun(\cJ,\cD)$ induces a functor $Ret(c^*_\cJ): Ret(\cD) \to Ret(Fun(\cJ,\cD))$. If $\cD$ admits small colimits, this functor has a left adjoint $Ret(\colim_\cJ):Ret(Fun(\cJ,\cD)) \to Ret(\cD)$.

In our case, $\cD = \cE(x)$ has small colimits and $\cJ = I$. In addition, $D_1 \in Ret(Fun(I,\cE(x)))$ and $D_2 \in Ret(\cE(x))$. The retract diagram for maps $f: A \to B$ and $g: X \to Y$ gives us a morphism $D_1 \to Ret(c^*_I)(D_2)$. Taking the adjoint to this map, we get a map of retract diagrams $Ret(\colim_I)(D_1) \to D_2$, which renders the relative latching map of $f$ at $x$, $$\Lat_x B \coprod_{\Lat_x A} A(x) \to B(x),$$ as a retract of the relative latching map of $g$ at $x$, $$\Lat_x Y \coprod_{\Lat_x X} X(x) \to Y(x).$$ Thus if the latter map is a (trivial) cofibration, then so is the former. For the relative matching maps, the proof is dual.
\endproof

\subsubsection*{Case of a direct category}
\label{subsectiondirectcategory}

We first consider the case when $\cR = \cR_+$ is a direct Reedy category. In this case $\cE \to \cR$ is an actual opfibration. Similarly, one can consider $\cR = \cR_-$, and work with a fibration over $\cR$.

\begin{prop}
	Reedy cofibrations, objectwise fibrations and objectwise weak equivalences form a model structure on $\Sect(\cR,\cE)$.
\end{prop}

First we need to address the limit-colimit axiom.
\begin{lemma}
\label{lemmadirectcategorysectionsbicomplete}
For $\cR = \cR_+$, the category $\Sect(\cR,\cE)$ admits limits and colimits.
\end{lemma}


\proof The existence of limits is Proposition \ref{prefibrationcolimits}.
The colimits are given by the dual of Proposition \ref{noetherlimits} since $\cR_+$ is
an Artin category. \endproof
\begin{lemma}
	Suppose given a diagram of sections
	\begin{diagram}[small,nohug]
	A	& \rTo	& S	 \\
	\dTo<f	 &			&	\dTo>p		\\
	B & \rTo & T	 \\
	\end{diagram}
	with $p$ and objectwise fibration (respectively trivial fibration). If for each $x \in \cR$, the map
	\begin{equation}
	\label{reedycof1}
	\Lat_x B  \coprod_{\Lat_x A } A(x) \to B(x)
	\end{equation}
	is a trivial cofibration (respectively a cofibration), then the diagram admits a lift.
\end{lemma}
\proof Proceed by induction on degree. For $deg x = 0$, $\Lat_x(A)$ is the initial object of $\cE(x)$ (and the same for $B$), so the map (\ref{reedycof1}) equals $A(x) \to B(x)$. The lift then exists simply because $\cE(x)$ is a model category.

For $deg x = n$, assume that we defined the lift for all lesser degrees. For each map $\alpha: y \to x$ with $deg y < n$, we have the assumed lift $h_y: B(y) \to S(y)$, and the composition $B(y) \to S(y) \to S(x)$ can be factored as $\alpha_! B(y) \to S(x)$, and that in turn induces the map $\Lat_x B \to S(x)$. We then get the following diagram,
\begin{diagram}[size=2.5em,nohug]
\Lat_x B  \coprod_{\Lat_x A } A(x)	& \rTo	& S(x)	 \\
\dTo<f	 &			&	\dTo>p		\\
B(x) & \rTo & T(x)	 \\
\end{diagram}
and we can find the necessary lift (by also remembering $A(x) \to \Lat_x B  \coprod_{\Lat_x A } A(x)	$). \endproof

\begin{lemma}
	\label{reedycofibrationsgivepointwise}
	Let $\alpha: A \to B$ be such that $ \Lat_x B  \coprod_{\Lat_x A } A(x) \to B(x)$ is a (trivial) cofibration for each $x \in \cR$. Then for any $y \in \cR$, the maps $\Lat_y(\alpha) : \Lat_y A \to \Lat_y B$ and $\alpha_y: A(y) \to B(y)$ are (trivial) cofibrations.
\end{lemma}

\proof For $y \in \cR$, note that the map $\alpha_y: A(y) \to B(y)$ equals
$$
A(y) \to \Lat_y B  \coprod_{\Lat_y A } A(y) \to B(y).
$$
The second map is a (trivial) cofibration by condition. It thus remains to examine the map $\Lat_y(\alpha): \Lat_y A \to \Lat_y B$. According to Definition \ref{admissibilitydef},
this map is a (trivial) cofibration, as required. \endproof

\begin{corr}
	Let $A \to B$ be such that $\Lat_x B  \coprod_{\Lat_x A } A(x) \to B(x)$ is a trivial cofibration for each $x \in \cR$. Then $A \to B$ is a Reedy cofibration and a weak equivalence. \endproof
\end{corr}

\begin{prop}
	\label{directfactor}
	Let $A \to C$ be a map in $\Sect(\cR,\cE)$. Then it can be factored as $A \to B \to C$ where
	\begin{itemize}
		\item the map $A \to B$ is such that $\Lat_x B  \coprod_{\Lat_x A } A(x) \to B(x)$ is a cofibration (respectively a trivial cofibration) for each $x \in \cR$,
		\item the map $B \to C$ is an objectwise trivial fibration (respectively a fibration).
	\end{itemize}
	The factorisations are functorial whenever this is the case for each $\cE(x)$.
\end{prop}
\proof Let us do the cofibration and trivial fibration part, the second part being dual. Factor $A(x) \to B(x)$ as $A(x) \to B(x) \to C(x)$ for each $x$ of degree zero. Assume now that the factorisation is there for each $y \in \cR$ of degree less than $n$. For $x$ with $deg x = n$, we have the diagram
\begin{diagram}[small]
\Lat_x A & \rTo & \Lat_x B \\
\dTo		&			&		\dTo		\\
A(x) &	\rTo		&	C(x)		\\
\end{diagram}
with $\Lat_x B \to C(x)$ defined with the use of the maps $B(y) \to C(y) \to C(x)$. We thus get a map $  A(x) \coprod_{\Lat_x A} \Lat_x B \to C(x)$, which we factor (if possible, functorially) as
$$
A(x) \coprod_{\Lat_x A} \Lat_x B \to B(x) \to C(x).
$$
The maps $\Lat_x B \to B(x)$ complete $B$ to a section on $\cR_{\leq n}$. Proceeding by induction, we get the desired factorisation. \endproof

\begin{corr}
	\label{trivialreedycof}
	A map $f:S \to T$ is a trivial Reedy cofibration iff the map $$\displaystyle \Lat_x T  \coprod_{\Lat_x S } S(x) \to T(x)$$ is a trivial cofibration for each $x \in \cR$.
\end{corr}
\proof
Take a trivial Reedy cofibration $f: S \to T$ and factor it using Proposition \ref{directfactor} as $S \stackrel g \to U \stackrel h \to T$ so that $  \Lat_x U  \coprod_{\Lat_x S } S(x) \to U(x)$ is a trivial cofibration. We then see that $f$ is a retract of $g$. \endproof

This proves the existence of the model structure on $\Sect(\cR,\cE)$ for a direct category $\cR$.

\subsubsection*{Finishing the Proof}

We now turn to the case when $\cR$ is an arbitrary Reedy category.

\begin{lemma}
The category $\Sect(\cR,\cE)$ is bicomplete.
\end{lemma}
\proof By Lemma \ref{lemmadirectcategorysectionsbicomplete} and its dual, we have that
both $\Sect(\cR_+,\cE)$ and $\Sect(\cR_-,\cE)$ are bicomplete. The result then follows
from Proposition \ref{janussectionslimits}. \endproof

\begin{lemma}
	\label{lemmatrivreedycofifffibrewise}
	A map $X \to Y$ is
	\begin{itemize}
		\item a trivial Reedy cofibration iff for each $x \in \cR$, the map $  \Lat_x Y  \coprod_{\Lat_x X } X(x) \to Y(x)$ is a trivial cofibration,
		\item a trivial Reedy fibration iff for each $x \in \cR$, the map $  X(x) \to Y(x) \prod_{\Mat_x Y } \Mat_x X$ is a trivial fibration.
	\end{itemize}
\end{lemma}
\proof For the first part, note that $X \to Y$ is a Reedy cofibration iff it is such when viewed as a morphism of sections in $\Sect(\cR_+,\cE)$, since the Reedy cofibration condition is formulated objectwise in $\cR$. It is, also, a weak equivalence iff it is such when restricted to a morphism of sections over $\cR_+$, for the same reason. We then use Corollary \ref{trivialreedycof} to get the result. The second part is proven in a dual manner. \endproof

\begin{prop}
	Suppose given a diagram of sections
	\begin{diagram}[small,nohug]
	A	& \rTo	& S	 \\
	\dTo<f	 &			&	\dTo>p		\\
	B & \rTo & T	 \\
	\end{diagram}
	where $f: A \to B$ is a Reedy cofibration and $p: S \to T$ is a Reedy fibration. Then a lift exists whenever $f$ or $p$ is trivial.
\end{prop}
\proof By induction we can assume having supplied a lift for $y \in \cR$ of degree less than $n$. Given an object $x$ of degree $n$, we can draw the following diagram
\begin{diagram}[size=2.5em,nohug]
A(x) & \rTo 			& A(x) \coprod_{\Lat_x A} \Lat_x B	& \rTo	& S(x)	&		&	 \\
&	\rdTo		&	\dTo	 &			&	\dTo											&	\rdTo	&\\
&					&	B(x) & \rTo & T(x) \prod_{\Mat_x T} \Mat_x S	& \rTo	& T(x). \\
\end{diagram}
Just as in the classical case, a lift in the middle square of this diagram (which exists whenever $f$ or $p$ is trivial) determines the looked-for lift $B \to S$ on objects of degree $n$.\endproof

\begin{prop}
	\label{reedyfactorisation}
	Let $A \to C$ be a map in $\Sect(\cR,\cE)$. Then it can be factored as $A \stackrel i \to B \stackrel p \to C$, with $i$ a Reedy cofibration and $p$ a Reedy fibration, such that either $i$ or $p$ is trivial. The factorisation is functorial whenever each $\cE(x)$ admits functorial factorisations.
\end{prop}
\proof
Assume again that, by induction, we have constructed the factorisation $A(y) \to B(y) \to C(y)$ for objects $y \in \cR$ of degree less than $n$. For $x$ of degree $n$, there is the following diagram
\begin{diagram}[small]
\Lat_x A & \rTo & A(x)		&	\rTo & \Mat_x B  \\
\dTo		&			&	\dTo		&			&		\dTo			\\
\Lat_x B &	 \rTo & C(x)		&	\rTo	&	\Mat_x C	\\
\end{diagram}
which exists because of the inductive assumption and provides us with the following map
$$
\Lat_x B \coprod_{\Lat_x A} A(x) \to C(x) \prod_{\Mat_x C} \Mat_x B.
$$
Factoring it (using the model structure of $\cE(x)$) as
$$
\Lat_x B \coprod_{\Lat_x A} A(x) \to B(x) \to C(x) \prod_{\Mat_x C} \Mat_x B.
$$
which, together with maps $\Lat_x B \to B(x)$ and $B(x) \to \Mat_x B$, yields the desired extension of the factorisation to the objects of degree $n$.\endproof

We have thus proven the existence of the Reedy model structure on $\Sect(\cR,\cE)$.

\begin{lemma}
	\label{reedycoffibarepointwisecoffib}
	Let $X \to Y$ be a Reedy cofibration (respectively a fibration). Then for each $x \in \cR$, the map $X(x) \to Y(x)$ is a cofibration (respectively a fibration).
\end{lemma}
\proof Direct consequence of Lemma \ref{reedycofibrationsgivepointwise}. \endproof

\subsection{Cofibrant generation of Reedy structures}

In this subsection we treat the situation when the semifibration $\cE \to \cR$
has cofibrantly generated model categories as fibres. One may ask in this case if the
model category of sections is cofibrantly generated.

The Reedy category case is unorthodox in the sense that we already know the model
structure, while the techniques of cofibrantly generated categories usually start
with a set of generating cofibrations and a well-behaved class (large set) of weak equivalences
in a presentable category to obtain a combinatorial model structure, as per Smith's theorem \cite[Proposition A.2.6.8]{LUHTT}.

For the purposes of this subsection, let us introduce the following definition.
As usual, write $[1]$ for the arrow category $0 \to 1$ and $[2]=0 \to 1 \to 2$.

\begin{opr}
A model category $\cM$ has \emph{accessible factorisations}, if the underlying
category of $\cM$ is presentable, and both (fibration-trivial cofibration) and
(trivial fibration-cofibration) factorisations are functorial and accessible,
when viewed as functors from $\Fun([1],\cM)$ to $\Fun([2],\cM)$.
\end{opr}

Any combinatorial model category has accessible factorisations \cite[Proposition 1.10]{BARLEFT}.
We would like to investigate when the converse implication holds. For our purposes, the
following shall suffice.

\begin{lemma}
\label{lemmareversecofgen}
Let $\cM$ be a model category with accessible factorisations, such that the
subcategories of trivial fibrations and fibrations are accessible and accessibly
embedded in $\Fun([1],\cM)$. Then $\cM$ is cofibrantly generated, and hence is combinatorial.
\end{lemma}
\proof Choose a regular cardinal $\lambda$ such that the following holds:
\begin{enumerate}
\item The categories $\cM, \Fun([1],\cM)$ and the subcategory $\cF ib \subset \Fun([1],\cM)$
consisting of fibrations (and a similar category for trivial fibrations) is $\lambda$-accessible,
and the functor $\cF ib \to \Fun([1],\cM)$ preserves $\lambda$-presented objects.

\item For each arrow $X \to Y$ where $X$ and $Y$ are $\lambda$-presentable
(or $\lambda$-compact in the terminology of \cite{LUHTT}), its factorisations
\begin{diagram}[small]
X & \rInto^\sim & Z & \rOnto Y & \text{and} & X & \rInto & Z' & \rOnto^\sim Y \\
\end{diagram}
have the property that both $Z$ and $Z'$ are $\lambda$-presentable. This is possible
due to the accessible factorisations condition, following the same reasoning as in
\cite[Proposition 7.2]{DUGGER}.
\end{enumerate}

The category $\Fun([1],\cM)$ is seen \cite[Proposition 5.4.4.3]{LUHTT} to be
generated by the (essentially small) subcategory of arrows $A \to B$ where $A$ and $B$ are $\lambda$-presentable.
We claim that the trivial cofibrations are generated by the subset $\cW \cC of_\lambda$ consisting of
all trivial cofibrations between $\lambda$-presentable objects.

It will suffice to check
that any map $f:X \to Y$ having the right lifting property with respect to $\cW \cC of_\lambda$
is a fibration. As an object of $\cM^{[1]}= \Fun([1],\cM)$, the morphism $f$ is a
colimit of $\lambda$-presentable objects over a small (modulo choice), $\lambda$-filtered diagram
$\cM^{[1]}_{\lambda}/f$. There is also a diagram $\cF ib_\lambda / f$ consisting of
all $\lambda$-presentable objects of $\cF ib$ together with a map in $\cM^{[1]}$ to $f$.
Since $\cF ib$ is accessibly embedded we have a natural fully faithful functor
$$
F:\cF ib_\lambda / f \to \cM^{[1]}_{\lambda}/f
$$
We claim that $F$ is cofinal. For this, just as in \cite[Theorem 4.8]{ADAMEK}, it
is enough to verify that for any $g \to f$ in $\cM^{[1]}_{\lambda}/f$ there exists
a map $g \to F(h)$ over $f$ for some $h \in \cF ib_\lambda / f$.

Factoring $g: A \to B$ as a trivial cofibration followed by fibration, we get the
diagram
\begin{diagram}[small]
A & \rInto^\sim & A' \\
\dTo<g	&				& \dOnto>{F(h)} \\
B & \rTo_= & B \\
\end{diagram}
the object $A'$ is $\lambda$-presentable just like $A$ and $B$, so the fibration
$F(h): A' \twoheadrightarrow B$ is $\lambda$-presentable in $\cF ib$ as well. Moreover any map
$g \to f$ can be extended in a compatible way to $F(h) \to f$, using the assumption
 that $f$ has the right
lifting property along trivial cofibrations between $\lambda$-presentable objects.
This concludes the proof of cofinality of $F$, and hence also of $\lambda$-filteredness
of $\cF ib_\lambda / f $. We then use the fact that $\cF ib$ is closed under $\lambda$-
filtered colimits to conclude that $f \in \cF ib$. The case of trivial fibrations is
treated similarly.
\endproof

\begin{corr}
Let $\cM$ be a model category with underlying category presentable. Assume that
both of the following hold:
\begin{itemize}
\item[i.] The (cofibration-trivial fibration) factorisation is functorial and accessible,
\item[ii.] The full subcategory $\cW \cF ib \subset \cM^{[1]}$ of trivial fibrations is accessible and accessibly
embedded.
\end{itemize}
Assume further that either of the following holds:
\begin{enumerate}
\item The (trivial cofibration- fibration) factorisation is functorial and accessible, and the full subcategory $\cF ib \subset \cM^{[1]}$ of fibrations is accessible and accessibly
embedded, or
\item The full subcategory $\cW \subset \cM^{[1]}$ of weak equivalences is accessible and accessibly
embedded.
\end{enumerate}
Then $\cM$ is combinatorial.
\end{corr}
\proof Combine the precedent lemma together with \cite[Corollary A.2.6.9]{LUHTT} \endproof

\begin{prop}
\label{propreedycofibrantgenerationimplicit}
Let $\cE \to \cR$ be an admissible model semifibration over a Reedy category $\cR$.
Assume that each $\cE(x)$ is combinatorial and all transition functors of the semifibration
are accessible. Then the model category of sections $\Sect(\cR,\cE)$ is combinatorial.
\end{prop}
\proof
The category $\cR$ is a directed colimit of its subcategories $\cR_{\leq n}$. Let us
first analyse the case of finite degree.

For $\cR_0$, the category $\Sect(\cR_0, \cE)$ is the product $\prod_x \cE(x)$
of the fibre categories over all objects of degree zero. It is, hence, combinatorial.

For $\cR_{\leq n}$, observe that we have the following pullback diagram
\begin{equation}
\label{diagramcofibrantgeneration}
\begin{diagram}[size=2.5em]
\Sect(\cR_{\leq n},\cE) & \rTo & \prod_{deg x = n} \cE(x)^{[2]} \\
\dTo									&					&							\dTo				\\
\Sect(\cR_{< n},\cE) & 		\rTo			& \prod_{deg x = n} \cE(x)^{[1]}. \\
\end{diagram}
\end{equation}
The upper horizontal functor is given by $S \mapsto \Lat_x S \to S(x) \to \Mat_x S$
for all objects $x$ of degree $n$. The bottom horizontal functor is similarly given
by $S \mapsto \Lat_x S  \to \Mat_x S$. The vertical functors are the evident restrictions.

Both categories $\prod_x \cE(x)^{[2]}$ and $\prod_x \cE(x)^{[1]}$ are presentable.
By induction, we can assume that $\Sect(\cR_{< n},\cE)$ are combinatorial and that
both $\Lat_x$ and $\Mat_x$ are accessible functors (the initialization is given by $\emptyset = \cR_{-1} \subset \cR_0$). This implies that
the horizontal functors in (\ref{diagramcofibrantgeneration}) are accessible.
We thus see that $\Sect(\cR_{\leq n},\cE)$ is accessible,
and hence presentable. We also infer, by induction, that sufficiently large
filtered colimits are computed fibrewise in $\Sect(\cR_{\leq n},\cE)$.

We shall not treat $\cE(x)^{[1]}$ and $\cE(x)^{[2]}$ as model category. However, consider the
full subcategory $ \cF ib(x)^{[2]} \subset \Fun([1],\cE(x)^{[2]})$ given by
all diagrams
\begin{diagram}[small]
A & \rTo & B & \rTo & C \\
\dTo &		& \dTo &  	& \dTo \\
X & \rTo & Y & \rTo & Z \\
\end{diagram}
such that both $C \to Z$ and $B \to Y \times_Z C$ are fibrations. The category
$ \cF ib(x)^{[2]}$ is the pullback, along the inclusion $[1] \subset [2]$ skipping $0$,
of the subcategory of fibrations for the inverse Reedy model structure on each $\cE(x)^{[1]}$. The
latter is well known to be combinatorial (it follows
for example from Proposition \ref{propreedycofibrantgenerationexplicit} below). Using
\cite[Corollary A.2.6.9]{LUHTT} and the closure of accessible categories under limits,
we obtain that $ \cF ib(x)^{[2]}$ is accessible and accessibly embedded.
In a similar vein, we can define an accessible and accessibly embedded subcategory
$\cF ib(x)^{[1]} \subset \Fun([1],\cE(x)^{[1]})$
to be given by all diagrams
\begin{diagram}[small]
A & \rTo   & C \\
\dTo &		&   \dTo \\
X   & \rTo & Z \\
\end{diagram}
such that $C \to Z$ is a fibration.

We then conclude, using the diagram (\ref{diagramcofibrantgeneration}), that the category of
fibrations of $\Sect(\cR_{\leq n}, \cE)$ (and, similarly, of trivial fibrations)
is accessible and accessibly embedded, since we can express it as a limit of
$\prod_x \cF ib(x)^{[2]}$, $\prod_x \cF ib(x)^{[1]}$ and the fibrations of
$\Sect(\cR_{< n},\cE)$ along accessible functors. This means that the fibrations
(and, similarly, trivial fibrations) form an accessible (and accessibly embedded)
subcategory of $\Sect(\cR_{\leq n},\cE)^{[1]}$. We also see that the factorisations
in $\Sect(\cR_{\leq n},\cE)$ are functorial and accessible: they are obtained through
the inductive procedure involving lesser degree latching and matching objects,
fibred products and pushouts, and factorisations in the fibres.

Now, given an object $y$ of the successive degree, consider the matching object functor
$\Mat_y: \Sect(\cR_{\leq n},\cE) \to \cE(y)$. It is obtained as the composition of
the restriction to the fibre $$\Sect(\cR_{\leq n},\cE) \to \Sect(Mat(y),\cE)  \to \Fun(Mat(y),\cE(y))$$
with the $Mat(y)$-limit functor.
Bearing in mind that we proved that (sufficiently large) filtered colimits
are computed fibrewise in $\Sect(\cR_{\leq n},\cE)$, all of these functors are accessible.
We thus obtain that $\Mat_y$ is accessible. One can prove similar for $\Lat_y$.

All this shows that for each degree $n$, the category $\Sect(\cR_{\leq n},\cE)$ is
a combinatorial model category. The restriction functors $\Sect(\cR_{\leq n},\cE)
\to \Sect(\cR_{\leq n-1},\cE)$ commute with colimits. Using the argumentation
as before, but with pullbacks of categories replaced by directed colimits,
we obtain that $\Sect(\cR,\cE) = \lim \Sect(\cR_{\leq n},\cE)$ satisfies the conditions
of Lemma \ref{lemmareversecofgen}. It is hence combinatorial. \endproof

We conclude this subsection by considering the case when $\cE \to \cR$ is a bifibration, that is, both
a Grothendieck fibration and an opfibration. The assumptions of the model semifibration (which in this
case is automatically admissible due to (1.) of Lemma \ref{admissibilitycriterion})
imply that for any $f: x \to y$, the induced adjunction
\begin{equation}
\label{eqqpresheafadjunction}
f_!: \cE(x) \rightleftarrows \cE(y):f^*
\end{equation}
is a Quillen pair. In the terminology of \cite{H-S}, we are dealing with a \emph{Quillen presheaf}.

Let $x \in \cR$ and $\cC$ be a subcategory of the comma category $x \backslash \cR$.
It comes equipped with a natural functor $p_\cC:\cC \to \cR$. This functor
induces the following adjunction:
$$
p_{\cC,!}: \Sect(\cC,\cE) \rightleftarrows \Sect(\cR,\cE): p_\cC^*.
$$
The left adjoint exists in the case of a fibrewise cocomplete bifibration
and sends $S$ to the section determined by
$$
p_{\cC,!}S(y) = \colim_{\cC/y} Res_y S|_{\cC/y}
$$
with $Res_y: \cE|_{\cC/y} \to \cE(y)$ being the usual restriction functor. On the other hand,
there is also an adjunction
$$
triv_\cC: \cE(x) \rightleftarrows \Sect(\cC,\cE): \Mat_\cC.
$$
The functor $triv_\cC$ sends $X \in \cE(x)$ to the section
$$
(f:x \to y) \in \cC \, \, \, \mapsto \, \, \, f_! X \in \cE(y).
$$
The functor $\Mat_\cC$ is the composition $\Sect(\cC,\cE) \to \Fun(\cC,\cE(x)) \to \cE(x)$
where the first functor is induced by restrictions along cartesian morphisms, and
the second is the limit functor. We thus have the composed adjunction
$$
p_{\cC,!}triv_\cC: \cE(x) \rightleftarrows \Sect(\cR,\cE): \Mat_\cC p_\cC^*.
$$
\begin{lemma}
For a Quillen presheaf $\cE \to \cR$, we have the following:
\begin{enumerate}
\item For each object $X \in \cE(x)$, the section
$$
i(X) = p_{x \backslash \cR, !} triv_{ x \backslash \cR} X
$$
satisfies the property
$$
\Sect(\cR,\cE)(i(X),S) \cong \cE(x)(X,S(x)).
$$
\item For each object $X \in \cE(x)$, the section
$$
m(X) = p_{Mat(x), !} triv_{ Mat(x)} X
$$
satisfies the property
$$
\Sect(\cR,\cE)(m(X),S) \cong \cE(x)(X,\Mat_x S).
$$
\item For each object $X \in \cE(x)$, we have a natural map $m(X) \to i(X)$.
\end{enumerate}
\end{lemma}
\proof Immediate. \endproof

\begin{prop}
\label{propreedycofibrantgenerationexplicit}
Let $\cE \to \cR$ be a Quillen presheaf. Assume that each model category $\cE(x)$
is cofibrantly generated, with generating cofibrations denoted $I_x$ and generating
trivial cofibrations denoted $J_x$. Write
$$I = \{ m(B) \coprod_{m(A)} i(A) \to i(B) | A \to B \in I_x, x \in \cR \},$$
$$J = \{ m(B) \coprod_{m(A)} i(A) \to i(B) | A \to B \in J_x, x \in \cR \}.$$
Then $I$ and $J$ are sets of generating cofibrations and trivial cofibrations for
the model structure of Theorem \ref{reedymodelstructuretheorem} on $\Sect(\cR,\cE)$.
\end{prop}
\proof
The usual adjunction observation, amounting to the fact that a diagram
\begin{diagram}[small]
m(B) \coprod_{m(A)} i(A) & \rTo & S \\
\dTo 								&			&		\dTo \\
i(B)								& \rTo & T \\
\end{diagram}
is the same data as the diagram
\begin{diagram}[small]
A     &	\rTo & S(x) \\
\dTo	&				& \dTo \\
B			& \rTo & T(x) \prod_{\Mat_x T} \Mat_x S \\
\end{diagram}
\endproof

%% file: rmsf-comparison.tex
\section{Comparing with higher sections}
\label{subsectionreedycomp}

\begin{opr}
A \emph{left model Reedy fibration} is a functor $p:\cE\to \cR$ to a Reedy category
that is a Grothendieck opfibration, a Grothendieck fibration over $\cR_-$, such that
the associated semifibration (cf Lemma \ref{semifibexm}) is admissible in the sense
of Definition \ref{admissibilitydef}.
\end{opr}

Recall the notion of an $(\infty,1)$-category, which can be modeled using quasicategories (called ``infinity-categories'' in \cite{LUHTT}).
Denote by $\cW \subset \cE$ the collection of those maps $\alpha: X \to Y$
such that $p\alpha$ is an isomorphism and the induced map $(p \alpha)_! X \to Y$ is
a fibrewise weak equivalence. Using the localisation in quasicategories \cite{HINICH},
localising along $\cW \subset \cE$ yields an infinity-functor $Lp: L \cE \to \cR$ (we make no distinction
between an ordinary category and its nerve in $\SSet$), which can be chosen to be
a categorical fibration. The associated quasicategory
of sections $\Sect(\cR,L\cE)$ is given by the sub-quasicategory of
$\Fun(\cR,L\cE)$ spanned by those $S$ such that $Lp \circ S$ is an identity; more
precisely $\Sect(\cR,L\cE)$ is the fibre of $Lp_*: \Fun(\cR,L \cE) \to \Fun (\cR, \cR)$
over the identity functor.

The goal of this subsection is to prove the strictification result, asserting
that the infinity-localisation $L \Sect (\cR,\cE)$
of the category of sections coincides with $\Sect(\cR,L\cE)$. After proving this,
we shall revisit the case of Quillen presheaves --- those $\cE \to \cC$ that are
bifibrations in model categories and Quillen adjunctions --- and show that
the strictification holds over an arbitrary base $\cC$.

\begin{rem}
The terminology employed in \cite{LUHTT} is different from the one which we have used
up to this point. Lurie uses the term ``coCartesian fibration'' \cite[Definition 2.4.2.1]{LUHTT} for what we would have
called ``opfibration of quasicategories'', and ``Cartesian fibration'' for fibrations.
We have chosen to stick to the terminology of \cite{LUHTT} when dealing with higher-categorical
fibrations, and we use our terminology for the 1-categorical objects.

We also specify that for us, the term infinity-localisation means a
functor $F: \cX \to \cY$ of quasicategories, such that for any quasicategory $\cZ$,
the induced infinity-functor $F^*:\Fun(\cY,\cZ) \to \Fun(\cX,\cZ)$ is full and faithful,
and its essential image consists of functors that send the $F$-equivalences of $\cX$ to equivalences
of $\cZ$. This is not the same meaning as the localisation used in the setting
of presentable infinity-categories \cite[Chapter 5]{LUHTT}.
\end{rem}

For any Reedy category, one can introduce the notion of a (transfinite) good filtration $\{ \cR_\beta \}$
\cite[Notation A.2.9.11]{LUHTT} for a Reedy category $\cR$ (a Reedy category in
our sense is a Reedy category in the sense of Lurie). The advantage of a good filtration
is that, for an ordinal $\beta$, one obtains $\cR_\beta$ by adjoining a single object
to $\cR_{< \beta} = \cup_{\gamma < \beta} \cR_\gamma$.

The key observation \cite[Proposition A.2.9.14]{LUHTT} for Reedy categories in
higher-categorical setting asserts that the diagram
\begin{diagram}[small]
(\cR_{<\beta}/x) \star (x \backslash \cR_{< \beta}) & \rTo & \cR_{< \beta} \\
\dTo																								&			&		\dTo		\\
(\cR_{<\beta}/x) \star \{x \} \star (x \backslash \cR_{< \beta}) & \rTo & \cR_\beta \\
\end{diagram}
is a homotopy pushout square of quasicategories. We modify this statement slightly.

\begin{prop}
\label{propreedyfiltrationstep}
Let $\cR$ be a Reedy category, and $\cR_\beta$ be a step of a good filtration obtained
from $\cR_{< \beta}$ by adding an object $x$. Then the following square
\begin{diagram}[small]
Lat(x) \star Mat(x) & \rTo & \cR_{< \beta} \\
\dTo									&			&		\dTo		\\
Lat(x) \star \{x \} \star Mat(x) & \rTo & \cR_\beta \\
\end{diagram}
is a homotopy pushout for the Joyal model structure.
\end{prop}
\proof
Observe that the inclusion $Lat(x) \subset \cR_{< \beta} / x$ is a right adjoint, with the
left adjoint obtained by using the factorisation system. This inclusion is
hence homotopy cofinal, as follows immediately from \cite[Theorem 4.1.3.1]{LUHTT}. Propositions 4.1.2.5 and 4.1.2.1 of \cite{LUHTT}
imply that $Lat(x) \to \subset \cR / x$ is right anodyne. Taking the pushout-join
with $\emptyset \subset \{x \}$ and applying \cite[Lemma 2.1.2.3]{LUHTT} yields that the diagram
\begin{diagram}[small]
Lat(x)   & \rTo &  (\cR_{< \beta} / x) \\
\dTo			&			&		\dTo		\\
Lat(x) \star \{x \} & \rTo & (\cR_{< \beta} /x) \star \{x \} \\
\end{diagram}
is a homotopy pushout for the Joyal model structure. Using the fact that joins
preserve connected homotopy colimits (see e.g. \cite[Lemma 4.14]{SHAH} for a relative version)
and a dual argument for the inclusion $Mat(x) \subset x \backslash \cR_{< \beta}$, we
conclude that we have a series of homotopy pushout squares
\begin{diagram}[small]
Lat(x) \star Mat(x)   & \rTo &  (\cR_{< \beta} / x ) \star Mat(x)  & \rTo & (\cR_{<\beta}/x) \star (x \backslash \cR_{< \beta})  \\
\dTo			&			&		\dTo	&  & \dTo 	\\
Lat(x) \star \{x \}  \star Mat(x) & \rTo & (\cR_{< \beta} /x) \star \{x \}  \star Mat(x) & \rTo & (\cR_{<\beta}/x) \star \{x \} \star (x \backslash \cR_{< \beta})\\
\end{diagram}
which together with \cite[Proposition A.2.9.14]{LUHTT} and pasting for homotopy pullbacks implies the result.
\endproof

It will be easiest to state the results in the infinity-category of all infinity-categories $\Cat_\infty$, with usual size issue remarks. The reason for this is that many of the proofs in this section will manipulate with diagrams which are canonically presented in $\SSet$, but in which some categorical equivalences
are pointing in the wrong direction. Such zig-zags of diagrams will induce well-defined
diagrams in $\Cat_\infty$:
\begin{lemma}
\label{lemmazigzagsofdiagramsincatinfty}
Assume given a zig-zag of diagrams $I \times [1] \to \SSet$,
valued in quasicategories and depicted schematically as
\begin{diagram}[small]
X_\bullet(0) 	& \lTo^\sim	& ...	&	\rTo & Y_\bullet(0)	\\
\dTo 					&						&	...	&      & \dTo					\\
X_\bullet(1) 	& \lTo^\sim	&	...	& \rTo & Y_\bullet(1)	\\
\end{diagram}
with the bottom index corresponding to the $I$-direction. Assume that each left-pointing arrow is a categorical equivalence. Then there exists an induced commutative diagram in $\Cat_\infty$,
\begin{diagram}[small]
X_\bullet(0) 	&  \rTo & Y_\bullet(0)	\\
\dTo 					&	 			& \dTo					\\
X_\bullet(1) 	&	\rTo & Y_\bullet(1),	\\
\end{diagram}
with each horizontal arrow induced by choosing inverses of left-pointing arrows and
compositions.
\end{lemma}
\proof
Direct consequence of \cite[Corollary 3.5.12]{CISBOOK} applied to $\Cat_\infty$ itself.
\endproof

We will often identify an $\SSet$-diagram with its image
in $\Cat_\infty$ if the details are clear from the context.

The following proposition is the expression of the Reedy induction as applied to
a suitably bicomplete infinity-category.

\begin{prop}
\label{proplurieremarkproven}
(Cf \cite[Remark A.2.9.16]{LUHTT}) let $\cR$ be a Reedy category and $\cY$
be a quasicategory. Assume that $\cY$ admits $Lat(y)$-colimits and $Mat(y)$-limits for
each $y \in \cR$. Then, in the usual notation for a good filtration, there is a diagram in $\Cat_\infty$
\begin{equation}
\label{dialurieremarkproven}
\begin{diagram}[small]
\Fun(\cR_\beta ,\cY) & \rTo & \Fun( \cR_{< \beta},\cY) \\
\dTo																													&						&				\dTo					\\
\Fun(Lat(x) \star \{x\} \star Mat(x),\cY) & \rTo & \Fun(Lat(x) \star Mat(x),\cY) \\
\dTo																													&						&				\dTo					\\
\Fun([2],\cY) 	&	\rTo	&	\Fun([1], \cY) \\
\end{diagram}
\end{equation}
with all squares pullbacks. The upper vertical functors are restrictions, the
bottom left vertical functor sends $S$ to $\colim_{Lat(x)} S \to S(x) \to \lim_{Mat(x)} S$,
and the bottom right vertical functor sends $S'$ to $\colim_{Lat(x)} S' \to \lim_{Mat(x)} S'$.
\end{prop}
\proof The nontrivial square of the diagram (\ref{dialurieremarkproven}) is the bottom one.
Let us correctly construct the functors appearing in it.
Consider the inclusion $Lat(x) \star \{x\} \star Mat(x) \subset Lat(x)^\triangleright
 \star \{x\} \star Mat(x)^\triangleleft$. The induced restriction functor
$$
\Fun(Lat(x) \star \{x\} \star Mat(x), \cY) \longleftarrow \Fun(Lat(x)^\triangleright
 \star \{x\} \star Mat(x)^\triangleleft, \cY)
$$
admits a section, given by putting a colimit in the cone vertex of $Lat(x)^\triangleright$
and a limit in the cone vertex of $Mat(x)^\triangleleft$. Formally, it is a composition
of a left then right Kan extension along the full inclusions
$$
Lat(x) \star \{x\} \star Mat(x) \subset Lat(x)^\triangleright
 \star \{x\} \star Mat(x) \subset Lat(x)^\triangleright
  \star \{x\} \star Mat(x)^\triangleleft.
$$
The results of \cite[4.3.2]{LUHTT} then imply that the induced functor
$$
\Fun(Lat(x) \star \{x\} \star Mat(x), \cY) \longleftarrow \Fun'(Lat(x)^\triangleright
 \star \{x\} \star Mat(x)^\triangleleft, \cY)
$$
is an equivalence of infinity-categories, where $\Fun'(Lat(x)^\triangleright
 \star \{x\} \star Mat(x)^\triangleleft, \cY)$ is the full subcategory consisting of
 functors $F:Lat(x)^\triangleright
  \star \{x\} \star Mat(x)^\triangleleft \to \cY$
which carry the cone vertex $l \in Lat(x)^\triangleright$ to the colimit of $F$
over $Lat(x)$, and the cone vertex $m \in Mat(x)^\triangleleft$ to the limit of
$F$ over $Mat(x)$. One can apply exactly the same analysis to the inclusion
$Lat(x) \star Mat(x) \subset Lat(x)^\triangleright \star Mat(x)^\triangleleft$.

The following will suffice. For any two categories $A,B$,
we need to show that the diagram
\begin{diagram}[small]
\Fun(A^\triangleright \star \{x\} \star B^\triangleleft, \cY) & \rTo & \Fun([2],\cY)	\\
\dTo							&																							&	\dTo				\\
\Fun(A^\triangleright \star B^\triangleleft, \cY) & \rTo & \Fun([1],\cY) \\
\end{diagram}
is a pullback square in $\Cat_\infty$. If we denote by $a \in A^\triangleright$
and $b \in B^\triangleleft$ the cone vertices, then the horizontal arrows are
induced by $[2] \cong \{a \} \star \{x \} \star \{b \} \subset A^\triangleright \star \{x\} \star B^\triangleleft$
and $[1] \cong \{a \} \star \{b \} \subset A^\triangleright \star B^\triangleleft$
respectively, with vertical arrows then also being the obvious restrictions.

In turn, it will suffice to prove that the diagram
\begin{equation}
\label{diagrampushoutinsertingmiddlex}
\begin{diagram}[small]
\{a\} \star \{b \} & \rTo & A^\triangleright \star B^\triangleleft \\
\dTo	& 	& \dTo \\
\{a\} \star \{ x \} \star \{b \} & \rTo & A^\triangleright \star \{x\} \star B^\triangleleft \\
\end{diagram}
\end{equation}
is a pushout in $\Cat_\infty$. To see that the latter is true, assume first that
$A = \emptyset$. Then the diagram (\ref{diagrampushoutinsertingmiddlex}) is a join
with $\{a \} \star -$ of
\begin{equation}
\label{diagrampushoutinsertingmiddlexaempty}
\begin{diagram}[small]
\{b \} & \rTo &  B^\triangleleft \\
\dTo	& 	& \dTo \\
\{ x \} \star \{b \} & \rTo &  \{x\} \star B^\triangleleft \\
\end{diagram}
\end{equation}
and this diagram is the pushout-join of $\emptyset \subset \{x\}$ with $\{b\} \subset B^\triangleleft$.
The latter is well-known to be left anodyne: it is a left adjoint, and
one can thus use the same reasoning as in Proposition \ref{propreedyfiltrationstep}. Hence the diagram (\ref{diagrampushoutinsertingmiddlexaempty})
 is infinity-coCartesian in
$\Cat_\infty$ by \cite[Lemma 2.1.2.3]{LUHTT}.
We now present the diagram (\ref{diagrampushoutinsertingmiddlex}) as
\begin{diagram}[small]
\{a\} \star \{b \} & \rTo & \{a \} \star B^\triangleleft & \rTo  & A^\triangleright \star B^\triangleleft \\
\dTo	& 	& \dTo & & \dTo \\
\{a\} \star \{ x \} \star \{b \} & \rTo & \{a \}  \star \{x\} \star B^\triangleleft & \rTo & A^\triangleright \star \{x\} \star B^\triangleleft. \\
\end{diagram}
By using the fact that joins preserve connected homotopy colimits \cite[Lemma 4.14]{SHAH}  and the argument for the diagram (\ref{diagrampushoutinsertingmiddlexaempty})
and its dual for $B = \emptyset$, we conclude that (\ref{diagrampushoutinsertingmiddlex})
is a pushout in $\Cat_\infty$.

To recapitulate, we have the following diagram in $\Cat_\infty$ (canonically induced
from a diagram in $\SSet$)
\begin{equation}
\label{dialurieremzigzag}
\begin{diagram}[small]
\Fun(Lat(x) \star \{x\} \star Mat(x),\cY) & \rTo & \Fun(Lat(x) \star Mat(x),\cY) \\
\uTo^\sim																													&						&				\uTo^\sim					\\
\Fun'(Lat(x)^\triangleright \star \{x\} \SEpbk \star Mat(x)^\triangleleft, \cY) & \rTo & \Fun'(Lat(x)^\triangleright \star Mat(x)^\triangleleft,\cY) \\
\dTo																													&						&				\dTo					\\
\Fun(Lat(x)^\triangleright \star \{x\} \SEpbk \star Mat(x)^\triangleleft,\cY) & \rTo & \Fun(Lat(x)^\triangleright \star Mat(x)^\triangleleft,\cY) \\
\dTo																													&						&				\dTo					\\
\Fun([2],\cY) 	&	\rTo	&	\Fun([1], \cY) \\
\end{diagram}
\end{equation}
and inverting in $\Cat_\infty$ the upper equivalences, we get the bottom square of (\ref{dialurieremarkproven}).
 \endproof

We would like to also address the functoriality of the diagram (\ref{dialurieremarkproven})
of Proposition \ref{proplurieremarkproven}. A typical functor which we would
be interested in is the projection to the infinity-localisation $p:\cM \to L \cM$
of a model category $\cM$. Given a diagram $X: Lat(x) \to \cM$, the image of
$\colim X$ in $L\cM$ ceases in general to be the colimit of $p \circ X$. It remains such,
however, if $X$ is Reedy-cofibrant \cite[Proposition 2.5.6]{HARPAZ2}. This motivates
the following definition.
\begin{opr}
Let $\cR$ be a Reedy category and $F: \cX \to \cY$ be an infinity-functor between quasicategories
both admitting $Lat(x)$-colimits and $Mat(x)$-limits for $x$. Let $\cR_\beta$ be the
step of a good filtration that adds $x$. A
diagram $X:\cR_{<\beta} \to \cX$ is called \emph{$(F,x)$-compatible} if the following holds:
\begin{itemize}
\item[i.] The colimit cone $Lat(x)^\triangleright \to \cX$
obtained by restricting $X$ to $Lat(x)$ and then taking the colimit, remains a colimit
cone after postcomposing with $F$,
\item[ii.] The limit cone $Mat(x)^\triangleleft \to \cX$
obtained by restricting $X$ to $Mat(x)$ and then taking the limit, remains a limit
cone after postcomposing with $F$.
\end{itemize}
A general diagram $\cR \to \cX$ is $(F,x)$-compatible iff its restriction to $\cR_{< \beta}$
is such. Denote by $\Fun_{(F,x)}(\cR_{<\beta},\cX)$, $\Fun_{(F,x)}(\cR_{\beta},\cX)$ and $\Fun_{(F,x)}(\cR,\cX)$
the full subcategories spanned by $(F,x)$-compatible functors.
\end{opr}

\begin{prop}
\label{proplurieremfunctoriality}
Let $\cR$ be a Reedy category and and $F: \cX \to \cY$ be an infinity-functor between quasicategories
both admitting $Lat(y)$-colimits and $Mat(y)$-limits for each $y$ in $\cR$. Then,
if $\cR_\beta$ is obtained from $\cR_{< \beta}$ by adding $x$,
the diagram (\ref{dialurieremarkproven}) of Proposition \ref{proplurieremarkproven}
induces the diagram $\cD ia(\cR_\beta,F)$:
\begin{equation}
\label{diaadmissibleconstantfibration}
\begin{diagram}[small]
\Fun_{(F,x)}(\cR_\beta, \cX) & \rTo & \Fun([2],\cX ) \\
\dTo 									&				& \dTo							\\
\Fun_{(F,x)}(\cR_{< \beta}, \cX) & \rTo		&	\Fun([1], \cX); 		\\
\end{diagram}
\end{equation}
the upper horizontal functor sends $S$ to $\colim_{Lat(x)} S \to S(x) \to \lim_{Mat(x)} S$,
and the bottom horizontal functor sends $S'$ to $\colim_{Lat(x)} S' \to \lim_{Mat(x)} S'$.

Moreover, we have an induced natural transformation of diagrams
$$
F_*:\cD ia(\cR_\beta,F) \longrightarrow \cD ia (\cR_\beta, Id_\cY)
$$
with $\cD ia (\cR_\beta, Id_\cY)$ denoting the diagram (\ref{dialurieremarkproven}) for $\cY$.
\end{prop}
\proof
The definition of compatibility implies that the following square is pullback in $\SSet$:
\begin{diagram}[small]
\Fun_{(F,x)}(\cR_\beta, \cX) & \rTo & \Fun(\cR_\beta, \cX) \\
\dTo 									&				& \dTo							\\
\Fun_{(F,x)}(\cR_{< \beta}, \cX) & \rTo		&	\Fun(\cR_{< \beta}, \cX); 		\\
\end{diagram}
the right functor is moreover a categorical fibration, since $\cR_{< \beta} \subset \cR_\beta$
is an inclusion. This diagram is hence a homotopy pullback, so the existence of
(\ref{diaadmissibleconstantfibration}) is clear. To see the functoriality, write
the diagram (\ref{dialurieremarkproven}) as a zig-zag, using the diagram (\ref{dialurieremzigzag}).
A close inspection shows that $F:\cX \to \cY$ induces a map between the zig-zags representing
(\ref{diaadmissibleconstantfibration}) and (\ref{dialurieremarkproven}): the nontrivial moment is in observing that
the restriction of the functor
$$\Fun'(Lat(x)^\triangleright \star \{x\} \star Mat(x)^\triangleleft, \cX)  \to  \Fun(Lat(x)^\triangleright \star \{x\} \star Mat(x)^\triangleleft, \cX')$$
to the diagrams $Lat(x)^\triangleright \star \{x\} \star Mat(x)^\triangleleft \to \cX$
which come from $(F,x)$-compatible diagrams, factors through $\Fun'(Lat(x)^\triangleright \star \{x\} \star Mat(x)^\triangleleft, \cX')$.
\endproof

\subsection{Induction for higher-categorical sections}

We need to introduce the notion of the higher-categorical restriction to
the fibre. For any object $x \in \cR$, the natural inclusion functor $Lat(x) \subset \cR$
extends to a functor
$$
C_x: Lat(x) \times [1] \to \cR
$$
which sends $(y \to x,0)$ to $y$ and $(y \to x,1)$ to $x$. There is a similar extension
for $Mat(x) \subset \cR$.

\begin{opr}
Let $\cX \to \cR$ be a coCartesian fibration and $x \in \cR$. A \emph{$x$-left restriction} of $S \in \Sect(\cR,\cX)$
is the infinity-functor $L_x S: Lat(x) \to \cX(x)$ defined as follows. In the diagram
in the diagram
\begin{diagram}[small,nohug]
Lat(x)						& \rTo^{S|_{Lat(x)}} &	\cX	\\
\dInto<{id \times \{0\}}				& 	 		&	\dTo		\\
Lat(x) \times [1]		&	\rTo_{C_x} & \cR		\\
\end{diagram}
choose a lifting $C_x S: Lat(x) \times [1]	\to \cX$ which sends each map
of the form $(y \to x,0) \to (y \to x,1)$ to a coCartesian morphism of $\cX$. Then
$L_x S$ is the restriction of $C_x S$ to $Lat(x) \times \{1\}$.

The definition of a \emph{$x$-right restriction} $R_x S:Mat(x) \to \cX(x)$
for a Cartesian fibration $\cX \to \cR$ is given dually.
\end{opr}

Arbitrary $x$-left and $x$-right restrictions always exist, and are defined up to
an equivalence:

\begin{lemma}
\label{lemmaexistenceofextensionsforfibrations}
For a coCartesian
fibration $\cX \to \cR$, the restriction functor $\Sect(Lat(x) \times [1],\cX) \to \Sect(Lat(x),\cX)$
induces an equivalence between $\Sect(Lat(x),\cX)$ and a full subcategory
$$\Sect_{[1]-cart}(Lat(x) \times [1],\cX) \subset \Sect(Lat(x) \times [1],\cX)$$
consisting of all sections that send maps of the form $(y \to x,0) \to (y \to x,1)$ to coCartesian morphisms of $\cX$. There is a dual statement for Cartesian fibrations.
\end{lemma}
\proof
Recall the notion of right marked anodyne map
(called marked anodyne in \cite[Definition 3.1.1.1]{LUHTT}), and the dual notion of left marked
anodyne map.
The inclusion $Lat(x)^\flat \hookrightarrow Lat(x)^\flat \times [1]^\sharp$
is left marked anodyne, which follows from the dual of \cite[Proposition 3.1.2.3]{LUHTT}. Observe
that the functor $\Sect_{[1]-cart}(Lat(x) \times [1],\cX) \to \Sect(Lat(x),\cX)$ is
precisely given by the equivalence \cite[Remark 3.1.3.4]{LUHTT}
$$\Map^\flat_\cR(Lat(x)^\flat \times [1]^\sharp,\cX^\natural) \stackrel \sim \to \Map^\flat_\cR(Lat(x)^\flat,\cX^\natural),$$
where $\cX^\natural$ means that we mark the coCartesian arrows in $\cX$.
Thus liftings $C_x S$ in the diagram
\begin{diagram}[small,nohug]
Lat(x)^\flat									& \rTo^{S|_{Lat(x)}} &	\cX^\natural	\\
\dInto<{id \times \{0\}}				& 	\ruDashto^{C_x S}			&	\dTo		\\
Lat(x)^\flat \times [1]^\sharp		&	\rTo_{C_x} & \cR^\sharp		\\
\end{diagram}
always exist and have the required property, and the space of such liftings is contractible. \endproof

\begin{lemma}
\label{lemmapreservationofleftrestr}
Let $F: \cX \to \cY$ be a coCartesian morphism between coCartesian fibrations
over $\cR$. Then $F$ preserves $x$-left restrictions for each $x \in \cR$.
\end{lemma}
\proof
The lemma follows from the existence and commutativity of the following diagram:
\begin{diagram}[small]
\Sect(Lat(x),\cX) & \lTo^\sim &  \Sect_{[1]-cart}(Lat(x) \times [1],\cX) & \rTo^{ev_1} & \Fun(Lat(x),\cX(x))	\\
\dTo							&						&				\dTo												&										&			\dTo									\\
\Sect(Lat(x),\cY) & \lTo^\sim &  \Sect_{[1]-cart}(Lat(x) \times [1],\cY) & \rTo^{ev_1} & \Fun(Lat(x),\cY(x)).	\\
\end{diagram}
\endproof

It will be illustrative to first start with the case of a direct Reedy category $\cD$ and a coCartesian
fibration of quasicategories $\cX \to \cD$.

\begin{lemma}
\label{lmmarkingofatrivialcofibration}
Let $S \in \SSet$ and $i:A \to B$ be a trivial cofibration (for the Joyal model structure)
in $\SSet/ S$. Then for any subset $\cE \in A(1)$ containing all degenerate edges,
the induced map of marked simplicial sets, $(A,\cE) \to (B,i\cE)$, is a trivial
cofibration in both Cartesian and coCartesian model structures on $\SSet_+ / S^\sharp$.
\end{lemma}
\proof
From \cite[Proposition 3.1.5.3]{LUHTT} we know that the functor $(-)^\flat:
\SSet/S \to \SSet/S^\sharp$ is left Quillen, for either Cartesian or coCartesian
model structure (one can check that the proof is self-dual).
We then have the following diagram in $\SSet/S^\sharp$
\begin{diagram}[small]
\cE \times [1]^\flat & \rTo & A^\flat & \rTo & B^\flat \\
\dTo									&			&	\dTo		&				&	\dTo 	\\
\cE \times [1]^\sharp	&	\rTo		&	(A,\cE) &	\rTo & (B,i\cE) \\
\end{diagram}
where for the leftmost arrow, we equip $\cE \times [1]$ with the map to $S$ induced by adjunction from $\cE \to S_1$.
The leftmost and the outer squares of this diagram are pushouts, hence the same
is true for the rightmost square. Thus $(A,\cE) \to (B,i\cE)$ is a pushout of
a trivial cofibration.
\endproof

\begin{lemma}
\label{lmjoinsofmarkedanodynes}
Let $f:A \to B$ be a left marked anodyne map and $X \in \SSet$. Then the join
$A \star X^\flat \to B \star X^\flat$ is left marked anodyne. Dually, a join
$X^\flat \star A \to X^\flat \star B$
of $X$ with a right marked anodyne map is right marked anodyne.
\end{lemma}
\proof We prove the left part. By \cite[Lemma 4.10]{SHAH}, the pushout-join of $f$ with
$\emptyset \to X^\flat$ is left marked anodyne. Unraveling the definition, the
pushout-join is given by
$$
B \coprod_{A} A \star X^\flat \to B \star X^\flat.
$$
Our map can be factored as a composition
$$
A \star X^\flat \to B \coprod_{A} A \star X^\flat \to B \star X^\flat,
$$
and the left map is a pushout of the left marked anodyne map $A \to B$. \endproof

\begin{lemma}
\label{lemmapushingtocoproductdirectedcase}
Let $\cD$ be a direct category and $\{ \cD_\beta \}$ denote a good filtration.
Then there are zig-zags of weak equivalences
in $\SSet_+/ \cD$ for the coCartesian model structure:
$$
Lat(x)^\flat \longrightarrow  Lat(x)^\flat \times [1]^\sharp
\longleftarrow (Lat(x)^\flat \times [1]^\sharp) \coprod_{Lat( x)_x^\flat} Lat( x)_x^\flat,
$$
$$
Lat(x)^\flat \star \{x\} \longrightarrow (Lat(x)^\flat \times [1]^\sharp) \star \{x\}
\longleftarrow (Lat(x)^\flat \times [1]^\sharp) \coprod_{Lat( x)_x^\flat} (Lat( x)_x^\flat \star \{x\}).
$$
Here the left maps are induced by the inclusion $\{0 \} \subset [1]^\sharp$ and
$Lat(x)^\flat_x$ denotes the category $Lat(x)^\flat$ together with a constant
functor to $\cD$ of value $x$.
\end{lemma}
\proof Denote by $M$ either $\{x\}$ or the empty simplicial set. The observations of Lemmas \ref{lemmaexistenceofextensionsforfibrations} and
\ref{lmjoinsofmarkedanodynes} instantly imply that $\Lat(x)^\flat \star M^\flat
\rightarrow (Lat(x)^\flat \times [1]^\sharp) \star  M^\flat $ is left marked
anodyne, and hence a coCartesian equivalence.

The remaining map is a pushout-join,
without any markings, of $\emptyset \subset M$ with $Lat(x) \subset Lat(x) \times [1]$,
with the latter map induced by the inclusion $\{1\} \subset [1]$, which is right
anodyne. The stability properties of right anodyne maps,
\cite[Corollary 2.1.2.7]{LUHTT}, and \cite[Lemma 2.1.2.3]{LUHTT}
again imply that the resulting pushout-join is inner anodyne. Lemma
\ref{lmmarkingofatrivialcofibration} then allows to conclude that
$$
(Lat(x)^\flat \times [1]^\sharp) \star M^\flat
\longleftarrow (Lat(x)^\flat \times [1]^\sharp) \coprod_{Lat( x)_x^\flat} (Lat( x)_x^\flat \star M)
$$
is a coCartesian eqiuvalence.
\endproof

\begin{prop}
\label{propinfinityreedyinductiondirect}
Let $\cX \to \cD$ be a fibrewise cocomplete coCartesian fibration over a direct
Reedy category $\cD$. Let ${\cD_\beta}$ denote a good filtration of $\cD$,
so that $\cD_\beta$ is obtained from
$\cD_{< \beta}$ by adding an object $x \in \cD$. Then there is a Cartesian
square in the quasicategory $\Cat_\infty$
\begin{diagram}[small]
\Sect(\cD_\beta, \cX) & \rTo & \Fun([1],\cX(x)) \\
\dTo 									&				& \dTo							\\
\Sect(\cD_{< \beta}, \cX) & \rTo		&	 \cX(x) 		\\
\end{diagram}
where
\begin{enumerate}
\item the left vertical arrow is induced by the restriction along
$\cD_{< \beta} \to \cD_\beta$,
\item the right vertical arrow is induced by the inclusion $[0] \to [1]$ skipping
$1 \in [1]$,
\item the bottom horizontal arrow sends $S \in \Sect(\cD_{< \beta}, \cX)$ to
$\colim_{Lat(x)} L_x S$. Here $L_x S$ is a $x$-left
restriction of $S$,
\item the upper horizontal arrow sends $S \in \Sect(\cD_\beta, \cX)$ to
$\colim_{Lat(x)} L_x S \to S(x)$.
\end{enumerate}
\end{prop}
\proof
Proposition \ref{propreedyfiltrationstep} implies that the square
\begin{diagram}[small]
\Sect(\cD_\beta, \cX) & \rTo & \Sect(Lat(x) \star \{x\}, \cX) \\
\dTo 									&				& \dTo							\\
\Sect(\cD_{< \beta}, \cX) & \rTo		&	 \Sect(Lat(x), \cX)  		\\
\end{diagram}
is Cartesian in $\Cat_\infty$. Lemma \ref{lemmapushingtocoproductdirectedcase}
then implies the existence of the following pullback square in $\Cat_\infty$:
\begin{diagram}[size=2.5em]
\Sect(Lat(x) \star \{x\}, \cX) & \rTo & \Sect((Lat(x)^\flat \times [1]^\sharp) \coprod_{Lat( x)_x^\flat} (Lat( x)_x^\flat \star \{x\}), \cX) \\
\dTo 									&				& \dTo							\\
\Sect(Lat(x), \cX) & \rTo		&	 \Sect((Lat(x)^\flat \times [1]^\sharp) \coprod_{Lat( x)_x^\flat} Lat( x)_x^\flat, \cX).  		\\
\end{diagram}
For the right arrow, the notation $\Sect(-,\cX)$ means $\Map^\flat_{\cD_\beta}(-,\cX^\natural)$. We notice that there is the following pullback square:
\begin{diagram}[size=2.5em]
\Sect((Lat(x)^\flat \times [1]^\sharp) \coprod_{Lat( x)_x^\flat} (Lat( x)_x^\flat \star \{x\}), \cX) & \rTo & \Fun(Lat(x) \star \{x\}, \cX(x)) \\
\dTo 									&				& \dTo							\\
\Sect((Lat(x)^\flat \times [1]^\sharp) \coprod_{Lat( x)_x^\flat} Lat( x)_x^\flat, \cX) & \rTo & \Fun(Lat(x),\cX(x))  		\\
\end{diagram}
which is induced by applying $\Sect(-,\cX)$ to a (homotopy) pushout square of
marked simplicial sets over $\cD$. Combining all the squares with Proposition \ref{proplurieremarkproven},
we finish the proof.
\endproof

To continue, we need to specify which fibrations we are going to consider over a
general Reedy category $\cR$. The higher-categorical generality that we choose to work with in this paper
is motivated by Lemma \ref{semifibexm}.
\begin{opr}
Let $\cR$ be a Reedy category. A \emph{left Reedy fibration} is a coCartesian fibration
of quasicategories $\cX \to \cR$ that is also a Cartesian (equivalently \cite[Corollary 5.2.2.4]{LUHTT} locally
Cartesian) fibration over the subcategory $\cR_-$. The notion of a right Reedy fibration
is given dually.
\end{opr}

Henceforth we shall stick with a left Reedy fibration $\cX \to \cR$. The corresponding
results for a right Reedy fibration can be obtained by dualisation.

\begin{lemma}
\label{lemmapushingtocoproduct}
Let $\cR$ be a Reedy category and $\{ \cR_\beta \}$ denote a good filtration as before.
Let $M$ denote either $\Mat(x)$ or $x \star Mat(x)$. Then there is a zig-zag of weak equivalences
in $\SSet_+/ \cR_\beta$ for the coCartesian model structure:
$$
\Lat(x)^\flat \star M^\flat \longrightarrow (\Lat(x)^\flat \times [1]^\sharp) \star M^\flat
\longleftarrow (Lat(x)^\flat \times [1]^\sharp) \coprod_{Lat( x)_x^\flat} (Lat( x)_x^\flat \star M).
$$
Here the left map is induced by the inclusion $\{0 \} \subset [1]^\sharp$ and
$Lat(x)^\flat_x$ denotes the category $Lat(x)^\flat$ together with a constant
functor to $\cR$ of value $x$. A dual result can be formulated for $Mat(x)$ and
the Cartesian model structure.
\end{lemma}
\proof Verbatim Lemma \ref{lemmapushingtocoproductdirectedcase}. \endproof

\begin{prop}
\label{propinfinityreedyinduction}
Let $\cX \to \cR$ be a fibrewise bicomplete left Reedy fibration. Let
${\cR_\beta}$ be a good filtration of $\cR$, so that $\cR_\beta$ is obtained from
$\cR_{< \beta}$ by adding an object $x \in \cR$. Then there is a Cartesian
square in the quasicategory $\Cat_\infty$
\begin{equation}
\label{diapropinfinityreedyinduction}
\begin{diagram}[small]
\Sect(\cR_\beta, \cX) & \rTo & \Fun([2],\cX(x)) \\
\dTo 									&				& \dTo							\\
\Sect(\cR_{< \beta}, \cX) & \rTo		&	\Fun([1],\cX(x))		\\
\end{diagram}
\end{equation}
where
\begin{enumerate}
\item the left vertical arrow is induced by the restriction along
$\cR_{< \beta} \to \cR_\beta$,
\item the right vertical arrow is induced by the inclusion $[1] \to [2]$ skipping
$1 \in [2]$,
\item the bottom horizontal arrow sends $S \in \Sect(\cR_{< \beta}, \cX)$ to
$\colim_{Lat(x)} L_x S \to \lim_{Mat(x)} R_x S$. Here $L_x S$ and $R_x S$ are $x$-left
and right restrictions of $S$,
\item the upper horizontal arrow sends $S \in \Sect(\cR_\beta, \cX)$ to
\begin{diagram}[small]
											&					& S(x)	&				&		\\
											& \ruTo		&				&	\rdTo	&		\\
\colim_{Lat(x)} L_x S	& 				&\rTo 	&				& \lim_{Mat(x)} R_x S. \\
\end{diagram}
\end{enumerate}
Moreover,
\begin{itemize}
\item[i.] For each ordinal $\beta$, the  $\Cat_\infty$-limit of $\{  \Sect(\cR_{\gamma},\cX) \}_{\gamma < \beta} $ is equivalent, via the evident map, to $\Sect(\cR_{ < \beta},\cX)$.
\item[ii.] The $\Cat_\infty$-limit of $\{\Sect(\cR_{\beta},\cX)\}_{\beta}$ is equivalent,
via the evident map, to $\Sect(\cR,\cX)$.
\end{itemize}
\end{prop}
\proof
The last two statements, (i.) and (ii.), follow immediately from the properties of the filtration $\{ \cR_\beta \}$.
For the rest, use again Proposition \ref{propreedyfiltrationstep}, and from the diagram
\begin{diagram}[small]
Lat(x) \star Mat(x) & \rTo & \cR_{< \beta} & \rTo & \cX\\
\dTo																								&			&		\dTo & & \dTo		\\
Lat(x) \star \{x \} \star Mat(x) & \rTo & \cR_\beta & \rTo & \cR, \\
\end{diagram}
get a Cartesian diagram in $\Cat_\infty$
\begin{diagram}[small]
\Sect(\cR_\beta, \cX) & \rTo & \Sect(Lat(x) \star \{x \} \star Mat(x),\cX) \\
\dTo								&				&				\dTo \\
\Sect(\cR_{<\beta}, \cX) & \rTo &  \Sect(Lat(x) \star Mat(x),\cX). \\
\end{diagram}
Lemma \ref{lemmapushingtocoproduct} implies that sections over $Lat(x) \star Mat(x)$ can
be replaced with sections over $(Lat(x) \times [1]) \coprod_{Lat( x)_x} (Lat(x)_x \star Mat(x))$
which are coCartesian along certain edges coming from $[1]$ in $Lat(x) \times [1]$.
We can depict this by asserting the existence of the $\Cat_\infty$-Cartesian diagram:
{\small
\begin{diagram}[size=2.5em]
\Sect(Lat(x) \star \{x\} \star Mat(x) , \cX) & \rTo & \Sect((Lat(x)^\flat \times [1]^\sharp) \coprod_{Lat( x)_x^\flat} (Lat( x)_x^\flat \star \{x\} \star Mat(x)^\flat), \cX) \\
\dTo 									&				& \dTo							\\
\Sect(Lat(x) \star Mat(x), \cX) & \rTo		&	 \Sect((Lat(x)^\flat \times [1]^\sharp) \coprod_{Lat( x)_x^\flat} (Lat( x)_x^\flat \star Mat(x)^\flat), \cX)
\end{diagram}
}%
obtained from a zig-zag in $\SSet^{[1]}$ of length two. For the right arrow, we notice that there is the following pullback square:
{\small
\begin{diagram}[size=2.5em]
\Sect((Lat(x)^\flat \times [1]^\sharp) \coprod_{Lat( x)_x^\flat} (Lat( x)_x^\flat \star \{x\} \star Mat(x)^\flat), \cX) & \rTo & \Sect(Lat( x)_x \star \{x\} \star Mat(x), \cX ) \\
\dTo 									&				& \dTo							\\
\Sect((Lat(x)^\flat \times [1]^\sharp) \coprod_{Lat( x)_x^\flat} (Lat( x)_x^\flat \star Mat(x)^\flat), \cX) & \rTo & \Sect (Lat( x)_x \star  Mat(x),\cX ).
\end{diagram}
}%
Note that the pullback of $\cX \to \cR$ to $Lat( x)_x  \star \{x\} \star Mat(x)$ is
a genuine Cartesian fibration. A similar argument concerning the replacement of $Mat(x)$ supplies us with the following pullback square in
$\Cat_\infty$:
\begin{diagram}[small]
\Sect(Lat( x)_x \star \{x\} \star Mat(x), \cX ) & \rTo & \Fun(Lat(x) \star \{x\} \star Mat(x),\cX(x)) \\
\dTo 									&				& \dTo							\\
\Sect (Lat( x)_x \star  Mat(x),\cX ) & \rTo & \Fun(Lat(x) \star Mat(x),\cX(x)) .  		\\
\end{diagram}
Proposition \ref{proplurieremarkproven} then concludes the proof. \endproof

We conclude our discussion by outlining the functoriality of the diagram (\ref{diapropinfinityreedyinduction})
that will be useful for the purposes of the comparison.

\begin{opr}
\label{oprfxcompsect}
Let $\cR$ be a Reedy category and $F: \cX \to \cY$ be an infinity-functor between fibrewise
bicomplete left Reedy fibrations over $\cR$.
Let $\cR_\beta$ be the
step of a good filtration that adds $x$. A
section $X:\cR_{<\beta} \to \cX$ is called \emph{$(F,x)$-compatible} if the following holds:
\begin{itemize}
\item[i.] For each $y \in \cR_{<\beta}$, the functor $F$ preserves coCartesian arrows
starting with $X(y)$ and Cartesian arrows over $\cR_-$ ending with $X(y)$,
\item[ii.] The cone $Lat(x)^\triangleright \to \cX(x)$ obtained by taking a $x$-left restriction of $X$ and then taking the colimit, remains a colimit cone in $\cY(x)$ cone after postcomposing with $F$,
\item[iii.] The cone $Mat(x)^\triangleleft \to \cX(x)$ obtained by
taking a $x$-right restriction of $X$ and then taking the limit, remains a limit
cone in $\cY(x)$ after postcomposing with $F$.
\end{itemize}
A general section $\cR \to \cX$ is $(F,x)$-compatible iff its restriction to $\cR_{< \beta}$
is such. Denote by $\Sect_{(F,x)}(\cR_{<\beta},\cX)$, $\Sect_{(F,x)}(\cR_{\beta},\cX)$ and $\Sect_{(F,x)}(\cR,\cX)$
the full subcategories spanned by $(F,x)$-compatible functors.
\end{opr}
An obvious compatibility is the following:
\begin{lemma}
\label{lemmapreservationofleftrestr2}
Let $F: \cX \to \cY$ be an infinity-functor between fibrewise
bicomplete left Reedy fibrations over $\cR$. Then $F$ preserves $x$-left restrictions and $x$-right restrictions of $(F,x)$-compatible sections.
\end{lemma}
\proof Immediate using (i.) of Definition \ref{oprfxcompsect}. \endproof

\begin{prop}
\label{propfunctorialityofreedyinduction}
Let $\cR$ be a Reedy category and and $F: \cX \to \cY$ be an infinity-functor between fibrewise
bicomplete left Reedy fibrations over $\cR$. Then,
if $\cR_\beta$ is obtained from $\cR_{< \beta}$ by adding $x$,
the diagram (\ref{diapropinfinityreedyinduction}) of Proposition \ref{propinfinityreedyinduction}
induces the diagram $\cD ia(\cR_\beta,F)$:
\begin{equation}
\label{diaadmissiblearbitraryfibration}
\begin{diagram}[small]
\Sect_{(F,x)}(\cR_\beta, \cX) & \rTo & \Fun([2],\cX(x) ) \\
\dTo 									&				& \dTo							\\
\Sect_{(F,x)}(\cR_{< \beta}, \cX) & \rTo		&	\Fun([1], \cX(x)). 		\\
\end{diagram}
\end{equation}

Moreover, post-composing with $F$ induces a natural transformation of diagrams
$$
F_*:\cD ia(\cR_\beta,F) \to \cD ia (\cR_\beta, Id_\cY)
$$
with $\cD ia (\cR_\beta, Id_\cY)$ denoting the diagram (\ref{diapropinfinityreedyinduction}) for $\cY \to \cR$.
\end{prop}
\proof
The definition of compatibility implies that the following square is pullback in $\SSet$:
\begin{diagram}[small]
\Sect_{(F,x)}(\cR_\beta, \cX) & \rTo & \Sect(\cR_\beta, \cX) \\
\dTo 									&				& \dTo							\\
\Sect_{(F,x)}(\cR_{< \beta}, \cX) & \rTo		&	\Sect(\cR_{< \beta}, \cX); 		\\
\end{diagram}
the right functor is moreover a categorical fibration, since $\cR_{< \beta} \subset \cR_\beta$
is an inclusion (and it becomes a cofibration in the coCartesian model structure on $\SSet_+/\cR$). This diagram is hence a homotopy pullback, so the existence of
(\ref{diaadmissiblearbitraryfibration}) is clear. To see the functoriality, write
the diagram (\ref{diapropinfinityreedyinduction}) as a zig-zag.
A close inspection shows that $F:\cX \to \cY$ induces a map between the zig-zags representing
(\ref{diaadmissiblearbitraryfibration}) and (\ref{diapropinfinityreedyinduction}): the nontrivial moments are
covered by Lemma \ref{lemmapreservationofleftrestr2} and observing that
the restriction of the functor (cf notation for the diagram (\ref{dialurieremzigzag}))
$$\Fun'(Lat(x)^\triangleright \star \{x\} \star Mat(x)^\triangleleft, \cX(x))  \to  \Fun(Lat(x)^\triangleright \star \{x\} \star Mat(x)^\triangleleft, \cY(x))$$
to the diagrams $Lat(x)^\triangleright \star \{x\} \star Mat(x)^\triangleleft \to \cX(x)$
which come from $(F,x)$-compatible sections, factors through $\Fun'(Lat(x)^\triangleright \star \{x\} \star Mat(x)^\triangleleft, \cY(x))$.
\endproof

\subsection{Families of relative categories}

\begin{opr}
Let $(\cX,\cW)$ be a relative category and $p:\cX \to \cC$ a functor sending $\cW$
to $\cC$-isomorphisms. An \emph{infinity-localisation} of $p$ is the factoring
of it as an infinity-localisation $\cX \to L \cX$ along $\cW$ followed by a categorical fibration
$L \cX \to \cC$ of quasicategories.

An infinity-localisation of $p$ is \emph{universal} iff for any functor $F:\cD \to \cC$,
the induced map $LF^* \cX \to F^* L \cX$ is a categorical equivalence.
\end{opr}
Using the model structure on marked simplicial sets \emph{over the point}, we can always (and even functorially)
factor $(\cX,\cW) \to (\cC,Iso_\cC) = \cC^\natural$ as a trivial cofibration followed by a fibration.
This shows that the infinity-localisation of functors always exists.

Let $\cE \to \cC$ be an opfibration. Assume that each fibre $\cE(c)$ is equipped
with weak equivalences $\cW(c)$ that are preserved by the transition functors: if
in a diagram
\begin{diagram}[small]
X					& \rTo & X' \\
\dTo^\sim		&				&	\dTo \\
Y					& \rTo & Y' \\
\end{diagram}
we have horizontal arrows opcartesian and $X \to Y \in \cW(c)$, then $X' \to Y' \in \cW(c')$.
The data of $\cW = Iso_\cE \cup (\cup_c \cW(c))$ provides (the nerve of) $\cE$ with the structure of a marked simplicial set.
Choose an infinity-localisation $(\cE,\cW) \stackrel \sim \hookrightarrow L \cE \twoheadrightarrow \cC^\natural$
of the functor $\cE \to \cC$.

Since $L\cE$ is fibrant in $\SSet_+$, it is a localisation of $\cE$.
One furthermore has \cite[Proposition 5.2.3]{MGEE1} that
$L \cE \to \cC$ is a coCartesian fibration classified by the functor $c \mapsto L \cE(c)$,
and the infinity-functor $\cE \to L \cE$ preserves coCartesian arrows.
Taking into account the naturality of the Grothendieck construction \cite[Remark 3.1.13]{MGEE1},
we also have that for any functor $F: \cD \to \cC$,
the natural morphism $L F^* \cE \to F^* L \cE$ is a coCartesian equivalence over $\cC$.
The localisation of $\cE \to \cC$ is thus universal.

\begin{rem}
\label{remdifferentlocalisers}
A result of Hinich \cite[2.1.4]{HINICH} is similar to that of Mazel-Gee \cite[Proposition 5.2.3]{MGEE1}
with a distinction: aside from inverting also the maps in the base, the fibration
localiser $\cW$ is assumed to be saturated. The way we construct
$\cW = Iso_\cE \cup (\cup_c \cW(c))$ does not usually produce a saturated localiser,
however in the Reedy setting this actually happens.
In many cases, $\cW$ saturates to a bigger class of weak equivalences: those maps
$X \to Y$ in $\cM$ which projection $f:x \to y$ is an isomorphism and the induced
map $f_! X \to Y$ is in $\cW(y)$.
\end{rem}

If we a dealing with, say, a Quillen presheaf $\cN \to \cC$, then the transition adjunctions
$$
\cN(c) \rightleftarrows \cN(c')
$$
do not in general preserve weak equivalences. Instead, they preserve weak
equivalences only if restricted to the subcategories of cofibrant or fibrant objects,
depending on the direction.
To study what happens if we localise $\cN$ in this case, let us first introduce some tools.

Consider a relative category
$(\cM,\cW)$. Given a subcategory $\cM_Q \subset \cM$ whose
objects we shall typically denote $QX$, for any composable sequence $X_0 \to ... \to X_n$
of morphisms in $\cM$, denote by $$\cM_Q /^{\cW} (X_0 \to ... \to X_n)$$ the category
whose objects are diagrams (with the bottom row always fixed)
\begin{diagram}[small]
QX_0 & \rTo & ... & \rTo & QX_n \\
\dTo^\sim 	&			&				&		&		\dTo^\sim \\
 X_0 & \rTo & ... & \rTo & X_n \\
\end{diagram}
with vertical arrows in $\cW$, and morphisms being weak equivalences between the
top rows. One can dually define an undercategorical analogue that we shall denote
$(X_0 \to ... \to X_n) \backslash^\cW \cM_Q.$
\begin{opr}
Let $(\cM,\cW)$ be a relative category.
\begin{enumerate}
\item A \emph{left approximation} is a full subcategory $\cM_Q \subset \cM$
such that for each sequence $X_0 \to ... \to X_n$, the category $\cM_Q /^{\cW} (X_0 \to ... \to X_n)$
is contractible.
\item A \emph{right approximation} is a full subcategory $\cM_R \subset \cM$
such that for each sequence $X_0 \to ... \to X_n$, the category $(X_0 \to ... \to X_n) \backslash^\cW \cM_R$
is contractible.
\end{enumerate}
\end{opr}

The following lemma justifies the choice of terminology.
\begin{lemma}
\label{lmleftapproximationsareapproximations}
Let $\cM_Q \subset \cM$ be a left approximation of a relative category $(\cM,\cW)$.
Then $$(\cM_Q,\cW|_{\cM_Q}) \to (\cM,\cW)$$ is a weak equivalence of relative categories.
There is a similar result that concerns right approximations.
\end{lemma}

Recall \cite{BARKAN} \cite[2.2]{HARPAZ2} that a functor $F:\cM \to \cM'$ of
relative categories is a \emph{weak equivalence}, or simply an equivalence,
if the induced map of the infinity-localisations $LF:L \cM \to L \cM'$
is an equivalence of quasicategories. This condition is equivalent to that, for each $n \geq 0$,
the functor
$$
F_*:\Fun_\cW([n],\cM) \longrightarrow \Fun_{\cW'}([n], \cM')
$$
is a weak equivalence of the associated nerves. Here $\Fun_\cW([n], \cM)$
denotes the subcategory of pointwise weak equivalences in $\Fun([n], \cM)$ and similarly
for $\cM'$.

\proof Denote by $\tilde \cM$ the full subcategory of
$\Fun([1],\cM)$ consisting of weak equivalences $QX \stackrel \sim \to X$ where
$QX$ belongs to $\cM_Q$. We have then a natural factorisation $\cM_Q \to \tilde \cM \to \cM$.
It is easy to see that $\cM_Q \to \tilde \cM$ is a weak equivalence of relative categories, since
it admits a homotopy inverse.

We observe that for each $n$, the functor $\Fun_\cW([n],\tilde \cM) \rightarrow \Fun_\cW([n], \cM)$
is an opfibration, hence by Quillen Theorem A it is sufficient to prove that its fibres are contractible.
However, the fibre of this functor over $X_0 \to ... \to X_n$ is exactly $\cM_Q /^{\cW} (X_0 \to ... \to X_n)$.
\endproof

As the choice of notation suggests, the fundamental examples \cite[1.3.7]{HINICH} \cite[Lemma 2.4.8]{HARPAZ2}
of left approximations are given by the subcategories $\cM_c \subset \cM$ of cofibrant
objects in a model category $\cM$. Fibrant objects, dually, serve as right approximations.
The identity functor is of course both left and right approximation. We will have further
examples below.

In the case of (co)fibrant objects in a model category, a stronger approximation
property holds.
\begin{opr} Let $(\cM,\cW)$ be a relative category.
A \emph{strong left approximation} is a full subcategory $\cM_Q \subset \cM$
such that
\begin{itemize}
\item[i.] for each $X \in \cM$, the category $\cM_Q /^\cW X$ is contractible, and
\item[ii.] for each $X \to Y \in \cM^{[1]}$, the projection $\cM_Q /^{\cW} (X \to Y) \to \cM_Q /^\cW X$
has contractible fibres.
\end{itemize}
the definition of a strong right approximation is given dually.
\end{opr}

\begin{lemma}
Any strong left approximation is a left approximation, and dually, any strong right
approximation is a right approximation.
\end{lemma}
\proof The projection functor $$\cM_Q /^{\cW} (X_0 \to ... \to X_n) \to \cM_Q /^{\cW} (X_0 \to ... \to X_{n-1})$$
is a Grothendieck fibration whose fibres are the same as the fibres of the projection
$$\cM_Q /^{\cW} (X_{n-1} \to X_n) \to \cM_Q /^{\cW}   X_{n-1}. $$ Quillen Theorem A
and induction imply the result. \endproof

\begin{lemma}
Let $\cM$ be a model category. The inclusion $\cM_c \subset \cM$ of cofibrant objects
is a strong left approximation. The inclusion $\cM_f \subset \cM$ of fibrant objects
is a strong right approximation.
\end{lemma}
\proof We prove the left part. The contractibility of $\cM_c /^\cW X$ is already established
\cite[Lemma 2.4.8]{HARPAZ2}. We now have to prove that for each $X \to Y$ and a weak
equivalence $QX \stackrel \sim \to X$ with $QX$ cofibrant, the category of fillers
of the upper right corner in the square
\begin{diagram}[small]
QX & \rDashto  & QY \\
\dTo^\sim 	&		 	&		\dDashto^\sim \\
 X  & \rTo & Y \\
\end{diagram}
is contractible. Note that this will follow if we prove that for any map $f:QX \to Y$
with $QX$ cofibrant, the category $Fact(f)$ of factorisations
\begin{equation}
\label{diafactorisationofamap}
\begin{diagram}[small,nohug]
		&				&	QY	&				&				\\
		& \ruTo	&			&	\rdTo^\sim	&	\\
QX 	&  			& \rTo^f	&				&	Y
\end{diagram}
\end{equation}
with $QY$ cofibrant and $QY \to Y$ a weak equivalence, is contractible.

Denote by $Fact_c(f) \subset Fact(f)$ the subcategory consisting of all the
factorisations in which $QX \to QY$ is a cofibration. The comma-fibre of this functor
over an object (\ref{diafactorisationofamap}) is the category of cofibrant
replacements of (\ref{diafactorisationofamap}) in the model category $f \backslash (\cM / Y)$.
Hence $Fact_c(f) \to Fact(f)$ is a homotopy equivalence. The category $Fact_c(f)$
is however contractible: an object
\begin{diagram}[small,nohug]
		&				&	QY	&				&				\\
		& \ruInto	&			&	\rdTo^\sim	&	\\
QX 	&  			& \rTo^f	&				&	Y
\end{diagram}
is a cofibrant replacement of $f$ in $QX \backslash \cM$. \endproof

Left and right approximations interact reasonably well in families. The following
propositions will be sufficient for our needs. First, consider a typical situation
when one has a family of model categories and functors that preserve cofibrations and
trivial cofibrations. Such functors do
not preserve weak equivalences on the nose, but we can always limit our attention
to the cofibrant objects:

\begin{lemma}
\label{lemmastrongleftapproximationsintegrate}
Let $p:\cM \to \cC$ be an opfibration such that each fibre $\cM(c)$ is equipped with
weak equivalences $\cW(c)$ and a \emph{strong} left approximation $\cM_Q(c)$, subject to the condition
that given a diagram
\begin{diagram}[small]
QX					& \rTo & X' \\
\dTo^\sim		&				&	\dTo \\
QY					& \rTo & Y' \\
\end{diagram}
with left vertical arrow a fibrewise weak equivalence, horizontal arrows opcartesian,
and $QX,QY \in \cM_Q(c)$ for some $c$, we have that $X' \to Y'$ is a weak equivalence in $\cM_Q(c')$.

Denote by $\cM_Q \subset \cM$ the subcategory spanned by all objects from all $\cM_Q(c)$.
Then $\cM_Q \to \cC$ is an opfibration whose transition functors
preserve weak equivalences. Moreover, the inclusion $\cM_Q \subset \cM$ is a strong left
approximation, where we endow $\cM$ with the totality of fibrewise weak equivalences
$\cW = \cup_c \cW(c)$.
\end{lemma}

\proof The only nontrivial part is proving that $\cM_Q \subset \cM$ is a strong left
approximation. For a single object $X$, the category $\cM_Q /^\cW X$ is the same as $\cM_Q(x) /^{\cW(x)} X$ for $pX = x$. It is hence contractible by assumption.

For $\alpha: X \to Y$ we study the fibres of the forgetful functor
$$
U:\cM_Q /^\cW (X \to Y) \to \cM_Q /^\cW X.
$$
Over $QX \stackrel \sim \to X$, the fibre is the category of fillers of the right upper corner in the following square:
\begin{diagram}[small]
QX & \rDashto  & QY \\
\dTo^\sim  	&		 	&		\dDashto^\sim  \\
X   & \rTo & Y; \\
\end{diagram}
This is the same category as the category of fillers in the diagram
\begin{diagram}[small]
QX & \rDashto  & QY \\
\dTo^=  	&		 	&		\dDashto^\sim  \\
QX   & \rTo & Y  \\
\end{diagram}
for the composition $QX \to X \to Y$. The latter is the same as the fibre of
$$
U':\cM_Q(y) /^\cW (p \alpha_! QX \to Y) \to \cM_Q(y) /^\cW p \alpha_! X
$$
over $p \alpha_! QX \to  p \alpha_! QX$. \endproof

Given a relative category $(\cX,\cW)$, we say that $\cW$ is saturated if any map
that becomes invertible in $L \cX$, belongs to $\cW$. Weak equivalences of model
categories are saturated.

\begin{corr}
\label{corollarylocalisationinfamilies}
In the situation of Lemma \ref{lemmastrongleftapproximationsintegrate}, consider
an infinity-localisation $\cM \to L \cM \to \cC$ of $\cM \to \cC$. Then $L \cM \to \cC$ is a coCartesian
fibration of quasicategories, and is a \emph{universal} localisation of $\cM \to \cC$.

Assume that each $\cW(c)$ is saturated. Then a map $\alpha:X \to Y$ of $\cM$ becomes coCartesian in $L \cM$
iff for any $\cM_Q$-approximation $QX \to X$, the induced map $p\alpha_! QX \to Y$
is a weak equivalence.
\end{corr}
Note that the last sentence implies that the infinity-functor $\cM \to L \cM$
preserves coCartesian arrows with $\cM_Q$-domain.

\proof
Since $\cM_Q \subset \cM$ is a weak equivalence of relative categories,
one has $L \cM_Q \cong L \cM$. Hence $L \cM \to \cC$ is a categorical fibration
that is equivalent to a coCartesian fibration, hence it is a coCartesian fibration
itself \cite[Corollary 3.4]{MGEE2}. Moreover, for any $F: \cD \to \cC$, we have
an equivalence $L F^* \cM_Q \cong L F^* \cM$, and since $L \cM_Q \to \cC$ is universal,
same result follows for $L \cM \to \cC$.

Furthermore, the functor $L \cM_Q \cong L \cM$ preserves
coCartesian arrows, since in this case being coCartesian coincides with the
invariant definition \cite[Theorem 3.3]{MGEE2}, and equivalences preserve invariantly
defined coCartesian arrows. The infinity-functor $\cM_Q \to L\cM_Q$ preserves coCartesian arrows
as well.

The last part of the corollary then follows from the observation that for
$\alpha: X \to Y$ and a $\cM_Q$-approximation $QX \stackrel \sim \to X$,
we have the following diagram in $\cM$:
\begin{diagram}[small,nohug]
QX				&	\rTo	 & 	p\alpha_! QX \\
\dTo<\sim	&					&		\dTo \\
X					& \rTo_\alpha & Y, \\
\end{diagram}
the map $\alpha$ being $L \cM$-coCartesian thus being equivalent to $p\alpha_! QX \to Y$
becoming invertible in $L \cM$.
\endproof

Our next result is a variation of \cite[Theorem 5.3.1]{MGEE1}.
\begin{lemma}
\label{lemmaadjunctionsandapproximations}
Let $(\cC,\cW_\cC)$ and $(\cD,\cW_\cD)$ be two relative categories and
$$
F:\cC \rightleftarrows \cD:G
$$
an adjunction such that
\begin{itemize}
\item[i.] there is a left approximation $\cC_Q \subset \cC$ such that $F$ preserves
weak equivalences between the objects of $\cC_Q$,
\item[ii.] there is a right approximation $\cD_R \subset \cD$ such that $G$ preserves
weak equivalences between the objects of $\cD_R$.
\end{itemize}
Let $\cM \to [1]$ be the bifibration corresponding to the adjunction $F \dashv G$
that we endow with the weak equivalences $\cW=\cW_\cC \cup \cW_\cD$.
Denote by $\cM_l$ the subcategory of $\cM$ spanned by $\cC_Q$ and $\cD$, and by
$\cM_r$ the subcategory spanned by $\cC$ and $\cD_R$. Then the inclusions
$\cM_l \subset \cM \supset \cM_r$ are, respectively, a left and a right approximation.
These approximations are strong whenever $\cC_Q \subset \cC$ and $\cD_R \subset \cD$ are.
\end{lemma}
\proof We prove the left part. A typical string of arrows in $\cM$ is given by
$$X_0 \to ... \to X_m \to Y_0 \to ... \to Y_n$$
with $X_i$ in $\cC$ and $Y_j$ in $\cD$. Formally, $m$ and $n$ can take value $-1$,
so that no $X$ or $Y$ appear in the string. When this happens, we are reduced
to proving that $\cC_Q \subset \cC$ or $\cD = \cD$ are left approximations, which is given.

For the middle case, consider the functor
$$
\cM_l /^\cW (X_0 \to ... \to X_m \to Y_0 \to ... \to Y_n) \to \cC_Q /^{\cW_\cC} (X_0 \to ... \to X_m)
$$
this is a Grothendieck fibration whose fibres have final objects. We thus conclude
what is needed. The treatment of the strong case is similar to the proof of Lemma
\ref{lemmastrongleftapproximationsintegrate}.
\endproof

\begin{corr}
\label{corollaryadjunctionsgiveadjunctions}
In the situation of Lemma \ref{lemmaadjunctionsandapproximations}, consider an infinity-localisation
$\cM \to L \cM \to [1]$. Then $L \cM \to [1]$ is a biCartesian fibration that classifies an
infinity-adjunction
$$
\bL F: L\cC \rightleftarrows L \cD: \bR G.
$$
The value of $\bL F$ at $X$ is represented by $F QX$, and the value of $\bR G$
at $Y$ is represented by $G RY$.

The value of the unit of the adjunction $\bL F \dashv \bR G$ at $X \in L \cC$ is equivalent to the image
of the following chain of $\cC$-morphisms:
$$
X \stackrel \sim \leftarrow QX \to G F Q X \to G (R F Q X)
$$
with $QX \stackrel \sim \to X$ and $FQX \stackrel \sim \to RFQX$ being the approximations. The expression for the counit is
given dually.
\end{corr}
\proof Lemma \ref{lemmaadjunctionsandapproximations} allows us to conclude
that $L \cM_l \cong L \cM \cong L \cM_r$, hence $L \cM \to [1]$ is a categorical fibration
that is equivalent to both a Cartesian and a coCartesian fibration. It is also universal, so
$L \cC \cong L \cM(0)$ and $L \cD \cong L \cM(1)$.
Same analysis as in Corollary \ref{corollarylocalisationinfamilies}
permits us to conclude that $\alpha: X \to Y$ of $\cM$ becomes coCartesian in
$L \cM$ if $p\alpha_! QX \to Y$ is a weak equivalence, and similarly for becoming Cartesian.

To see the unit statement, recall \cite[Proposition 5.2.2.8]{LUHTT}. Given a biCartesian
fibration $\cB \to [1]$ and $x \in \cB(0)$, first choose a coCartesian arrow
$x \to f_! x$ covering $0 \to 1$, and then a Cartesian arrow $f^* f_! x \to f_! x$
covering the same arrow. The induced arrow $x \to f^* f_! x$ is (equivalent to)
the value of the unit of the adjunction $f_! \dashv f^*$ at $x$.

We now observe that for $QX \in \cC_Q$, the composition $QX \to F QX \stackrel \sim \to RFQX$
becomes coCartesian in $L \cM$. Since Cartesian maps with $\cD_R$-codomains remain
Cartesian, the resulting map $QX \to G(RFQX)$ represents the unit of the infinity-adjunction
$\bL F \dashv \bR G$ evaluated at $QX$. However, it admits a factorisation $QX \to G F Q X \to G (R F Q X)$
in $\cC$. \endproof

We conclude this subsection with a couple of lemmas that will permit us to conclude
that the higher-categorical Reedy induction holds for the localisations of model
categories of sections.
For a relative category $(\cM,\cW)$, denote by $L^H \cM$ its Dwyer-Kan ``hammock''
localisation.
\begin{lemma}
\label{lmrecognitionhomotopypullbacks}
Let
\begin{diagram}[small]
(\cM_1, \cW_1) & \rTo^F & (\cM_2,\cW_2) \\
\dTo<G					&				& 		\dTo>K \\
(\cM_3, \cW_3) & \rTo_H & (\cM_4,\cW_4) \\
\end{diagram}
be a diagram of relative categories and equivalence preserving functors,
such that
\begin{itemize}
\item[i.] The induced square
\begin{diagram}[small]
\Ob \cM_1 & \rTo &  \Ob \cM_2   \\
\dTo					&				& 		\dTo \\
\Ob \cM_3 & \rTo & \Ob \cM_4 \\
\end{diagram}
is a pullback.
\item[ii.] for each $X,Y \in \cM_1$, the induced square of simplicial sets
\begin{diagram}[small]
L^H \cM_1(X,Y) & \rTo &  L^H \cM_2(FX,FY)   \\
\dTo					&				& 		\dTo \\
L^H \cM_3(GX,GY) & \rTo &  L^H \cM_4(KFX,KFY) \\
\end{diagram}
is a homotopy pullback.
\item[iii.] The functor $\Ho K : \Ho \cM_2 \to \Ho \cM_4$ is an isofibration.
\end{itemize}
Then the induced square
\begin{equation}
\label{diapullbacklemma}
\begin{diagram}[small]
L^H \cM_1 & \rTo  &  L^H \cM_2    \\
\dTo 				&				& 		\dTo  \\
L^H \cM_3  & \rTo  &  L^H \cM_4  \\
\end{diagram}
\end{equation}
is a homotopy pullback of simplicial categories.
\end{lemma}
\proof We use the key observation of \cite[Lemma 3.1.11]{HARPAZ} asserting
that a functor $F: \cC \to \cD$ of simplicial categories can be factored
as a weak equivalence and an isomorphism on object-sets followed by a map
that is a fibration on mapping spaces. Using this observation we replace the diagram
(\ref{diapullbacklemma}) by a weakly equivalent diagram
\begin{diagram}[small]
L^H \cM_1 & \rTo  &  (L^H \cM_2)'    \\
\dTo 				&				& 		\dTo  \\
L^H \cM_3  & \rTo  &  (L^H \cM_4)'  \\
\end{diagram}
where $(L^H \cM_2)' \to (L^H \cM_4)'$ is now a genuine fibration of simplicial categories.
It is still true that the diagram of object-sets is a pullback, and an analogue of
(ii.) holds for this diagram as well. It thus follows that
$$L^H \cM_1 \to L^H \cM_3  \times_{(L^H \cM_4)'} (L^H \cM_2)'$$
is bijective on objects and homotopically fully faithful. Since the model structure
on simplicial categories is right proper, our proof is complete. \endproof

\begin{lemma}
\label{lmrecognitionhomotopyinvlims}
Let $\alpha$ be an ordinal and assume given, for each $\beta < \alpha$, a simplicial
category $A_\beta$, and for each $\beta' < \beta$, a simplicial functor $A_\beta \to A_{\beta'}$
rendering the assignment $\beta \mapsto A_\beta$ functorial.
Assume furthermore given,
for each $\beta$, a simplicial functor $A_\beta \to A_{< \beta}$, and for each
$\gamma < \beta$, a commutative diagram
\begin{diagram}[small]
		&					&	A_\beta	& 				&		\\
		& \ldTo		&					&	\rdTo	&			\\
A_{< \beta}		&					&	\rTo	& 				&	A_\gamma	\\
\end{diagram}
exhibiting $A_{< \beta}$ as a cone of $\{ A_{\gamma} \}_{\gamma < \beta}$.

We require that the following holds.
\begin{enumerate}
\item[0.] For each successor ordinal the functor $A_{\beta+1} \to A_{<\beta+1}$ is
equal to the functor $A_{\beta + 1} \to A_\beta$,
\item For each ordinal $\beta$, the functor $A_\beta \to A_{<\beta}$ is an isofibration
after applying $\Ho$,
\item For each ordinal $\beta$, the functors $A_{< \beta} \to A_\gamma$
exhibit $\Ob A_{< \beta}$ as a limit of $\{\Ob A_\beta \}_{\gamma < \beta}$,
\item For each ordinal $\beta$ and $x,y \in A_{< \beta}$, denote by $x_\gamma,y_\gamma$ their image in
$A_\gamma$. Then the maps $A_{< \beta} \to A_\gamma$ exhibit the simplicial set $A_{< \beta}(x,y)$ as
a homotopy limit of $\{ A(x_\gamma,y_\gamma)\}_{\gamma < \beta} $
\end{enumerate}
Let $A$ be a simplicial category and $A \to A_\beta$ be a compatible family of
functors such that
\begin{itemize}
\item[i.] The maps $A \to A_\beta$ exhibit $\Ob A$ as a limit of $\Ob A_\beta$,
\item[ii.] For each $x,y \in A$, denote by $x_\beta,y_\beta$ their image in
$A_\beta$. Then the maps $A \to A_\beta$ exhibit the simplicial set $A(x,y)$ as
a homotopy limit of $A(x_\beta,y_\beta)$.
\end{itemize}
Then for any ordinal $\beta$, $A_{<\beta}$ is the homotopy limit
of $\{A_\gamma \}_{\gamma < \beta}$, and $A$ is the homotopy limit of $\{A_\beta \}_{\beta < \alpha}$.
\end{lemma}
Note that due to (0.), the conditions (2.) and (3.) are nontrivial only for a limit
ordinal $\beta$.

\proof
Use transfinite induction to construct a new diagram of fibrations
between fibrant simplicial categories that we shall denote $\{B_\beta \}$,
together with suitable maps $B_\beta \to B_{< \beta}$.
For $\beta = 0$, the map $A_0 \to B_0$ is the weak equivalence that is an identity on
objects followed by a fibration $B_0 \to [0]$. For the inductive step,
set $B_{< \beta} = \lim_{\gamma < \beta} B_\gamma$. An
elementary lifting argument shows that for each $\gamma < \beta$, the projections
$B_{< \beta} \to B_\gamma$ are fibrations. Moreover $B_{< \beta}$ is
the limit of a diagram of fibrations between fibrant simplicial categories, hence
a homotopy limit as well. Conditions (2.) and (3.) imply that $A_{< \beta} \to B_{< \beta}$
is bijective on objects and weakly fully faithful, hence an
equivalence of simplicial categories.

We now factor the composition
$A_\beta \to A_{< \beta} \to B_{< \beta}$
as an object-bijective weak equivalence $A_\beta \to B_\beta$ followed by a local
fibration $B_\beta \to B_{< \beta}$. Since bijective-on-objects equivalences are
isofibrations, the simplicial functor $B_\beta \to B_{< \beta}$ will be a $\Ho$-isofibration
and hence a fibration as well. The induction is thus complete and the limit of
$\{B_\beta \}_{\beta < \alpha}$ is also the homotopy limit.

To complete the proof, note that $A \to B_\beta$ satisfy an analogue of (i.) and (ii.), hence the functor
$A \to \lim_\beta B_\beta$ is bijective on objects and weakly fully faithful, hence an
equivalence of simplicial categories.
\endproof

The following proposition will be of use in our induction later on, and
is a good illustration of one of the lemmas we just proved.

\begin{prop}[Cisinski, Hirschowitz-Simpson]
\label{propcomparisonforasimplex}
Let $\cM$ be a model category and $[n]$ a finite
ordinal. Then the functor $\Fun([n],\cM) \to \Fun([n],L \cM)$ induces an equivalence
of quasicategories
$$L\Fun([n],\cM) \cong \Fun([n],L \cM).$$
\end{prop}
\proof Equip $\Fun([n],\cM)$ with projective model structure. It will be sufficient to
prove that the induced functor $L\Fun([n],\cM)_{cf} \cong \Fun([n],L \cM_{cf})$ is
an equivalence.

The map $[n-1] \coprod_{[0]} [1] \to [n]$ (we include $[n-1]$ as a left interval in $[n]$)
is inner anodyne (one can express
it as a pushout-join and apply \cite[Lemma 2.1.2.3]{LUHTT}). This implies that for
any quasicategory $\cX$, the induced square
\begin{diagram}[small]
\Fun([n],\cX) & \rTo & \Fun([1],\cX) \\
\dTo 					&				&		\dTo 				\\
\Fun([n-1],\cX) &	\rTo		& \cX					\\
\end{diagram}
is a pullback in $\Cat_\infty$. In particular, it is true for $L \cM_{cf}$.

Note that we have a similar pullback square of 1-categories of fibrant-cofibrant objects:
\begin{diagram}[small]
\Fun([n],\cM)_{cf} & \rTo & \Fun([1],\cM)_{cf} \\
\dTo 					&				&		\dTo 				\\
\Fun([n-1],\cM)_{cf} &	\rTo		& \cM_{cf}					\\
\end{diagram}
where, again, we equip all functor categories with projective model structure. This
implies that for any $X \in \Fun([n],\cM)_{cf}$ and a simplicial resolution
$Y_\bullet \in \Fun([n],\cM)^{\Delta^\op}_{cf}$, the following square is a pullback in $\SSet$:
\begin{diagram}[small]
\Fun([n],\cM)_{cf}(X,Y_\bullet) & \rTo & \Fun([1],\cM)_{cf}(X,Y_\bullet) \\
\dTo 					&				&		\dTo 				\\
\Fun([n-1],\cM)_{cf}(X,Y_\bullet) &	\rTo		& \cM_{cf}(X,Y_\bullet),					\\
\end{diagram}
where we abuse the notation and identify $X,Y_\bullet$ and their images in the different
categories. Note that $Y_\bullet$ remains a simplicial resolution even when restricted to
$\Fun([n-1],\cM)$ and all other categories.
All simplicial sets participating in the diagram above are fibrant, and moreover,
the map
\begin{equation}
\label{eqcompsimplex}
\Fun([1],\cM)_{cf}(X,Y_\bullet) \to \cM_{cf}(X,Y_\bullet)
\end{equation}
is a Kan fibration:
denoting by $ev$ the evaluation-at-zero functor $\Fun([1],\cM)  \to \cM$ and by $ev_!$ its left adjoint,
we see that the map (\ref{eqcompsimplex}) is given by precomposing with $ev_! ev  X \to X$, which
is a projective cofibration by inspection.

Recall \cite[Lemma 6.1]{MANDELL} that given any model category $\cN$,
there exists a span of weak equivalences
of simplicial sets relating the mapping spaces $L^H \cN(X,Y)$ of the Dwyer-Kan localisation and homotopy
function complexes $\cN(X,Y_\bullet)$. Note that the span of \cite[Lemma 6.1]{MANDELL} is functorial,
both with respect to maps in $\cN$ and along functors $\cN_{cf} \to \cN'$ that preserve
weak equivalences.

 Applying the above observation to establish the equivalence between $\Fun([n],\cM)_{cf}(X,Y_\bullet)$ and $L^H \Fun([n],\cM)_{cf}(X,Y_0)$, we see that
 the conditions (i.) and (ii.) of Lemma \ref{lmrecognitionhomotopypullbacks} are satisfied.
 Since every isomorphism in $\Ho \cM_{cf}$ can be realised as a weak equivalence between
 fibrant-cofibrant objects, it is easy to check that the functor $\Fun([1],\cM)_{cf}  \to \cM_{cf} $
 is a $\Ho$-isofibration as well. Thus Lemma \ref{lmrecognitionhomotopypullbacks} (and a passage through
 homotopy coherent nerve) implies that
 \begin{diagram}[small]
 L\Fun([n],\cM)_{cf}  & \rTo & L\Fun([1],\cM)_{cf}  \\
 \dTo 					&				&		\dTo 				\\
 L\Fun([n-1],\cM)_{cf}  &	\rTo		& L\cM_{cf} ,					\\
 \end{diagram}
is a pullback in $\Cat_\infty$. Everything will thus follow by induction from the case $n=1$.

The infinity-functor $L\Fun([1],\cM)_{cf} \cong \Fun([1],L \cM_{cf})$ is essentially
surjective. This is true since for any cofibrant-fibrant $x,y$, the map
$$\cM_{cf}(x,y) \to \pi_0 L \cM_{cf}(x,y) \cong \Ho \cM_{cf}(x,y)$$
is surjective, as per well known fact for model categories. It remains to prove
fully faithfulness. Let $f:A \to B$ and $g_\bullet:C_\bullet \to D_\bullet$ be, respectively,
an object of $L\Fun([1],\cM)_{cf}$ and a simplicial resolution in $L\Fun([1],\cM)_{cf}^{\Delta^\op}$.
The following square is a pullback in $\SSet$:
\begin{diagram}[small]
 \Fun([1],\cM)_{cf}(f,g_\bullet) & \rTo &   \cM_{cf}(A,C_\bullet) \\
\dTo										&				&				\dTo					\\
 \cM_{cf}(B,D_\bullet) &		\rTo	&		  \cM_{cf}(A,D_\bullet), \\
\end{diagram}
and the bottom map is a fibration, being induced as a pullback along the cofibration
$A \to B$. Using \cite[Lemma 6.1]{MANDELL}, this implies again that
the mapping space $L^H \Fun([1],\cM)_{cf}(f,g_0)$ is equivalent to the homotopy
end of the functor
$$
[1]^\op \times [1] \stackrel{(f,g_0)}{\longrightarrow} L^H \cM_{cf}^\op \times L^H \cM_{cf}
\stackrel{\Hom}{\longrightarrow} \SSet,
$$
however, up to an application of the coherent nerve and a choice of a model of the
hom-functor \cite[Proposition 5.2.1.11]{LU}, it is exactly the same presentation
that is valid for the mapping spaces in $\Fun([1],L \cM_{cf})$, \cite[Proposition 5.1]{GEPNERNIK}.
\endproof

\subsection{The comparison}

Recall that a left model Reedy fibration is a functor $p:\cM \to \cR$ to a Reedy category
that is a Grothendieck opfibration, a Grothendieck fibration over $\cR_-$, such that
the associated semifibration (cf Lemma \ref{semifibexm}) is admissible in the sense
of Definition \ref{admissibilitydef}.

\begin{corr}
\label{corollarylocalisationofabifibration}
The infinity-localisation along the fibrewise weak equivalences $L p: L \cM \to \cR$ is universal, and is a left Reedy fibration of quasicategories.
The functor $\cM \to L \cM$ preserves coCartesian arrows with cofibrant domain,
and Cartesian arrows over $\cR_-$ with fibrant codomain.
\end{corr}
\proof An immediate specialisation of Corollary \ref{corollarylocalisationinfamilies}. \endproof

\begin{rem}[Cf Remark \ref{remdifferentlocalisers}]
Using the calculus of zig-zags, it is possible to show that the saturated
localiser giving rise to $L \cM$ coincides in fact with the totality of weak equivalences:
it suffices to see that a fibrewise map that gets inverted in $L \cM$, is inverted
in the fibre model category localisation $L \cM(c)$. By choosing a zig-zag representative
for the inverse, the statement becomes trivial, since a Reedy category has no non-identity
isomorphisms.
\end{rem}

The infinity-functor $\cM \to L \cM$ induces the functor
$\Sect(\cR,\cM) \to \Sect(\cR,L \cM)$ from (the nerve of) the 1-category $\Sect(\cR,\cM)$
to the sections of the localisation $L \cM \to \cR$.
Since fibrewise weak equivalences are inverted when projected to $L \cM$, we get a well defined
infinity-functor $$L \Sect(\cR,\cM) \to \Sect(\cR,L \cM)$$ from the quasicategorical
localisation of the model structure on $\Sect(\cR,\cM)$.

\begin{thm}
\label{thmcomparison}
Let $\cM \to \cR$ be a left model Reedy fibration. Then
the induced infinity-functor $L \Sect(\cR,\cM) \to \Sect(\cR,L \cM)$ is a categorical equivalence.
\end{thm}

The proof will rely on the inductive description of sections over Reedy categories.
For $\Sect(\cR,L \cM)$, we have Proposition \ref{propinfinityreedyinduction} and,
evidently, the fact that for a good filtration $\cR_\beta$, we have that
$\Sect(\cR,L \cM) \cong \lim \Sect(\cR_\beta,\cM)$, where the limit is taken in
$\Cat_\infty$ along the filtration of $\cR$. Let us see if the similar description
holds for $L \Sect(\cR,\cM)$.

\smallskip

Given any model category $\cX$, denote by $\cX_{c \to f} \subset \Fun([1],\cX)$
the full subcategory given by $X \to Y$ with $X$ cofibrant and $Y$ fibrant.
Similarly, denote $\cX_{c \hookrightarrow cf \twoheadrightarrow f} \subset \Fun([2],\cX)$
the full subcategory of $X \to Y \to Z$ with $X$ cofibrant, $Z$ fibrant, $X \to Y$
a cofibration, $Y \to Z$ a fibration. These categories appear naturally as targets
for latching-matching decompositions.

\begin{prop}
\label{propcofibranttofibrant}
Let $\cX$ be a model category.

\begin{enumerate}
\item The inclusions $\cX_{c \to f} \subset \Fun([1],\cX)$ and $\cX_{c \hookrightarrow cf \twoheadrightarrow f} \subset \Fun([2],\cX)$ induce weak equivalences of infinity-localisations.

\item The functor $\cX_{c \hookrightarrow cf \twoheadrightarrow f} \to \cX_{c \to f}$
sending $X \to Y \to Z$ to $X \to Z$ induces an isofibration after taking (ordinary)
localisations, $\Ho \cX_{c \hookrightarrow cf \twoheadrightarrow f} \to \Ho \cX_{c \to f}$.
\end{enumerate}
\end{prop}
\proof \newline
1. We will do the case of $\cX_{c \to f} \subset \Fun([1],\cX)$, with the second
case being similar but more cumbersome. Denote by $\cY$ the category which objects are
diagrams
\begin{equation}
\label{diadefinitionoftriples}
\begin{diagram}[small]
X_1 					& 	\rTo 	& 		Y_1 \\
\uTo<\sim	&					&		\uTo>\sim \\
QX_2					& \rTo 		& 		Y_2 \\
\dTo^\sim			&					&		\dTo>\sim \\
QX_3					& \rTo 		& 		RY_3 \\
\end{diagram}
\end{equation}
where the vertical maps are weak equivalences, the objects whose name starts with
$Q$ are cofibrant, and the objects whose name starts with $R$ are fibrant. Various
projections from $\cY$ define a diagram of weak equivalence preserving functors
fitting together in a 2-diagram
\begin{diagram}[small,nohug]
							&				&  \cY 					&					&						\\
							&	\ldTo>\Leftarrow	&		\dTo				&	\rdTo<\Rightarrow		&							\\
\cX_{c \to f} & \lTo 	& \cX_{c \to *} & \rTo & \cX_{* \to *}.		\\
\end{diagram}
Here, as the names suggest, $\cX_{* \to *} = \Fun([1],\cX)$ and $\cX_{c \to *}$ is
its full subcategory of arrows with cofibrant domain. Both depicted natural transformations
are valued in $\cW$. Moreover, the inclusion $\cX_{c \to f} \subset \cX_{* \to *}$
can be factored as $\cX_{c \to f} \to \cY \to \cX_{* \to *}$. Here, the first functor
is a homotopy inverse (in the sense of \cite{BARKAN0}) to the depicted $\cY \to \cX_{c \to f}$,
given by putting $QX \to RY$ in all three rows of (\ref{diadefinitionoftriples}).

In general, given two localisers $(\cM,\cW_\cM),(\cN,\cW_\cN)$ one has the following diagram
of (higher) categories
$$
\Fun'(\cM,\cN) \longrightarrow \Fun'(\cM, L\cN) \stackrel \sim \longleftarrow \Fun(L \cM, L \cN)
$$
with primes indicating (weak) equivalence preserving functors. As a consequence, a
natural transformation $\alpha: F \Rightarrow G$ of weak equivalence preserving functors induces a natural transformation
$\bar \alpha: \bar F \to \bar G$
of infinity-functors from $L \cM$ to $L \cN$; if $\alpha$ took values in $\cW_\cN$
then $\bar \alpha$ is invertible. Thus everything
will follow if we prove that both functors
$$
\cX_{c \to f} \leftarrow \cX_{c \to *} \rightarrow \cX_{* \to *}
$$
are weak equivalences in the model structure of \cite{BARKAN0}. Let us first consider the
case of $\cX_{c \to *} \rightarrow \cX_{* \to *}$.
We shall prove that it is a left approximation. Given $X_\bullet \to Y_\bullet$
in $\Fun([n],\cX_{* \to *})$, we have to prove
the contractibility of the category $\cC(X_\bullet \to Y_\bullet)$ of diagrams
\begin{diagram}[small]
QX_\bullet					& \rTo 		& 		Y'_\bullet \\
\dTo^\sim			&					&		\dTo>\sim \\
X_\bullet			& \rTo 		& 		Y_\bullet \\
\end{diagram}
with bottom row fixed. There is a functor $\cC(X_\bullet \to Y_\bullet) \to \cC(X_\bullet)$,
where $\cC(X_\bullet)$ is the category of objectwise cofibrant replacements of $X_\bullet$,
and hence it is contractible. The functor $\cC(X_\bullet \to Y_\bullet) \to \cC(X_\bullet)$
is furthermore a Grothendieck fibration, so it will suffice to show that its fibres
have contractible nerve. The fibres are seen to have final objects, given by diagrams
\begin{diagram}[small]
QX_\bullet					& \rTo 		& 		Y_\bullet \\
\dTo^\sim			&					&		\dTo>= \\
X_\bullet			& \rTo 		& 		Y_\bullet. \\
\end{diagram}

In the case of $\cX_{c \to f} \leftarrow \cX_{c \to *}$, we are, using the right
approximation argument, presented with studying the category $\cD(QX_\bullet	\to Y_\bullet)$ of
diagrams
\begin{equation}
\label{diacoftof}
\begin{diagram}[small]
QX_\bullet					& \rTo 		& 		Y_\bullet \\
\dTo^\sim			&					&		\dTo>\sim \\
QX'_\bullet			& \rTo 		& 		RY_\bullet \\
\end{diagram}
\end{equation}
with fixed first row, and we instead consider the projection $\cD(QX_\bullet	\to Y_\bullet)
\to \cD(Y_\bullet)$ to the category of fibrant replacements of $Y_\bullet$. The fibres
of this projection are given by diagrams (\ref{diacoftof}) with top row and left column fixed,
and again,
\begin{diagram}[small]
QX_\bullet					& \rTo 		& 		Y_\bullet \\
\dTo^=			&					&		\dTo>\sim \\
QX_\bullet			& \rTo 		& 		RY_\bullet \\
\end{diagram}
serves as an initial object.

2. Consider first the following situation:
\begin{diagram}[small]
X 		& \rInto & Y & \rOnto & Z \\
\dTo^\sim &			&		&				&		\dTo>\sim \\
X'		&					&	\rTo	&			& Z' \\
\end{diagram}
with the top row in $\cX_{c \hookrightarrow cf \twoheadrightarrow f}$ and the bottom
row in $\cX_{c \to f}$. Taking the pushout of the left corner, we factor the bottom
map as follows:
\begin{diagram}[small]
X 		& \rInto & Y & \rOnto & Z \\
\dTo^\sim &			& \dTo>\sim		&				&		\dTo>\sim \\
X'		&		\rInto	& Y''		&	\rTo		& Z'. \\
\end{diagram}
Since $X, X', Y$ are cofibrant and $X \to Y$ is a cofibration, we see that
the induced map $Y \to Y''$ is a weak equivalence. Factoring $Y'' \to Z'$ as a
trivial cofibration $Y'' \hookrightarrow Y'$ followed by a fibration $Y' \twoheadrightarrow Z'$,
we get a modified factorisation of the map $X' \to Z'$, together with the diagram
\begin{diagram}[small]
X 		& \rInto & Y & \rOnto & Z \\
\dTo^\sim &			& \dTo>\sim		&				&		\dTo>\sim \\
X'		&		\rInto	& Y'		&	\rOnto		& Z'  \\
\end{diagram}
in $\cX_{c \hookrightarrow cf \twoheadrightarrow f}$.

Using a dual argument with pullbacks, one can treat the case of weak equivalences pointing up:
\begin{diagram}[small]
X 		& \rInto & Y & \rOnto & Z \\
\uTo^\sim &			&		&				&		\uTo>\sim \\
X'		&					&	\rTo	&			& Z'. \\
\end{diagram}
In general, a weak equivalence in $\Ho \cX_{c \to f}$ can be represented as a
zig-zag of weak equivalences in $\cX_{c \to f}$. Using the above arguments,
the isofibration property of $\Ho \cX_{c \hookrightarrow cf \twoheadrightarrow f} \to \Ho \cX_{c \to f}$
is readily verified step by step. \endproof

\begin{rem}
The main difficulty of Proposition \ref{propcofibranttofibrant}.1 comes from
the non-assumption of the functoriality of factorisations. If the factorisations
in $\cX$ are functorial, it is easy to construct the (weakly) inverse functors to the inclusions $\cX_{c \to f} \subset \Fun([1],\cX)$ and $\cX_{c \hookrightarrow cf \twoheadrightarrow f} \subset \Fun([2],\cX)$.
\end{rem}

Denote by $ \Sect(\cR,\cM)_{\cf}$
the full subcategory of cofibrant-fibrant sections for the Reedy model structure.
Since $L \Sect(\cR,\cM)_{\cf} \cong L \Sect(\cR,\cM)$, Theorem \ref{thmcomparison}
will follow if we prove that the induced functor $L \Sect(\cR,\cM)_{\cf} \to \Sect(\cR,L \cM)$
is a categorical equivalence.

\begin{prop}
\label{propinductionformodelcategoriesofsections}
Let $\cM \to \cR$ be a left model Reedy fibration. Choose a
good filtration $\{ \cR_\beta \}_{\beta < \alpha}$. Then
\begin{enumerate}
\item if $\cR_\beta$ is obtained from $\cR_{< \beta}$ by adding $x \in \cR$, then
the assignments
$$X \mapsto (\Lat_x X \to X(x) \to \Mat_x X) \text{  and  } X \mapsto (\Lat_x X \to \Mat_x X)$$
determine a pullback square in $\Cat_\infty$
\begin{diagram}[small]
L \Sect(\cR_\beta, \cM)_{\cf} & \rTo & L \Fun([2],\cM(x)) \\
\dTo 									&				& \dTo							\\
L \Sect(\cR_{< \beta}, \cM)_{\cf} & \rTo		&	L \Fun([1],\cM(x)).		\\
\end{diagram}
\item For each ordinal $\beta$, the  $\Cat_\infty$-limit of $\{ L \Sect(\cR_{\gamma},\cM)_{cf}\}_{\gamma < \beta} $ is equivalent, via the evident map, to $L \Sect(\cR_{ < \beta},\cM)_{cf}$.
\item The $\Cat_\infty$-limit of $\{L \Sect(\cR_{\beta},\cM)_{cf}\}_{\beta}$ is equivalent,
via the evident map, to $L\Sect(\cR,\cM)_{cf}$.
\end{enumerate}
\end{prop}

\proof Considering the relation between quasicategorical and Dwyer-Kan simplicial localisation
as explained in \cite{HINICH}, it will be sufficient to prove the corresponding statements
for simplicial localisations. For this, we shall use Lemmas \ref{lmrecognitionhomotopypullbacks}
and \ref{lmrecognitionhomotopyinvlims}.

Using the same notations as before, there is a pullback square of categories
\begin{equation}
\label{diastartingpbck}
\begin{diagram}[small]
 \Sect(\cR_\beta, \cM)_{\cf} & \rTo &  \cM(x)_{c \hookrightarrow cf \twoheadrightarrow f} \\
\dTo 									&				& \dTo							\\
 \Sect(\cR_{< \beta}, \cM)_{\cf} & \rTo		&	 \cM(x)_{c \to f}		\\
\end{diagram}
\end{equation}
with horizontal functors given by latching-matching factorisations. This square has the
property that all functors preserve weak equivalences.

Take a cofibrant-fibrant section $X \in  \Sect(\cR_\beta, \cM)_{\cf}$ and a(n always existing)
simplicial resolution $\mathbf Y \in \Sect(\cR_\beta, \cM)^{\Delta^\op}_{\cf}$.
Note that the restriction $X|_{\cR_{< \beta}}$ remains cofibrant-fibrant and $\mathbf Y |_{\cR_{< \beta}}$
remains a resolution valued in cofibrant-fibrant objects, since the condition for
$\mathbf Y$ to be Reedy fibrant (as a simplicial object) is required objectwise. We thus have
that the natural map
$$
 \Sect(\cR_\beta, \cM)(X,\mathbf Y ) \to   \Sect(\cR_{< \beta}, \cM)(X|_{\cR_{< \beta}},\mathbf Y |_{\cR_{< \beta}})
$$
is a map of Kan complexes representing the corresponding map $L^H \Sect(\cR_\beta, \cM)(X,\mathbf Y(0)) \to L^H  \Sect(\cR_{< \beta}, \cM)(X|_{\cR_{< \beta}},\mathbf Y (0)|_{\cR_{< \beta}})$ (we omit the subscript $(-)_{cf}$
as no information is lost).

\smallskip
Put the projective model structure on $\Fun([1],\cM(x))$ (the unique map of $[1]$ raises
the degree). For $\Fun([2],\cM(x))$, consider
the Reedy structure such that $0 \to 1$ and $0 \to 2$ are
treated as degree-raising morphisms, and $1 \to 2$ as a degree-lowering. The vertical functor
$\Fun([2],\cM(x)) \to \Fun([1],\cM(x))$ then preserves all classes of the model structure.
The functor
$$
\Sect(\cR_\beta, \cM)_{\cf} \rTo \Fun([2],\cM(x)), \, \, \, X \mapsto X_{012} = \Lat_x X \to X(x) \to \Mat_x X,
$$
(we shall henceforth write $X_0 = \Lat_x X$ and so on) however, does not send fibrant-cofibrant objects to fibrant-cofibrant objects.
Same observation applies to the functor
$$
\Sect(\cR_{< \beta}, \cM)_{\cf} \rTo \Fun([1],\cM(x)), \, \, \, X \mapsto X_{02} = \Lat_x X \to \Mat_x X.
$$
Let us see how to correct this. For cofibrant-fibrant $X$, its value $X_{012} = X_0 \to X_1 \to X_2$
has the property that $X_0$ is cofibrant and $X_0 \to X_1$ is a cofibration. For our
purposes, it will be sufficient to work with a cofibrant replacement $Q X_{012}$
that we obtain by factoring $X_0 \to X_2$ as a cofibration $X_0 \to Q X_2$ followed
by a trivial fibration $Q X_2 \to X_2$. Since $X_0 \to X_1$ is a cofibration,
the dotted arrow exists in the diagram
\begin{diagram}[small]
		&				&	Q X_2	&			&	\\
		&	\ruTo &		\uDashto			&	\rdTo	&		\\
X_0 & \rTo & X_1 				& \rTo & X_2 \\
\end{diagram}
thus defining $Q X_{012} \to X_{012}$. Note that the map resulting from composing
these sequences, $Q X_{02} \to X_{02}$, is a cofibrant replacement as well.

Similarly, taking a simplicial resolution $\mathbf Y$, we would like to arrange
for its Reedy fibrant replacement. Since the transition functors of $\cM \to \cR$
along the ``matching'' maps $\cR_-$ are right adjoints, and since the limits in
$\Sect(-,\cM)$ are calculated fibre by fibre, we have the commutativity of $\Mat_x$
and the simplicial matching functor $\Mat^\Delta_n$. In particular, we have that
the assignment
$$
\Sect(\cR_\beta, \cM)^{\Delta^\op}_{\cf} \rTo \Fun([2],\cM(x))^{\Delta^\op}, \, \, \, \mathbf Y  \mapsto \mathbf Y _{012}  = \Lat_x \mathbf Y  \to \mathbf Y (x) \to \Mat_x \mathbf Y ,
$$
induces for each $[n] \in \Delta$ the diagram
\begin{diagram}[small]
\mathbf Y _0(n) & \rTo & \mathbf Y _1(n) & \rTo & \mathbf Y_2(n) \\
\dTo			&			&		\dTo	&		&	\dTo	\\
\Mat^\Delta_n \mathbf Y_0  & \rTo & \Mat^\Delta_n \mathbf Y_1  & \rTo & \Mat^\Delta_n \mathbf Y_2  \\
\end{diagram}
such that the maps
$ \mathbf Y_2(n) \to \Mat^\Delta_n \mathbf Y_2$ and $\mathbf Y_1(n) \to  \mathbf Y_2(n) \prod_{\Mat^\Delta_n \mathbf Y_2} \Mat^\Delta_n \mathbf Y_1$
(and hence $ \mathbf Y_1(n) \to \Mat^\Delta_n \mathbf Y_1$) are fibrations.
As a consequence, to get a fibrant replacement
of $\mathbf Y_{012}$ it will suffice to factor $\mathbf Y_0 \to \mathbf Y_1$ as
a trivial (simplicial) Reedy cofibration
$\mathbf Y_0 \to R{\mathbf Y}_1$ followed by a Reedy fibration
$R \mathbf Y_1 \to \mathbf Y_1$. Denote the result
by $\mathbf Y_{012} \to R {\mathbf Y}_{012}$;
the associated map $\mathbf Y_{02} \to R {\mathbf Y}_{02}$ is also
a Reedy fibrant replacement in $\Fun([1],\cM(x))^{\Delta^\op}$.

Observe now that, starting from  $X \in  \Sect(\cR_\beta, \cM)_{\cf}$ and $\mathbf Y \in \Sect(\cR_\beta, \cM)^{\Delta^\op}_{\cf}$ as above, we have the following pullback squares of simplicial sets:
\begin{equation}
\label{diasimpsonlemma}
\begin{diagram}[small]
 \Sect(\cR_\beta, \cM)(X,\mathbf Y) & \rTo &  \Sect(\cR_{< \beta}, \cM)(X|_{\cR_{< \beta}},\mathbf Y|_{\cR_{< \beta}}) \\
\dTo 									&				& \dTo							\\
  \Fun([2],\cM(x))(X_{012},\mathbf Y_{012})  & \rTo		&	 \Fun([1],\cM(x))(X_{02},\mathbf Y_{02}) 		\\
	\dTo 									&				& \dTo							\\
	  \Fun([2],\cM(x))(QX_{012},\mathbf Y_{012})  & \rTo		&	 \Fun([1],\cM(x))(QX_{02},\mathbf Y_{02}) 		\\
		\dTo 									&				& \dTo							\\
		  \Fun([2],\cM(x))(QX_{012},R\mathbf Y_{012})  & \rTo		&	 \Fun([1],\cM(x))(QX_{02},R\mathbf Y_{02}).		\\
\end{diagram}
\end{equation}
We thus have that the outer square is likewise a pullback. Moreover, observe that
$$\Fun([2],\cM(x))(QX_{012},R\mathbf Y_{012}) \longrightarrow \Fun([1],\cM(x))(QX_{02},R\mathbf Y_{02})$$
is a Kan fibration: it is induced by applying $\Fun([2],\cM(x))(-,R\mathbf Y_{012})$
to the map
\begin{diagram}[small]
X_0 & \rTo & X_0 & \rTo & QX_2 \\
\dTo	&			&	\dTo	&		&		\dTo \\
X_0 & \rTo & X_1 & \rTo & QX_2 \\
\end{diagram}
which is a cofibration in the chosen model structure on $\Fun([2],\cM(x))$.
This implies that
\begin{equation}
\label{requiredmapisakanfib}
 \Sect(\cR_\beta, \cM)(X,\mathbf Y)\longrightarrow  \Sect(\cR_{< \beta}, \cM)(X|_{\cR_{< \beta}},\mathbf Y|_{\cR_{< \beta}})
\end{equation}
is a Kan fibration and that the outer square of (\ref{diasimpsonlemma}) is a homotopy
pullback.

Let us
examine the induced pull-back square of simplicial localisations,
\begin{equation}
\label{squaredk}
\begin{diagram}[small]
L^H \Sect(\cR_\beta, \cM)_{\cf} & \rTo &  L^H \cM(x)_{c \hookrightarrow cf \twoheadrightarrow f} \\
\dTo 									&				& \dTo							\\
L^H \Sect(\cR_{< \beta}, \cM)_{\cf} & \rTo		&	 L^H  \cM(x)_{c \to f}.		\\
\end{diagram}
\end{equation}
The right vertical functor is a $\Ho$-isofibration by (2.) of Proposition \ref{propcofibranttofibrant},
and for each $X,Y \in L^H \Sect(\cR_\beta, \cM)_{\cf}$,
the induced square of simplicial hom-spaces is a homotopy pullback, as follows from the
argument above and \cite[Lemma 6.1]{MANDELL}.
Lemma \ref{lmrecognitionhomotopypullbacks} thus tells us that (\ref{squaredk})
is a homotopy pullback. Combining with the equivalences
$$
L^H \cM(x)_{c \hookrightarrow cf \twoheadrightarrow f} \stackrel \sim \to
L^H \Fun([2],\cM(x)) \text{ and } L^H  \cM(x)_{c \to f} \stackrel \sim \to L^H \Fun([1],\cM(x)),
$$
we prove the first statement of the proposition.

The $1$-category $\Sect(\cR,\cM)_{cf}$ is equivalent (even isomorphic) to the limit of
the inverse system $\{\Sect(\cR_{\beta},\cM)_{cf}\}_{\beta < \alpha}$.
Let us verify that the conditions of Lemma \ref{lmrecognitionhomotopyinvlims} are satisfied
for $$L^H \Sect(\cR,\cM)_{cf} \to L^H \Sect(\cR_\beta,\cM)_{cf} \to L^H \Sect(\cR_{<\beta},\cM)_{cf}.$$
The conditions (0.), (2.) and (i.) of
Lemma \ref{lmrecognitionhomotopyinvlims} are readily verified.

Using the pullback diagram (\ref{diastartingpbck}) and (2.) of Proposition \ref{propcofibranttofibrant},
we conclude that for each $\beta$,
the functor $L^H \Sect(\cR_{\beta},\cM)_{cf} \to L^H \Sect(\cR_{< \beta},\cM)_{cf}$
is a $\Ho$-isofibration. This corresponds to the condition (1.) of Lemma \ref{lmrecognitionhomotopyinvlims}.
The limit ordinals and $\cR$ itself are treated by the same sort of argument.
If $\beta$ is an ordinal, then $\Sect(\cR_{< \beta}, \cM)_{cf}(X,\mathbf Y)$
coincides with the limit of $\{\Sect(\cR_\gamma, \cM)_{cf}(X|_{\cR_\gamma},\mathbf Y|_{\cR_\gamma}) \}_{\gamma < \beta}$.
Using the Kan fibration remark around (\ref{requiredmapisakanfib}) and induction, we observe that
$\Sect(\cR_{< \beta}, \cM)_{cf}(X,\mathbf Y)$ is also the homotopy limit of its restrictions.
Inductively we also get that $\Sect(\cR, \cM)_{cf}(X,\mathbf Y)$ is the homotopy limit of its
restrictions, and using \cite[Lemma 6.1]{MANDELL} again, we verify (3.) and (ii.).
\endproof

\proof[Proof of Theorem \ref{thmcomparison}.] In light of Propositions \ref{propinfinityreedyinduction} and
\ref{propinductionformodelcategoriesofsections}, everything will follow by transfinite induction from the
comparison for model categories as described in Proposition \ref{propcomparisonforasimplex}, if we show
that for each $\beta$, there is a map in $\Cat_\infty$,
naturally induced on each term by $F:\cM \to L \cM$, from the
square
\begin{equation}
\label{diacomp1}
\begin{diagram}[small]
 \Sect(\cR_\beta, \cM)_{\cf} & \rTo &  \Fun([2],\cM(x)) \\
\dTo 									&				& \dTo							\\
 \Sect(\cR_{< \beta}, \cM)_{\cf} & \rTo		&	 \Fun([1],\cM(x)) 	\\
\end{diagram}
\end{equation}
to the square
\begin{equation}
\label{diacomp2}
\begin{diagram}[small]
\Sect(\cR_\beta, L\cM)  & \rTo &   \Fun([2],L\cM(x)) \\
\dTo 									&				& \dTo							\\
\Sect(\cR_{< \beta}, L\cM)  & \rTo		&	  \Fun([1],L\cM(x)). 		\\
\end{diagram}
\end{equation}
For this, we observe that cofibrant-fibrant
sections are $(F,x)$-compatible in the sense of Definition \ref{oprfxcompsect}.
Indeed, the functor $\cM \to L \cM$ preserves coCartesian arrows with cofibrant domain
and Cartesian arrows over $\cR_-$ with fibrant codomain; given a cofibrant-fibrant
section $S$, its right restriction $R_x S: Mat(x) \to \cM(x)$
is an injectively fibrant diagram, and hence its limit remains \cite[Proposition 2.5.6]{HARPAZ2} a limit in $L \cM(x)$.
A dual observation is true for the left restrictions. Therefore, the square (\ref{diacomp1})
factors as
\begin{diagram}[small]
 \Sect(\cR_\beta, \cM)_{\cf} & \rTo &   \Sect(\cR_\beta, \cM)_{(F,x)} & \rTo &   \Fun([2],\cM(x)) \\
\dTo 									&				& \dTo				&				& \dTo						\\
 \Sect(\cR_{< \beta}, \cM)_{\cf} & \rTo		&	 \Sect(\cR_{< \beta}, \cM)_{(F,x)} & \rTo		&	 \Fun([1],\cM(x)).	\\
\end{diagram}
Using Proposition \ref{propfunctorialityofreedyinduction}, we get the desired map.
\endproof

\begin{corr}
\label{cormodelcatbicomplete}
Let $\cM$ be a model category. Then the infinity-category $L \cM$ is bicomplete.
\end{corr}
\proof
By duality, it is enough to prove the colimit part.
Theorem \ref{thmcomparison} establishes an equivalence $L \Fun(\cR,\cM) \cong \Fun(\cR, L \cM)$
for any Reedy category $\cR$. Moreover, if $\cR$ is direct, then we have a Quillen
adjunction
$$
\colim_\cR: \Fun(\cR,\cM) \rightleftarrows \cM: c_\cR
$$
that by Corollary \ref{corollaryadjunctionsgiveadjunctions} gives an infinity-adjunction
$$
\bL\colim_\cR: \Fun(\cR,L\cM) \rightleftarrows L \cM: \bR c_\cR.
$$
We thus conclude that $L \cM$ has colimits of shape $\cR$ for each direct Reedy
category $\cR$. But in particular that means the existence of pushouts and arbitrary
coproducts, and by \cite[Proposition 4.4.2.6]{LUHTT} $L \cM$ thus has all small colimits.
\endproof

\subsection{On Quillen presheaves}

We conclude this paper by applying Theorem \ref{thmcomparison} to the presheaves
of model categories and Quillen adjunctions.

It is known \cite[Theorem 7.9.8]{CISBOOK} that for a model category $\cM$ and
a small category $\cC$, one
has the equivalence $L \Fun( \cC , \cM) \to \Fun(\cC, L \cM)$, even if
$\Fun( \cC , \cM)$ has no obvious model structure if $\cM$ is not cofibrantly
generated. The theory of Quillen presheaves is rich in structure, so it is perhaps
of no surprise that the same observation applies in this setting:

\begin{prop}(Cf \cite[Conjecture 18.3]{H-S})
\label{propquillenpresheafstrict}
Let $\cM \to \cC$ be a Quillen presheaf over a small category $\cC$. Then localising $\cM$ along the
fibrewise weak equivalences yields an equivalence of infinity-categories
$$
L\Sect(\cC,\cM) \stackrel \sim \longrightarrow \Sect(\cC, L \cM).
$$
In other words, sections of Quillen presheaves can be strictified over an
arbitrary 1-categorical base.
\end{prop}
It will be convenient to prove this proposition by the end of this subsection.
\begin{rem}
The proof of Proposition \ref{propquillenpresheafstrict} and of all subsequent
results depends only on one adjoint. Thus everything can be generalised to the
case of an opfibration in model categories $\cM \to \cC$ which transition
functors preserve cofibrations, trivial cofibrations, and colimits. We have
chosen to work with Quillen presheaves as almost all the examples of
colimit-preserving functors between model categories admit right adjoints.
\end{rem}

Given a Quillen presheaf $\cM \to \cC$, consider a map $f: c \to c'$ and the
induced adjunction
$f_!:\cM(c) \rightleftarrows \cM(c'):f^*.$
Corollary \ref{corollaryadjunctionsgiveadjunctions} supplies us with the infinity-adjunction
$\bL f_!:L\cM(c) \rightleftarrows L\cM(c'):\bR f^*,$
and we shall use the same notation to denote the induced adjunction on the homotopy
categories:
$\bL f_!:\Ho \cM(c) \rightleftarrows \Ho \cM(c'):\bR f^*.$

\begin{lemma}
\label{lemmacartcocartsectionsofaqpresheaf}
Let $\cM \to \cC$ be a Quillen presheaf, and $S_l,S_r$ two sets of arrows of $\cC$.
Then the equivalence of Proposition \ref{propquillenpresheafstrict} induces an equivalence
$$
L\Sect_{(S_l,S_r)}(\cC,\cM) \stackrel \sim \longrightarrow \Sect_{(S_l,S_r)}(\cC, L \cM)
$$
to the infinity category $\Sect_{(S_l,S_r)}(\cC, L \cM)$ of $S_l$-coCartesian and
$S_r$-Cartesian sections, from the localisation of the full subcategory
$\Sect_{(S_l,S_r)}(\cC,\cM) \subset \Sect(\cC,\cM)$ consisting of all sections
$S: \cC \to \cM$ such that
\begin{itemize}
\item[i.] For each $f:c \to c' \in S_l$, the induced map $\bL f_! S(c) \to S(c')$ is an
isomorphism in $\Ho \cM(c')$,
\item[ii.] For each $f:c \to c' \in S_r$, the induced map $S(c) \to \bR f^* S(c')$ is an
isomorphism in $\Ho \cM(c)$.
\end{itemize}
\end{lemma}
\proof Follows directly from the observations of Corollaries \ref{corollarylocalisationinfamilies} and
\ref{corollaryadjunctionsgiveadjunctions}. \endproof

\begin{corr}
Let $\cM \to \cC$ be a Quillen presheaf. Then
the higher-categorical limit of the covariant (respectively contravariant)
diagram $x \mapsto L \cM(x)$ is given by
the infinity-localisation $L \Sect_{(\cC,\emptyset)}(\cC,\cM)$, (respectively
$L \Sect_{(\emptyset,\cC)}(\cC,\cM)$)
consisting of those sections $S$ such that, given any $f:c \to c'$,
the induced map $\bL f_! S(c) \to S(c')$ (respectively  $S(c) \to \bR f^* S(c')$) is an
isomorphism in the homotopy category.
\end{corr}
\proof Direct consequence of \cite[Corollary 3.3.3.2]{LUHTT}, which is itself
a consequence of the naturality of the Grothendieck construction \cite[Remark 3.1.13]{MGEE1}.
\endproof

Quillen presheaves are closely related to the notion of descent. For an exemplary statement,
let $B^\bullet: \Delta \to \cC$ be a cosimplicial diagram in $\cC$, and
$A \to B^\bullet$ a natural transformation from the (constant diagram given by)
$A \in \cC$. For a Quillen presheaf $\cM \to \cC$, denote by $\Sect(B^\bullet,\cM)=\Sect(\Delta,(B^\bullet)^*\cM)$.

\begin{opr}
Let $\cM \to \cC$ be a Quillen presheaf. A morphism $A \to B^\bullet$ in the notation above
satisfies the descent property with respect to $\cM \to \cC$ if the induced functor
$$\cM(A) \to \Sect_{\bL cart}(B^\bullet,\cM)$$
is a weak equivalence of relative categories, where $\Sect_{\bL cart}(B^\bullet,\cM)$
is the sub-category of sections $X:\Delta \to \cM$ such that for each $\alpha:[n] \to [m]$,
the induced map $\mathbb L \alpha_! X(n) \to X(m)$ is an equivalence.
\end{opr}

\begin{corr}
In the situation above, if $A \to B^\bullet$
satisfies the descent property with respect to $\cM \to \cC$,
then the map $A \to B^\bullet$ exhibits $L\cM(A)$ as a $\Cat_\infty$-limit of $L \cM(B^\bullet)$.
\end{corr}
\proof Immediate. \endproof

\begin{rem}
\label{remdescent}
Observe that the functor $\cM(A) \to \Sect_{\bL cart}(B^\bullet,\cM)$ factors as
$$\cM(A) \to \Fun(\Delta,\cM(A)) \to \Sect_{\bL cart}(B^\bullet,\cM)$$
where the first functor is the constant diagram inclusion, and the second one
is induced by applying $\Sect(-,\cM)$ to the map of the cosimplicial diagrams
$f:A \to B^\bullet$. By Theorem \ref{thmcomparison} and Corollary \ref{cormodelcatbicomplete}, there is an adjunction
$$
\bL const:L\cM(A) \rightleftarrows  \Fun(\Delta,L\cM(A)) \cong L\Fun(\Delta,\cM(A)): \bR \lim,
$$
even though $\bR \lim$ does not in general come from a Quillen functor.
We also have an induced infinity-adjunction
$$
\bL f_!: L\Fun(\Delta,\cM(A)) \rightleftarrows L\Sect(B^\bullet,\cM): \bR f^*
$$
where $f_! \dashv f^*$ is the Quillen adjunction induced by $f:A \to B^\bullet$. Thus to check
that  $L\cM(A) \to L\Sect_{\bL cart}(B^\bullet,\cM)$ is an equivalence, it is enough
to check that the infinity-adjunction $\bL f_! \bL const \dashv \bR \lim \bR f^*$
is an adjoint equivalence. This, however, can be done on the level of homotopy categories.
Thus the already known theory of descent for Quillen presheaves
(see e.g. \cite[Lemma 2.2.2.13]{HAG2}) generalises to our setting.
\end{rem}

There are many examples of Quillen presheaves that live over a base which is
equipped with a subcategory of weak equivalences. Following \cite{HAPRASMA}, we introduce
\begin{opr}
A \emph{relative Quillen
presheaf} over a category with weak equivalences $(\cC,\cW_{\cC})$ is a Quillen presheaf
$\cM \to \cC$ such that for each $f: c \to c'$ in $\cW_{\cC}$ the induced Quillen adjunction
$f_!:\cM(c) \rightleftarrows \cM(c'):f^*$ is a Quillen equivalence.
\end{opr}

Note that if we take the associated ``straightened'' functor $St(\cM):c \mapsto \cM(c)$ and postcompose with
the localisation, the induced infinity-functor $L St(\cM): \cC \to \Cat_\infty$ sends the
maps of $\cW_{\cC}$ to equivalences of infinity-categories. It thus comes from an
essentially unique functor $\overline{L St(\cM)}: L \cC \to \Cat_\infty$ by applying
pullback along $\cC \to L \cC$.

If $(\cC,\cW_{\cC})$ is small, we can consider the category $\Sect_{\cW_{\cC}}(\cC,\cM)$
consisting of those sections $S$ such that for each $f: c \to c'$ in $\cW_{\cC}$,
both $\bL f_! S(c) \to S(c')$ and $ S(c) \to \bR f^* S(c')$ are isomorphisms in
the homotopy category.

\begin{lemma}
\label{lemmasectionslocalisation}
Given a relative category $(\cC,\cW_\cC)$ and a categorical fibration $\cX \to L \cC$,
the pull-back operation induces a categorical equivalence $\Sect(L\cC,\cX) \cong \Sect_{\cW_\cC}(\cC,\cX)$, with the latter
denoting the subcategory of sections $\cC \to \cX$ which send $\cW_\cC$ to
equivalences of $\cX$.
\end{lemma}
\proof We have a ($\Cat_\infty$-pullback) diagram
\begin{equation}
\label{diasectionsoflocalisation}
\begin{diagram}[small]
\Fun(L\cC,\cX) & \rTo^\sim & \Fun_{\cW_\cC}(\cC,\cX) \\
\dTo 					&						&		\dTo						\\
\Fun(L\cC, L\cC)	&	\rTo^{\sim}& \Fun_{\cW_\cC}(\cC,L \cC)			\\
\end{diagram}
\end{equation}
with $\Fun_{\cW_\cC}(\cC,\cX)$ (and similarly $\Fun_{\cW_\cC}(\cC,L \cC)$)
denoting the infinity-category of functors $\cC \to \cX$
that send $\cW_\cC$ to equivalences of $\cX$. Both vertical maps
in (\ref{diasectionsoflocalisation}) are categorical fibrations, so taking strict $\SSet$-pullbacks
over $id_\cC$ and $\cC \to L \cC$
induces the sought-after equivalence $\Sect(L\cC,\cX) \cong \Sect_{\cW_\cC}(\cC,\cX)$.

\begin{prop}
Let $\cM \to \cC$ be a Quillen presheaf over a small localiser $(\cC,\cW_{\cC})$.
Then the infinity-category $L\Sect_{\cW_{\cC}}(\cC,\cM)$ is naturally equivalent to the
sections
$\Sect(L \cC, \int \overline{L St(\cM)})$ of the unstraightening $\int \overline{L St(\cM)}) \to L \cC$.
\end{prop}
\proof
Localising $\cM$ along the fibrewise weak equivalences, we see that there is a
following diagram of $\Cat_\infty$-pullbacks:
\begin{equation}
\label{diasectionsoflocalisation0}
\begin{diagram}[small]
L \cM & \rTo^\sim & \int L St(\cM) & \rTo &  \int \overline{L St(\cM)}) \\
\dTo 	&						&				\dTo 					&			&								\dTo 							\\
\cC		&		\rTo^=	&		\cC					&	\rTo	&			L \cC. 									\\
\end{diagram}
\end{equation}
The second square is a pullback due to the naturality of the Grothendieck construction.
The first one is the adjoint of the map of $\Cat_\infty$-functors
$St L (\cM) \to  L St(\cM)$; the latter being an equivalence is a reformulation
of the universality of the localisation $L \cM \to \cC$.

Combining Lemma \ref{lemmasectionslocalisation} with (\ref{diasectionsoflocalisation0}), we see that
$\Sect(L \cC, \int \overline{L St(\cM)})$ is canonically identified with
the subcategory of those $S$ in $\Sect(\cC, L \cM)$ that send $\cW_\cC$
to Cartesian (or, equivalently, coCartesian) maps of $L \cM$. The conclusion is
then reached using Lemma \ref{lemmacartcocartsectionsofaqpresheaf}.
\endproof
\begin{rem}
The above proposition shows that $L\Sect_{\cW_{\cC}}(\cC,\cM)$ serves as a strict
model for the lax limit of the $\Cat_\infty$-diagram naturally associated to
the data of a relative Quillen presheaf. The situation with the lax colimit is more
intricate. In effect, one would like to localise $\cM$ along the following class
of weak equivalences $\cW(\cM)$: those maps $\alpha : X \to Y$ such that the corresponding
$\cC$-map $f:x \to y$ belongs to $\cW$, and either $\bL f_! X \to Y$ or
$X \to \bR f^* Y$ is an equivalence. However, even with $\cW_\cC$ saturated,
we cannot a priori guarantee that the described class $\cW(\cM)$ is saturated as well.
Thus one cannot conclude that there is an equivalence between $L_{\cW(\cM)} \cM$ and $\overline{L St(\cM)})$.

The work \cite{HAPRASMA} proves that when $\cC$ is a model category, the localiser
$\cW(\cM)$ is saturated. More generally, it seems that whenever the base admits a certain calculus of
fractions, one can verify the saturation property of $\cW(\cM)$. Conjecturally,
$\cW(\cM)$ should be saturated whenever $\cW_\cC$ is.
\end{rem}

Let us finally prove Proposition \ref{propquillenpresheafstrict}. Our strategy will
be to replace any small category by a suitable Reedy category. Given a category
$\cC$, denote by $\Delta \cC$ its category of simplicies: the objects of $\Delta \cC$
are given by functors $\sigma:[n] \to \cC$, and the morphisms are natural transformations
$[n] \to [n']$ compatible with maps to $\cC$.
The category $\Delta \cC$ is a Reedy category
\cite[22.10]{DHKS}, and
the assignment $\sigma \mapsto \sigma(n)$ defines a functor
$p:\Delta \cC \to \cC$, called $p_t$ in $\cite[22.11]{DHKS}$. The comma-fibres
of this functor are moreover related to the comma-fibres of the original category:
one has $p/c = \Delta(\cC / c)$. We shall henceforth write $\Delta \cC / c$ to denote
both categories. The inclusion functor
$\Delta \cC / c \to \Delta \cC$ is seen to identify the latching categories:
$Lat(\sigma,\sigma(n) \to c)$ computed in $\Delta \cC / c$, is isomorphic to $Lat(\sigma)$
computed in $\Delta \cC$, with the map to $c$ being automatically supplied. This fact
will be useful later on.

Our first step is the proof that $p: \Delta \cC \to \cC$ presents $\cC$ as the higher-categorical
localisation of $\Delta \cC$ along the subset $\cW$ consisting of those maps that are sent to identities by $p$.
The proof that we give below is the corrected version \cite[Proposition A.1]{TV}: there, the authors
in particular state that $p$ is a coCartesian fibration, and such a claim is false (unfortunately,
we committed the same mistake in an earlier version of this paper). For a different write-up,
we invite the reader to consult \cite[5.3]{STEVENSON}.
\begin{lemma}
\label{lemmacriterionlocalisation}
The functor $p: \Delta \cC \to \cC$  is an infinity-localisation along the $p$-identities $\cW$,
meaning that for any infinity-category $\cX$, the infinity-functor $p^*:\Fun(\cC,\cX) \to \Fun(\Delta \cC,\cX)$ is full and faithful and its essential image consists of all those functors $F: \Delta \cC \to \cX$ that send $\cW$ to equivalences in $\cX$.
\end{lemma}
\proof
Factoring $p$ as $\Delta \cC \stackrel l \longrightarrow L_\cW \Delta \cC \stackrel \pi \longrightarrow \cC$
with $L_\cW \Delta \cC$ being the infinity-localisation, it remains to prove that $\pi$ is
an equivalence of (infinity)-categories. Applying the homotopy category functor to this
diagram would give the usual factorisation $\Delta \cC \to \Ho \Delta \cC \to \cC$
through the one-categorical localisation of $\Delta \cC$; this observation permits to
conclude that $\pi$ is essentially surjective.

To show that $\pi$ is fully faithful, one can use Yoneda lemma which tells us that it is enough
\cite[Proposition 4.5.2 (iii)]{CISBOOK} to verify that the induced left Kan extension
$\pi_!: \Fun(L_\cW \Delta \cC, \cS) \to  \Fun( \cC, \cS)$ is fully faithful, with $\cS$
denoting the infinity-category of spaces. We will show that $\pi_!$ is part of an adjoint equivalence.

We observe the existence of the following diagram,
\begin{diagram}[small,nohug]
                          &     & \Fun_\cW (\Delta \cC, \cS)    &              &           \\
                        & \ruTo<{l^*}    &                      & \rdTo>{p^*}  &           \\
\Fun(L_\cW \Delta \cC, \cS) &     &           \lTo_{\pi^*}                   &              &  \Fun( \cC, \cS),  \\
\end{diagram}
with $\Fun_\cW (\Delta \cC, \cS)$ meaning those infinity-functors that send $\cW$ to
equivalences. The functor $l^*$ is an equivalence by definition. The functor $\pi^*$
is right adjoint to $\pi_!$; it will thus suffice to show that $p^*$ is an equivalence.

The left Kan extension along $p$ restricts to $\Fun_\cW (\Delta \cC, \cS)$, giving
an infinity-adjunction
$$
p_!: \Fun_\cW (\Delta \cC, \cS) \rightleftarrows \Fun(\cC,\cS): p^*.
$$
The surjectivity of $p$ on objects implies conservativity of $p^*$. Thus if we prove
that the unit map $id \to p^* p_!$ is an equivalence, the rest will follow from the
triangle
identities (\cite[Theorem 6.1.23]{CISBOOK}; alternatively one can pass to the homotopy adjunction
and use the associated triangle identities). The pointwise expression for Kan extensions,
\cite[Proposition 6.4.9]{CISBOOK}, implies that for each $X \in \Fun_\cW (\Delta \cC, \cS)$,
the unit map evaluated at $\tau: [n] \to \cC$ with $\tau(n) = c$ takes the form
\begin{equation}
\label{unitmapevaluated}
X(\tau) \to \colim_{\Delta \cC/c} X|_{\Delta \cC /c}.
\end{equation}
As per \cite[23.5]{DHKS},
denote by $p^{-1} c$ the subcategory of $\Delta \cC$ consisting of those simplices
$\sigma: [m] \to \cC$ such that $\sigma(m) = c$ and those maps $[m] \to [m']$ that
map $m$ to $m'$. The notation is from \cite{DHKS} and is abusive: the category
$p^{-1} c$ is \emph{not} the fibre of $p$ at $c$. There is however an obvious functor
$i:p^{-1} c \to \Delta \cC/c$ sending $\sigma$ to $(\sigma, id: \sigma(n)=c)$ and this functor
admits a left adjoint, sending $(\sigma: [n] \to \cC, f: \sigma(n) \to c)$ to the
concatenated simplex $\sigma * f: [n+1] \to \cC$ that is $\sigma$ on first $n+1$ elements and $f$ on the
remaining edge (if we took the naive fibre we would not have an adjunction).
The functor $i$ is thus cofinal. Considering then the composition
$$
[0] \stackrel {\{\tau\}} \longrightarrow p^{-1} c \stackrel i \longrightarrow \Delta \cC/c
$$
and taking the colimit of $X|_{\Delta \cC /c}$ and its pullbacks at each category, we factor
the map (\ref{unitmapevaluated}) as
$$
X(\tau) \to \colim_{p^{-1}c }X|_{p^{-1} c} \stackrel \cong \to \colim_{\Delta \cC/c} X|_{\Delta \cC /c}.
$$
Since $X$ belongs to $\Fun_\cW (\Delta \cC, \cS)$, its restriction to $p^{-1} c$
sends all maps to equivalences in $\cS$. The remaining fact that $X(\tau) \to \colim_{p^{-1}c }X|_{p^{-1} c}$ is an equivalence then follows from Lemma \ref{lemmalocconstfun} below since the category $p^{-1}c$ is
contractible, possessing the initial object given by the zero-simplex $c$. \endproof

\begin{lemma}
\label{lemmalocconstfun}
Let $\cX$ be a cocomplete quasicategory, $K$ a contractible simplicial set and
$x \in K$ a vertex. Then for any $F: K \to \cX$ sending all edges to equivalences,
the natural map $F(x) \to \colim_K F$ is an equivalence.
\end{lemma}
\proof
In $\SSet_+$ factor $K^\flat \to \Delta^0$ as a trivial cofibration $K^\flat \to LK^\natural$
followed by a fibration $LK^\natural \to \Delta^0$. The quasicategory $LK$ is
the localisation of $K$ with respect to all edges, hence $LK$ is a Kan complex.
Moreover the functor $F$ induces an infinity-functor $\bar F: LK \to \cX$.
Both the map $K \to LK$ and $x: [0] \to K \to LK$ are moreover cofinal, as per
\cite[Corollary 4.1.2.6]{LUHTT}. We thus have $\colim_K F \cong \colim_{LK} \bar F$
and $\colim_{LK} \bar F \cong \bar F(x) \cong F(x)$.
\endproof

\proof[Proof of Proposition \ref{propquillenpresheafstrict}.]
By Lemma \ref{lemmacriterionlocalisation}, the functor $p: \Delta \cC \to \cC$ is the localisation along the maps $\cW$
that are sent to identities by $p$. We can assume the validity of Proposition \ref{propquillenpresheafstrict} for
Quillen presheaves over Reedy categories. Given a general Quillen presheaf $\cM \to \cC$, the diagram
\begin{diagram}[small]
p^* L \cM			&	\rTo & L \cM \\
\dTo					&				&	\dTo 	\\
\Delta \cC		&	\rTo &		\cC,  \\
\end{diagram}
together with Lemmas \ref{lemmasectionslocalisation} and \ref{lemmacartcocartsectionsofaqpresheaf}
imply the canonical identification of $\Sect(\cC, L \cM)$ with the infinity-localisation of
the category $\Sect_{ \cW}(\Delta \cC,\cM)$ consisting of those sections
$S:\Delta \cC \to \cM$ that send $\cW$ to (fibrewise) equivalences
of $\cM$. Note that taking pull-backs of sections induces a functor
$$
p^*:\Sect(\cC,  \cM) \to \Sect_{  \cW}(\Delta \cC,\cM)
$$
that preserves pointwise weak equivalences. It remains to show that $p^*$ is a weak equivalence
of relative categories.

Observe that $p^*$ possesses a left adjoint, given by
$$
p_! S(c) = \colim_{\Delta \cC /c} Res_c S|_{\Delta \cC /c}
$$
where $Res_c: \cM|_{\Delta \cC /c} \to \cM(c)$ is the usual covariant restriction
to the fibre. Since $\Delta \cC /c \subset \Delta \cC$ identifies the latching categories,
if $S$ is a Reedy-cofibrant section, then $ Res_c S|_{\Delta \cC /c}: \Delta \cC /c \to \cM(c)$
can be checked to be a Reedy-cofibrant section as well (the verification is the same
as in the proof of (1.) of Lemma \ref{reedycofibrationsgivepointwise}). Thus we see
that the functor $p_!$ preserves weak equivalences between Reedy-cofibrant sections.
Since $p^*$ preserves weak equivalences, Corollary \ref{corollaryadjunctionsgiveadjunctions}
implies the derived adjunction $\bL p_! \dashv \bR p^*$ between the localisations.
It remains to check that the restriction
$$
\bL p_! : L \Sect_{  \cW}(\Delta \cC,\cM) \rightleftarrows L \Sect(\cC,  \cM): \bR p^*
$$
is an adjoint equivalence. For this, it is enough to verify that the corresponding
homotopy category adjunction
$$
\bL p_! : \Ho \Sect_{ \cW}(\Delta \cC,\cM) \rightleftarrows \Ho \Sect(\cC,  \cM): \bR p^*
$$
is an equivalence.

Since the functor $p: \Delta \cC \to \cC$ is surjective on objects,
we have that the functor $\bR p^*$ is conservative. Everything will thus follow
from the triangle identities if we verify that
the unit map $id \to \bR p^* \bL p_!$ is an isomorphism. Given a fibrewise-weakly-constant
$S \in \Ho \Sect_{  \cW}(\Delta \cC,\cM)$, the unit map evaluated
at $\tau:[n] \to \cC$ becomes
$$
S(\tau) \to \bL \colim_{\Delta \cC/c} Res_{c} S|_{\Delta \cC/c}
$$
where $c = \tau(n)$. This map is induced by the canonical inclusion $\{\tau \} \to \Delta \cC/c$. Since $\Delta \cC/c$ is a Reedy category, Theorem \ref{thmcomparison}
implies the equivalence $$L \Fun(\Delta \cC/c , \cM(c)) \cong \Fun(\Delta \cC/c, L \cM(c)),$$
and we recall that the quasicategory $L\cM(c)$ is cocomplete, Corollary \ref{cormodelcatbicomplete}.
Given the identification of left adjoints along equivalences,
it will be sufficient to prove the following. For any $X:\Delta \cC/c \to L \cM(c)$ sending
to equivalences those maps of $\Delta \cC/c$ that project to identities under the composition
$\Delta \cC/c \to \Delta \cC \to \cC$, one has that
$$
X(\sigma) \to \colim_{\Delta \cC/c} X
$$
is an equivalence. This is done exactly as in the proof of Lemma \ref{lemmacriterionlocalisation}.
\endproof

%% file: rmsf-appendices.tex
\begin{appendix}
\section{Appendix}
\subsection{Over the simplex category}

In what follows, we will identify partially ordered sets, henceforth referred as posets, with small categories having at most one morphism between each two objects.

\begin{opr}
Denote by $[n]$ the category
$$
[n ] = 0 \to 1 \to 2 \to ... \to n
$$
with exactly one morphism from $i \to j$ when $i \leq j$, and no other morphisms.
Denote by $\Delta$ the full subcategory of $\Cat$ consisting of categories $[n]$ for $n \geq 0$.
\end{opr}

\begin{lemma} Each morphism in $\Delta$ can be factored as a surjection (in the poset sense) followed by an injection (in the poset sense). Surjections and injections form a factorisation system $(\Delta_s,\Delta_i)$ on $\Delta$ which, together with the natural choice of a degree, $deg [n] = n$, makes it into a Reedy category.
\end{lemma}
\proof Clear. \endproof
\begin{corr}
	The category $\Delta^\op$ is a Reedy category for the factorisation system $(\Delta^\op_-,\Delta^\op_+)$ consisting of (the opposites of) injections and surjections.
\end{corr}

\begin{opr}
	A map $\rho:[m] \to [n]$ of $\Delta$ is \emph{a Segal inclusion}, or simply \emph{Segal} iff it is an interval inclusion of $[m]$ as first $m+1$ elements of $[n]$, i.e. $\rho(i) = i$ for $0 \leq i \leq m$. In particular, $m$ should be less or equal than $n$.

	A map $\zeta: [n] \to [m]$ of $\Delta$ is \emph{anchor} iff it preserves the endpoints: $\zeta(n) = m$.
\end{opr}

We denote by $\Alpha$, $\Sigma$ the subcategories of anchor and Segal maps in $\Delta$. It is easy to see that $(\Alpha,\Sigma)$ is a factorisation system on $\Delta$.

\begin{opr}
	A \emph{Segal factorisation system} on $\Delta^\op$ consists of the pair $(\sS,\sA)$ where $\sS$ is the subcategory of Segal maps induced from $\Sigma^\op$, and $\sA$ is the subcategory of anchor maps induced from $\Alpha^\op$.
\end{opr}

\begin{lemma}
	The identity functor sends $\Delta^\op_+$  to $\sA$. \endproof
\end{lemma}

\begin{opr}
	A $\Delta$-indexed category is a discrete Grothendieck opfibration $\pi: \cX \to \Delta^\op$ (that is, every map of $\cX$ is $\pi$-opcartesian). In particular, there exist a unique, up to isomorphism, simplicial set representing $\pi$ through the Grothendieck construction.
\end{opr}
We shall often write $\cX$ instead of $\pi$, when this abuse of notation leads to no confusion.


\begin{lemma}
	\label{deltaindexedtwosystems}
	Let $\pi:\cX \to \Delta^\op$ be a $\Delta$-indexed category. Then
	\begin{enumerate}
		\item there is a factorisation system $(\cX_-,\cX_+)$ which $\pi$- projects to $(\Delta^\op_-,\Delta^\op_+)$. We call it the Reedy factorisation system of $\cX$.
		\item There is a factorisation system $(\sS_\cX,\sA_\cX)$, which $\pi$-projects to $(\sS,\sA)$. We call it the Segal factorisation system on $\cX$.
		\item The identity functor $id: \cX \to \cX$ preserves the maps of the right class: $id(\cX_+) \subset \sA_\cX$.
	\end{enumerate}
\end{lemma}
\proof Immediate. \endproof

\begin{corr}
	Let $\cE \to \cX$ be an admissible model semifibration over
	the Reedy factorisation system $(\cX_-,\cX_+)$, then the category
	$\Sect(\cX,\cE)$ is a model category.
\endproof
\end{corr}

\subsection{Normalised model structure}
For a variation of the argument of Theorem \ref{reedymodelstructuretheorem}, let us look
at the following.
Let $\cX$ be a $\Delta$-indexed category. The subcategory $\cX_+ \subset \sA_\cX$ controls degenerations. Recasting the usual definition,

\begin{opr}
	An object $x \in \cX$ is \emph{degenerate} if there exists a non-identity degree-raising map $y \to x$ of $\cX_+$. An object $x$ is thus non-degenerate iff $\cX_+ / x = \{ id: x \to x\}$, or, equivalently, $Lat(x) = \emptyset$.
\end{opr}

For the purposes of this subsection, consider a functor $\cE \to \cX$ such that
\begin{enumerate}
	\item It is a semifibration over the Segal factorisation system.
	\item The induced semifibration over the Reedy factorisation system $(\cX_-,\cX_+)$ is an admissible model Reedy semifibration.
	\item It is normalised, that is its restriction $\cE \to \sA_\cX$ is a locally constant fibration,
	 a fibration for which all the transition functors are equivalences.
\end{enumerate}

Given a section $X \in \Sect(\cX,\cE)$ of such $\cE \to \cX$, we can thus conclude that $\Lat_x X$ is the initial object of $\cE(x)$ for each non-degenerate $x$.


\begin{opr}
	A section $X$ is \emph{normalised} iff it takes any arrow of $\cX_+$ to an opcartesian arrow of $\cE \to \cX$.
\end{opr}

\begin{lemma}
	A section $X$ is \emph{normalised} iff for any degenerate object $y \in \cX$, the latching map $\Lat_y X \to X(y)$ is an isomorphism.
\end{lemma}
\proof Given the definition of a normalised section, we have that for each $f:x \to y$ in $\cX_+ / y$, the map $f_! X(x) \to X(y)$ is an isomorphism. One then checks that the latching category $Lat(y) \subset \cX_+ / y$ is connected
and so the colimit of a constant $Lat(y)$-diagram with value $X(y)$ gives $X(y)$.  \endproof

\begin{rem}
	If we take $x \to y$ to be an ordinary degeneracy (if projected to $\Delta$) then $X(x) \to X(y)$ is an isomorphism (note that $\cE(x) \cong \cE(y)$.
\end{rem}

Denote by $\Sect_{N}(\cX,\cE) \subset \Sect(\cX,\cE)$ the full subcategory of normalised sections.

\begin{lemma}
	\label{nondegeneratelimits}
	The category $\Sect_N(\cX,\cE)$ admits limits and colimits, which are calculated in $\Sect(\cX,\cE)$.
\end{lemma}
\proof The colimit part is trivial and is left to the reader. For the limit part, we will use the Segal factorisation system on $\cX$ to calculate limits. For the proof, recall also the functor $\pi: \cX \to \Delta^\op$.

Let $x \in \cX$, and consider the category $x \backslash \sS_\cX$. Given that on the level of $\Delta$, the maps of $\sS_\cX$ are interval inclusions, and so we have an equivalence $x \backslash \sS_\cX \cong \pi(x) \in \Delta \subset \Cat$.

Now, consider a morphism $f: x \to x'$ in $\cX_{+}$. It also means that $f$ belongs to $\sA_\cX$, but in any case, the factorisation system $(\sS_\cX,\sA_\cX)$ defines a functor
$
\bar f: x' \backslash \sS_\cX \to x \backslash \sS_\cX
$
by projecting to $\Delta$, one can examine and check that $f^*$, after the equivalences $x' \backslash \sS_\cX \cong \pi x'$ and $x \backslash \sS_\cX \cong \pi x$, is just the map
$$
\pi(f): \pi x' \to \pi x
$$
corresponding to $f$ by projection to $\Delta^\op$.
In all, we constructed the following diagram
\begin{diagram}[small]
x' \backslash \sS_\cX & \rTo^{\bar f} & x \backslash \sS_\cX \\
\uTo>\cong					&						&		\uTo>\cong \\
\pi x'						&	\rTo^{\pi(f)}	&		\pi x.
\end{diagram}
If we note by $p_x: x \backslash \sS_\cX \to \cX$ and $p_{x'} :x' \backslash \sS_\cX$ the natural projections, then the map $\bar f^* p_{x}^* \cE \to p_{x'}^* \cE$ (cf Proposition \ref{janus2functoriality}) of prefibrations over  $x' \backslash \sS_\cX $ is in fact an equivalence due to the normalisation condition, since the natural transformation which induces it, $p_x \bar f \to p_{x'}$, lies in $\cX_+$ and not just in $\sA_\cX$. Hence there is no confusion about lifting $\cE \to \cX$ to this diagram.
When computing limits in $\Sect(\cX,\cE)$, it is done by taking limits of certain sections over categories like $x \backslash \sS_\cX$. It will thus suffice to check that the functor
$$
\Sect(x \backslash \sS_\cX , p_x^* \cE) \stackrel{\bar f ^*}{\longrightarrow} \Sect(x' \backslash \sS_\cX, \bar f ^* p_x^* \cE) \cong \Sect(x' \backslash \sS_\cX, p^*_{x'} \cE)
$$
preserves limits, and the resulting section will then be normalised. But this is equivalent to showing that the functor
$$
\pi (f) ^* : \Sect (\pi x, \cE) \to \Sect( \pi x , \cE)
$$
preserves limits. This is sufficient to test when $\pi f$ is an elementary degeneracy, and in this case $\pi f: \pi x' \to \pi x$ admits both left and right adjoints. All this suffices to show that, when we compute a limit of a diagram of normalised sections, the values of the limit on degeneracies are isomorphisms.
\endproof

Denote by $\cX_{nd}$ the subcategory of $\cX$ consisting of nondegenerate objects. Clearly, $\cX_{nd} \subset \cX_-$, and moreover it is naturally an inverse Reedy category. Consequently, for each $x \in \cX_{nd}$ and a section $X: \cC \to \cE$, we can define $\Mat^{nd}_x X$, the matching object of $X$ at $x$ in the category $\Sect(\cX_{nd},\cE)$. It is defined as the limit
$$
\Mat^{nd}_x X = \lim_{Mat^{nd}(x)} R_x X|_{Mat^{nd}(x)}
$$
where $Mat^{nd}(x) \subset x \backslash \cX_{nd}$ is the subcategory of all maps out of $x$ in $\cX_{nd}$ safe the identity.

The inclusion $x \backslash \cX_{nd} \subset x \backslash \cX_-$ induces the functor $Mat^{nd}(x) \subset Mat(x)$, and thus a map $\Mat_x X \to \Mat^{nd}_x X$.

\begin{lemma}
	\label{nondegeneratetodegeneratematching}
	Let $X$ be a normalised section. Then the map $\Mat_x X \to \Mat^{nd}_x X$ is an isomorphism for each $x \in \cX_{nd}$.
\end{lemma}
\proof One has to observe, that in $x \backslash \cX_-$, there are objects
$x \to y$ such that $y$ may be degenerate. For such each $y$, choose a non-degenerate $\bar y$ and a map $\bar y \to y$ in $\cX_+$ degenerating $y$. Each such map admits a section $y \to \bar y$, which lies in $\cX_-$.

Moreover, if $y \to z$ is a map in $\cX_-$ to a non-degenerate object, there exists a factorisation $y \to \bar y \to z$ in $\cX_-$, where $\bar y$ is non-degenerate as before. One can see that such observations are sufficient to prove that the functor $Mat^{nd}(x) \to Mat(x)$ is final (or right cofinal in the sense of \cite{HIRSCHHORN}), and this implies the isomorphism of limits.
\endproof

\begin{thm}
	\label{normalisedreedymodelstructure}
	The category $\Sect_{N}(\cX,\cE) $ possesses a model structure with limits and colimits created by the inclusion to $\Sect(\cX,\cE)$. The classes of cofibrations, fibrations and weak equivalences are given as follows:
	\begin{itemize}
		\item a map $A \to B$ of $\Sect_{N}(\cX,\cE) $ is a cofibration iff it is a Reedy cofibration in $\Sect(\cX,\cE)$,
		\item a map $A \to B$ of $\Sect_{N}(\cX,\cE) $ is a weak equivalence iff it is such in $\Sect(\cX,\cE)$,
		\item a map $X \to Y$ of $\Sect_{N}(\cX,\cE) $ is a fibration iff for each \emph{non-degenerate} object $x \in \cX$, the relative matching map $X(x) \to Y(x) \prod_{\Mat_x Y } \Mat_x X$ is a fibration in $\cE(x).$
	\end{itemize}
	Moreover $\Mat_x X \cong \Mat^{nd}_x X$ for each nondegenerate $x \in \cX_{nd}$.
\end{thm}

\begin{lemma}
	\label{normalisedtrivialcoffib}
	In $\Sect_{N}(\cX,\cE)$,
	\begin{itemize}
		\item a map $A \to B$ is a cofibration and a weak equivalence iff for each $x \in \cX$, the relative latching map $\Lat_x B  \coprod_{\Lat_x A } A(x) \to B(x)$ is a trivial cofibration in $\cE(x)$,
		\item  a map $X \to Y$ is a fibration and a weak equivalence iff for each non-degenerate object $x \in \cX$, the relative matching map $X(x) \to Y(x) \prod_{\Mat_x Y } \Mat_x X$ is a trivial fibration in $\cE(x).$
	\end{itemize}
\end{lemma}
\proof The first is done by restricting to $\Sect(\cX_+,\cE)$ and using Corollary \ref{trivialreedycof}, just as for Lemma \ref{lemmatrivreedycofifffibrewise}. The second is done by restricting to $\Sect(\cX_{nd},\cE)$, and using the dual of Corollary \ref{trivialreedycof} together with Lemma \ref{nondegeneratetodegeneratematching}. \endproof

\proof[Proof of Theorem \ref{normalisedreedymodelstructure}.]
\begin{enumerate}
	\item The limits and colimits axiom is clear, see Lemma \ref{nondegeneratelimits}.
	\item The weak equivalences of $\Sect(\cX,\cE)$ satisfy three-for-two, hence the same property applies for the weak equivalences between non-degenerate sections.
	\item The retract stability for the three classes of maps is verified just as in Lemma \ref{retracts}.
	\item The lifting is proven analogously to the Reedy case. Consider a diagram
	\begin{diagram}[small,nohug]
	A	& \rTo	& S	 \\
	\dTo<f	 &			&	\dTo>p		\\
	B & \rTo & T	 \\
	\end{diagram}
	with, say, $f$ a cofibration and $p$ a trivial fibration, and we keep in mind the result of Lemma \ref{normalisedtrivialcoffib}. We observe that each degree zero object $x$ of $\cX$ has empty latching and matching categories, and is moreover non-degenerate. Hence in this case the relative latching map reduces to a cofibration $A(x) \to B(x)$, the relative matching map reduces to a trivial fibration $S(x) \to T(x)$, and finding a lifting is trivial. For the induction step, consider, for a non-degenerate $x \in \cX_{nd}$, the diagram
	\begin{diagram}[size=2.5em,nohug]
	A(x) & \rTo 			& A(x) \coprod_{\Lat_x A} \Lat_x B	& \rTo	& S(x)	&		&	 \\
	&	\rdTo		&	\dTo	 &			&	\dTo											&	\rdTo	&\\
	&					&	B(x) & \rTo & T(x) \prod_{\Mat_x T} \Mat_x S	& \rTo	& T(x). \\
	\end{diagram}
	which admits a lifting as in Reedy case. If $y \in \cX$ is, however, a degenerate object, then $\Lat_y A \cong A(x)$ and $\Lat_y B \cong B(x)$, and the square
	\begin{diagram}[small,nohug]
	\Lat_y A	& \rTo^\sim	& A(y)	 \\
	\dTo	 &			&	\dTo>{f(y)}		\\
	\Lat_y B & \rTo^\sim & B(y)	 \\
	\end{diagram}
	is a pushout, hence the relative latching map is isomorphic to $B(y) \to B(y)$, and finding the lift in
	\begin{diagram}[size=2.5em,nohug]
	A(y) & \rTo 			& B(y) & \rTo	& S(y)	&		&	 \\
	&	\rdTo		&	\dTo>=	 &			&	\dTo											&	\rdTo	&\\
	&					&	B(y) & \rTo & T(y) \prod_{\Mat_y T} \Mat_y S	& \rTo	& T(y). \\
	\end{diagram}
	is trivial, whichever the property the map on the right of the square possesses.
	\item Assume given a map of normalised section $A \to C$. Degree zero objects $x$ are non-degenerate and have no matching-latching categories, so we simply factor our map as $A(x) \to B(x) \to C(x)$ using the model structure of $\cE(x)$. So far, $B$ is trivially a normalised section.

	By induction, we have constructed the factorisation $A(y) \to B(y) \to C(y)$ for objects $y \in \cX$ of degree less than $n$, and $B: \cX_{<n} \to \cE$ is non-degenerate. For $x$ of degree $n$, there is the following diagram
	\begin{diagram}[small]
	\Lat_x A & \rTo & A(x)		&	\rTo & \Mat_x B  \\
	\dTo		&			&	\dTo		&			&		\dTo			\\
	\Lat_x B &	 \rTo & C(x)		&	\rTo	&	\Mat_x C	\\
	\end{diagram}

	which exists due to the inductive assumption and provides us with the following map
	$$
	\Lat_x B \coprod_{\Lat_x A} A(x) \to C(x) \prod_{\Mat_x C} \Mat_x B.
	$$
	If $x$ is non-degenerate, we factor it as
	$$
	\Lat_x B \coprod_{\Lat_x A} A(x) \to B(x) \to C(x) \prod_{\Mat_x C} \Mat_x B.
	$$
	which, together with maps $\Lat_x B \to B(x)$ and $B(x) \to \Mat_x B$, yields the desired extension of the factorisation to $x$. For a degenerate object $y$, we simply put $B(y) = \Lat_y B \coprod_{\Lat_y A} A(y) $. Then the natural map $\Lat_y B \to B(y)$ is an isomorphism (since $A$ is normalised) and the factorisation
	$$
	\Lat_y B \coprod_{\Lat_y A} A(y) = B(y) \to C(y) \prod_{\Mat_y C} \Mat_y B.
	$$
	is as needed, given the first map satisfies lifting along any map of $\cE(y)$ and the second map is not forced to any condition. \endproof
\end{enumerate}
\end{appendix}

%% file: rmsf-master.bbl
\begin{thebibliography}{99}

	\bibitem{ADAMEK} J. Adamek, J. Rosicky, \textit{Locally Presentable and Accessible Categories},
	Cambridge University Press 1994





	\bibitem{HARPAZ2} Ilan Barnea, Yonatan Harpaz, Geoffroy Horel, \textit{Pro-categories in homotopy theory},  Algebraic and Geometric Topology, 17 (1), 2017, p. 567-643

	\bibitem{BARLEFT} Clark Barwick, \textit{On (Enriched) Left Bousfield Localization of Model Categories}, preprint \url{https://arxiv.org/abs/0708.2067}



	\bibitem{BARKAN0} Clark Barwick, Daniel Kan, \textit{Relative categories: Another model for the homotopy theory of homotopy theories}, preprint \url{https://arxiv.org/abs/1011.1691}

	\bibitem{BARKAN} Clark Barwick, Daniel Kan, \textit{In the category of relative categories the Rezk equivalences are exactly the DK-equivalences}, preprint \url{https://arxiv.org/abs/1012.1541}








	\bibitem{CISBOOK}  Denis-Charles Cisinski, \textit{Higher Categories and Homotopical Algebra,} on-line book, available at
	\url{http://www.mathematik.uni-regensburg.de/cisinski/CatLR.pdf}


	\bibitem{DUGGER} Daniel Dugger, \textit{Combinatorial model categories have presentations}, Advances in Mathematics
Volume 164, Issue 1, 1 December 2001, Pages 177-201 \url{https://arxiv.org/abs/math/0007068}

	\bibitem{DHKS} William G. Dwyer, Philip S. Hirschhorn, Daniel M. Kan, and Jeffrey H. Smith, \textit{Homotopy Limit Functors on Model Categories and Homotopical Categories}, AMS 2004

	\bibitem{GEPNERNIK} David Gepner, Rune Haugseng, Thomas Nikolaus \textit{Lax colimits and free fibrations in $\infty$-categories}, \url{https://arxiv.org/abs/1501.02161}


	\bibitem{GROTH} Moritz Groth, \textit{On the theory of derivators}, doctoral dissertation, Bonn 2011, \url{http://www.math.uni-bonn.de/people/grk1150/DISS/dissertation-groth.pdf}


	\bibitem{SGA1} Alexander Grothendieck, Mich\`ele Raynaud et al., \textit{Rev\^etements \'etales et groupe fondamental (SGA I)}, Lecture Notes in Mathematics 224, Springer 1971

	\bibitem{HINICH} Vladimir Hinich, \textit{Dwyer–Kan localization revisited},
	Homology, Homotopy and Applications, Volume 18 (2016), Number 1

	\bibitem{HIRSCHHORN} Philip S. Hirschhorn, \emph{Model Categories and Their Localisations},  No. 99. American Mathematical Soc., 2009.

	\bibitem{HARPAZCOMP} Yonatan Harpaz, \textit{Lax limits of model categories}, preprint

	\bibitem{HARPAZ} Yonatan Harpaz, Joost Nuiten, Matan Prasma, \textit{The abstract cotangent complex and Quillen cohomology of enriched categories}, preprint  \url{https://arxiv.org/abs/1612.02608}

\bibitem{HAPRASMA} Yonatan Harpaz, Matan Prasma, \textit{The Grothendieck construction for model categories},
preprint \url{https://arxiv.org/abs/1404.1852}


	\bibitem{H-S} Andr\'e Hirschowitz, Carlos Simpson, \textit{Descente pour les n-champs (Descent for n-stacks)}, preprint \url{http://arxiv.org/abs/math/9807049}


	\bibitem{HOVEY} Mark Hovey, \textit{Model Categories}, Mathematical Surveys and Monographs, vol. 63, American Mathematical Society, Providence, Rhode Island, 1999.












	\bibitem{LUHTT} Jacob Lurie, \textit{Higher Topos Theory}, book, Annals of mathematics studies ; no. 170

	\bibitem{LU} Jacob Lurie, \textit{Higher Algebra}, on-line book, available at \url{http://www.math.harvard.edu/~lurie/}

	\bibitem{MGEE1} Aaron Mazel-Gee, \textit{Goerss–Hopkins obstruction theory via model $\infty$-categories},
	PhD Thesis, 2016, \url{etale.site}

	\bibitem{MGEE2} Aaron Mazel-Gee, \textit{A user's guide to co/cartesian fibrations}, preprint
	\url{https://arxiv.org/abs/1510.02402}


	\bibitem{ML} Saunders Mac Lane, \textit{Categories for The Working Mathematician,} Graduate Texts in Mathematics 5 (second ed.). Springer, 1998

	\bibitem{MANDELL} Michael A. Mandell, \textit{Equivalence of simplicial localizations of closed model categories},
	Journal of Pure and Applied Algebra 142 (1999) 131–152


	\bibitem{QHA} Daniel Quillen, \textit{Homotopical Algebra}, Lecture Notes in Mathematics book series (LNM, volume 43)


	\bibitem{REEDY} Chris L. Reedy, \textit{Homotopy theory of model categories}, preprint 1974

	\bibitem{RIEHL} Emily Riehl, \textit{Categorical Homotopy Theory,} Cambridge University Press, 2014


	\bibitem{SHAH} Jay Shah, \textit{Parametrized higher category theory}, PhD thesis \url{http://math.mit.edu/~jshah/thesis.pdf}

	\bibitem{SPITZCOMP} Markus Spitzweck, \textit{Homotopy limits of model categories over inverse index categories}, Journal of Pure and Applied Algebra Volume 214, Issue 6, June 2010, Pages 769-777

	\bibitem{STEVENSON} Danny Stevenson, \textit{Covariant Model Structures and Simplicial Localization}, North-West. Eur. J. Math. Vol. 3 (2017) pp. 141-202, \url{https://arxiv.org/abs/1512.04815}



	\bibitem{HAG2} Bertrand To\"en, Gabriele Vezzosi, \textit{Homotopical Algebraic Geometry II: geometric stacks and applications}, Memoirs of the American Mathematical Society (Book 92), April 27, 2008

	\bibitem{TV} Bertrand To\"en, Gabriele Vezzosi, \textit{Caract\`eres de Chern, traces \'equivariantes et g\'eom\'etrie alg\'ebrique d\'eriv\'ee} G. Sel. Math. New Ser. (2015) 21: 449, \url{https://arxiv.org/abs/0903.3292}


	\bibitem{VIST} Angelo Vistoli, \textit{Notes on Grothendieck topologies, fibered categories and descent theory,} \url{http://arxiv.org/abs/math/0412512}


\end{thebibliography}
